\documentclass[times,sort&compress,3p]{elsarticle}
\usepackage{graphicx,psfrag,epsf}
\usepackage{enumerate}
\usepackage{natbib}
\usepackage{url} % not crucial - just used below for the URL 
\usepackage{amssymb,amsmath,amsthm,tikz,ifthen,subfigure,booktabs,slashbox}
\usepackage{titlesec,setspace,mathptmx}
\usepackage[ngerman,english]{babel}
\usepackage[babel,german=quotes]{csquotes}
\usepackage{wasysym}
\usepackage[weather]{ifsym}
\usepackage{multirow}
\usepackage[ansinew]{inputenc} % to include Umlaute
\usepackage{multibib}
\usepackage{scalerel}
\usepackage{enumitem}   
\usepackage{listings}
\usepackage{subfigure}

\usepackage[mathscr]{euscript}
\usetikzlibrary{plotmarks,calc,arrows}
\usepackage{rotating}
\usepackage{units}
\usepackage{amsmath,amsthm}

\usepackage{hyperref}
\hypersetup{
	pdfstartview={FitH},
	colorlinks=true,
	linkcolor=black,
	urlcolor=black,
	citecolor=black,
	bookmarksopen,
	bookmarksopenlevel={1},
	bookmarksnumbered
}

\newtheorem{satz}{Theorem} 
\newtheorem{prop} {Proposition}
\newtheorem{lemma} {Lemma}
\newdefinition{defi} {Definition}
\theoremstyle{definition}

\newtheorem{bem}{Remark}

\newtheorem{assumpx}{Assumption} 
\newenvironment{assump}
{\pushQED{\qed}\assumpx}
{\popQED\endassumpx}

\theoremstyle{plain}

\usepackage{marginnote}

\newcommand\R{\mathbb R} % K?rper der reellen Zahlen
\newcommand\N{\mathbb N} % nat?rliche Zahlen	
\newcommand\E{{\rm E}} % Erwartungswert
\newcommand{\argmax}[1]{\underset{#1}{\operatorname{arg}\,\operatorname{max}}\;} % argmax
\newcommand\Hr{\mathbb{H}} %  H
\newcommand\M{\mathbb{M}}
\newcommand\veps{\varepsilon}

\newcommand\abs[1]{ |\,#1 \,|} % Absolutbetrag
\newcommand\norm[1]{\| #1 \|}  
\newcommand\scalarp[1]{\langle\,#1\,\rangle} % 2-Norm
\newcommand\norminf[1]{||\,#1\,||_\infty} % Infinity Norm
\newcommand\inv{^{-1}} % inverse
\newcommand\lb{\Big(}
\newcommand\rb{\Big)}

\newcommand\Pto{\stackrel{\text{P}}{\to}}

\newcommand\dx{\mathrm{d} x }
\newcommand\dy{\mathrm{d} y }
\newcommand\dz{\mathrm{d} z }
\newcommand\dt{\mathrm{d} t }
\newcommand\epi{\text{epi}}

\newcommand\sumieiz{ \sum\limits_{i_1,i_2=1}^n }
\newcommand\ionetwo{i_1,i_2}
\newcommand\yionetwo{Y_{\ionetwo}}
\newcommand\xionetwo{\pvec{x_{\ionetwo}}}
\newcommand\eionetwo{  \boldsymbol{{\veps}}_{\ionetwo}}
\newcommand\pvec[1]{\mathbf{#1}}

\renewcommand\dz{\mathrm{d} \pvec{z} }

\newcommand\dzone{\mathrm{d} z_1 }
\newcommand\dztwo{\mathrm{d} z_2 }
\newcommand\du{\mathrm{d} \pvec{u} }

\newcommand\psix{\psi(x)}
\newcommand\empcrit{\hat\M_n}
\newcommand\detcrit{\M_n}
\newcommand\asympcrit{ \M}
\newcommand\estpsi{\hat{\psi}_n}
\newcommand\estphi{\hat{\phi}_n}
\newcommand\estw{\hat{w}_n}
\newcommand\esttau{\hat{\tau}_n}
\newcommand{\diag  }{{\rm diag}}

\newcommand{\px}{\pvec{p}(x)}

\newcommand{\rotmat}{ \pmb{\mathscr R}}

\newcommand{\paramset}{\Theta_n}
\newcommand{\maxset}{\tilde \Theta_{n,x}}
\newcommand{\xset}{I}
\newcommand{\levyanti}{\mathbb{L}}

\newcommand{\zbootphi}{Z^\xi_{n, \widehat \phi}}
\newcommand{\zbootpsi}{Z^\xi_{n, \widehat \psi}}
\newcommand{\zboottau}{Z^\xi_{n, \widehat \tau}}

\newcommand{\zscorephi}{Z^\varepsilon_{n, \phi}}
\newcommand{\zscorepsi}{Z^\varepsilon_{n, \psi}}
\newcommand{\zscoretau}{Z^\varepsilon_{n, \tau}}

\newcommand{\bootprocessphi}{\mathscr{Z}^\xi_{n,\widehat\phi}}
\newcommand{\bootprocesspsi}{\mathscr{Z}^\xi_{n,\widehat \psi}}
\newcommand{\bootprocesstau}{\mathscr{Z}^\xi_{n,\widehat \tau}}

\newcommand{\orgprocessphi}{{\mathscr{Z}}_{n,\phi}^\varepsilon}
\newcommand{\orgprocesspsi}{{\mathscr{Z}}_{n,\psi}^\varepsilon}
\newcommand{\orgprocesstau}{{\mathscr{Z}}_{n,\tau}^\varepsilon}

\newcommand{\remainderpartphipsi}{R_n^{\phi,\psi}}

\newcommand{\remainderparttau}{R_n^{\tau}}

\newcommand{\scorepartphipsi}{\pmb{\mathscr S}_{n}^{\phi,\psi}}
\newcommand{\scorepartphi}{\pmb{\mathscr S}_{n}^{\phi}}
\newcommand{\scorepartpsi}{\pmb{\mathscr S}_{n}^{\psi}}
\newcommand{\scoreparttau}{\pmb{\mathscr S}_{n}^{\tau}}

\newcommand{\biaspartphipsi}{\pmb{\mathscr B}_n^{\phi,\psi}}
\newcommand{\biaspartphi}{\pmb{\mathscr B}_n^{\phi}}
\newcommand{\biaspartpsi}{\pmb{\mathscr B}_n^{\psi}}
\newcommand{\biasparttau}{\pmb{\mathscr B}_n^{\tau}}

\newcommand{\varphisco}{V_{\phi}^S}
\newcommand{\varpsisco}{V_{\psi}^S}
\newcommand{\vartausco}{V_{\tau}^S}
\newcommand{\varphihess}{V_{\phi}^H}
\newcommand{\varpsihess}{V_{\psi}^H}

\newcommand{\varphiscoest}{V_{\widehat\phi}^S}

\newcommand{\varphihessest}{V_{\widehat\phi}^H}
\newcommand{\varpsihessest}{V_{\widehat\psi}^H}

\newcommand{\IGNORE}[1]{}

\usepackage{lipsum}

\newcommand\blfootnote[1]{%
	\begingroup
	\renewcommand\thefootnote{}\footnote{#1}%
	\addtocounter{footnote}{-1}%
	\endgroup
}

\makeatletter
\def\ps@pprintTitle{%
	\let\@oddhead\@empty
	\let\@evenhead\@empty
	\def\@oddfoot{\centerline{\thepage}}%
	\let\@evenfoot\@oddfoot}
\makeatother

\usepackage{xr}
\usepackage{pdfpages}

\begin{document}
	
	\allowdisplaybreaks

	\begin{frontmatter}
		\title{Asymptotic confidence sets for the jump curve \\
		in bivariate regression problems}
		\author{Viktor Bengs, Matthias Eulert}
		\author{Hajo Holzmann\corref{mycorrespondingauthor}}
		\address{Fb. 12 - Mathematik und Informatik, Philipps-Universit\"at Marburg, Hans-Meerwein-Straße 6, 35032 Marburg, Germany}
		
		\cortext[mycorrespondingauthor]{Corresponding author}
		\ead{holzmann@mathematik.uni-marburg.de}
		
\begin{abstract}
We construct uniform and {point-wise} asymptotic confidence sets for the single edge in an otherwise smooth image function which are based on  rotated differences of two one-sided kernel estimators. Using methods from M-estimation, we show consistency of the estimators of location, slope and height of the edge function and develop a uniform linearization of the contrast process. The uniform confidence bands then rely on a Gaussian approximation of the score process together with anti-concentration results for suprema of Gaussian processes, while {point-wise} bands are based on asymptotic normality. The finite-sample performance of the {point-wise} proposed methods is investigated in a simulation study. An illustration to real-world image processing is also given.  
\end{abstract}
		
\noindent
		
\vfill
\end{frontmatter}
		
\section{Introduction}

\blfootnote{\copyright 2019. This manuscript version is made available under the CC-BY-NC-ND 4.0 license \href{http://creativecommons.org/licenses/by-nc-nd/4.0/}{http://creativecommons.org/licenses/by-nc-nd/4.0/}. }

Statistical methodology, in particular nonparametric regression, plays an important role nowadays in image reconstruction and denoising. In two-dimensional image functions, interest often focuses on the location of edges, i.e., discontinuity curves of the intensity function of the image. Edge estimation has been studied in the statistics literature from a minimax point-of-view in the monograph by \citet{korostelev1993minimax}. \citet{qiu2005image} gives an overview of more practical reconstruction methods, which are partially motivated from the literature on computer vision. 
	
Apart from mere estimation, statistical modeling allows for the construction of confidence regions, in which the object of interest is located with high probability. In this paper, we construct asymptotic confidence sets for the single edge in an otherwise smooth image function based on the rotational difference kernel method by \citet{qiu1997nonparametric} and \citet{muller1994maximin}, which is related to the Sobel edge detector from the image processing literature \citep{qiu2002nonparametric}. 
	
Confidence sets are by now well-developed for various problems in nonparametric statistics, e.g., for nonparametric density estimation \citep{bickel1973some, gine2010confidence, chernozhukov2014anti}, smooth regression functions \citep{eubank1993confidence, neumann1998simultaneous}, and deconvolution and errors-in-variables problems \citep{bissantz2007, proksch2015confidence, delaigle2015confidence}. \citet{mammenpolonik2013} and \citet{qiaopolonik2016} focus on more geometrical features, and construct confidence regions for density level sets and the density ridge, respectively. 
		
When studying nonparametric regression problems with discontinuities, the focus may either be on the detection of potential jumps, or, assuming the existence of jumps, on estimation and the construction of confidence sets. In the univariate setting with a regression curve having a single jump, the latter problem was studied under various design assumptions by \citet{muller1992change}, \citet{loader1996change}, \citet{gijbels2004interval} and \citet{seijo2011change}, among others. The problem of jump detection was additionally addressed, e.g., in \cite{muller1999discontinuous, wu1993kernel, wang1995jump, porter2015regression}. 
	
In bivariate problems, \citet{korostelev1993minimax}, \citet{muller1994maximin}, \citet{qiu1997jump}, \citet{qiu2002nonparametric}, and   \citet{garlipp2007robust} focus on the estimation of a given discontinuity curve, while \citet{wang1998change}, \citet{qiu1997jump}, \citet{kang2014jump}  and \citet{qiu2002nonparametric} also investigate edge detection. In the latter three papers, the set of jump-location curves is estimated as a point-set, where a certain criterion function exceeds a threshold value.  
	
However, currently there seem to be no methods available to construct a confidence set for the edge in the bivariate case. Our approach to this problem uses the criterion function of the rotational difference kernel method \citep{qiu1997nonparametric, muller1994maximin}, for which we obtain necessary additional kernel conditions for asymptotic normality. The uniform confidence bands then rely on a Gaussian approximation of the score process together with anti-concentration results for suprema of Gaussian processes from \citet{chernozhukov2014anti}, while point-wise bands are based on asymptotic normality. A technical difficulty in the problem are the distinct rates of the estimators of location and slope. As a byproduct of our analysis, we obtain uniform rates of convergence for the estimators of the jump-location-curve, the slope-curve as well as the height-curve which are optimal up to a logarithmic factor.

The paper is structured as follows. In Section~\ref{sec:model_notations} we introduce the model as well as the estimators for location and slope of the edge. The main theoretical results can be found in Section~\ref{sec:main_results}. Outlines of the proofs for these results are provided in Section~\ref{sec:proofs}. Section~\ref{sec:simulation} contains a simulation study and an illustrative application of the proposed method to a real-world image, while Section~\ref{sec:discuss} concludes. Detailed proofs are provided in the Online Supplement \citep{beh2018supplement}.  
	
We shall use the following notation. If a sequence of random variables $X_1,X_2,\ldots$ converges in probability to a random variable $X$, we write $X_n \Pto X$ or $X_n = X +o_P(1)$. If the convergence is in distribution, we write $X_n \rightsquigarrow X$. Let $\lambda_d$ denote the $d$-dimensional Lebesgue measure. For $d\in \N,$ $\mu\in \R^d, \Sigma \in \R^{d\times d}$ symmetric and positive definite, let  $\mathcal{N}_d(\mu,\Sigma)$ be the normal distribution with expectation $\mu$ and covariance matrix $\Sigma.$ For a vector $\pvec{a} \in \R^d$, we denote by $\diag  (\pvec{a})$ the $d \times d$ diagonal matrix with diagonal entries $\pvec{a}$. Given $\alpha\in(0,1)$, we denote by $q_\alpha(X)$ the $\alpha$-quantile of a random variable $X$ resp.\  by  $q_\alpha(Q)$ the $\alpha$-quantile of a distribution $Q$. Let $\scalarp{\cdot, \cdot }$ denote the Euclidean scalar product on $\R^d$.
 
 \section{Model and notation} \label{sec:model_notations}
 
We consider a boundary fragment of a noisy gray scale image, in which real-valued random variables $Y_{i_1,i_2}$ are observed according to the model defined, for all $(i_1,i_2) \in \{1,\ldots,n\} \times \{1,\ldots,n\}$, by
\begin{align} \label{model}
\yionetwo = m_\phi(\xionetwo) + \eionetwo,
\end{align}
where $\xionetwo=(x_{i_1},x_{i_2})^\top$ form a deterministic, regular rectangular grid in $[0,1]^2,$ and $m_\phi$ is an unknown square-integrable function on $[0,1]^2,$ which is sufficiently smooth besides a discontinuity curve $\phi$. Specifically, we assume that  $m_\phi$ is of the form
\begin{align}\label{assumption_regression_function}
\begin{split}
m_\phi(x,y) & = m(x,y) + j_{\tau,\phi}(x,y),\\
j_{\tau,\phi}(x,y) & = j_\tau(x,y) = \tau(x) \mathbf{1}_{[0,\phi(x)]}(y),
\end{split}
\end{align}
where  $m:[0,1]^2 \to \R$ is the smooth part of the image,  $\tau:[0,1] \to \R_+$ the jump-height curve and $\phi:[0,1] \to (0,1)$  the jump-location curve.

To estimate $\phi(x)$, define the rotation matrix
\begin{align}  \label{eq:rotation_matrix}
\rotmat_{\psi}  = \left(  \begin{array}{cc}
\cos(\psi) & \sin(-\psi) \\ \sin(\psi) & \cos(\psi)
\end{array}	 	\right),
\end{align}
where $\psi \in \R$. Further, for a bandwidth $h>0$ and $\pvec{z}=(z_1,z_2)^\top \in \R^2$, consider the rotated difference kernel
\begin{align*}
K(\pvec{z};\psi,h) =  K(h\inv  \rotmat_{-\psi} \pvec{z})/h^2,
\end{align*}
where $ K(\pvec{z}) = K(z_1,z_2)=K_1(z_1)K_2(z_2)$ is a product kernel of univariate kernel functions $K_1$ and $K_2$, and $K_1$ is even while $K_2$ is odd and hence corresponds to the difference of two one-sided kernels.  	

Define for $\pvec{z} \in[0,1]^2$ and $\psi \in[-\pi/2,\pi/2]$ the contrast process with Priestley--Chao-type weights as
\begin{align} \label{defi:contrast_emp}
\hat{\text{M}}_n(\pvec{z};\psi,h) = n^{-2} \sum\limits_{i_1,i_2=1}^n \yionetwo 	K(\pvec{z}-\xionetwo;\psi,h).
\end{align} 
For an $x\in(0,1)$ denote by $\psix= \arctan \{ \phi'(x)\}$ the slope of the tangent at $\phi(x)$. An estimator for the bivariate parameter $(\phi(x), \psix )$ is then given by
\begin{align} \label{definition_estimator_phi}
(\estphi(x), \estpsi(x) )
\in  \argmax{   y \in [h,1-h], \psi \in[-\pi/2,\pi/2] }  \hat{\text{M}}_n((x,y)^\top;\psi,h).	
\end{align}
Our estimator for the jump-height at $x$ is given by
$$
\esttau(x) =    \hat{\text{M}}_n ((x,\estphi(x))^\top;\estpsi(x),h)  .
$$
	
\begin{figure} [t!]
% FIGURE 1
	\centering
	\includegraphics[width=0.9\linewidth]{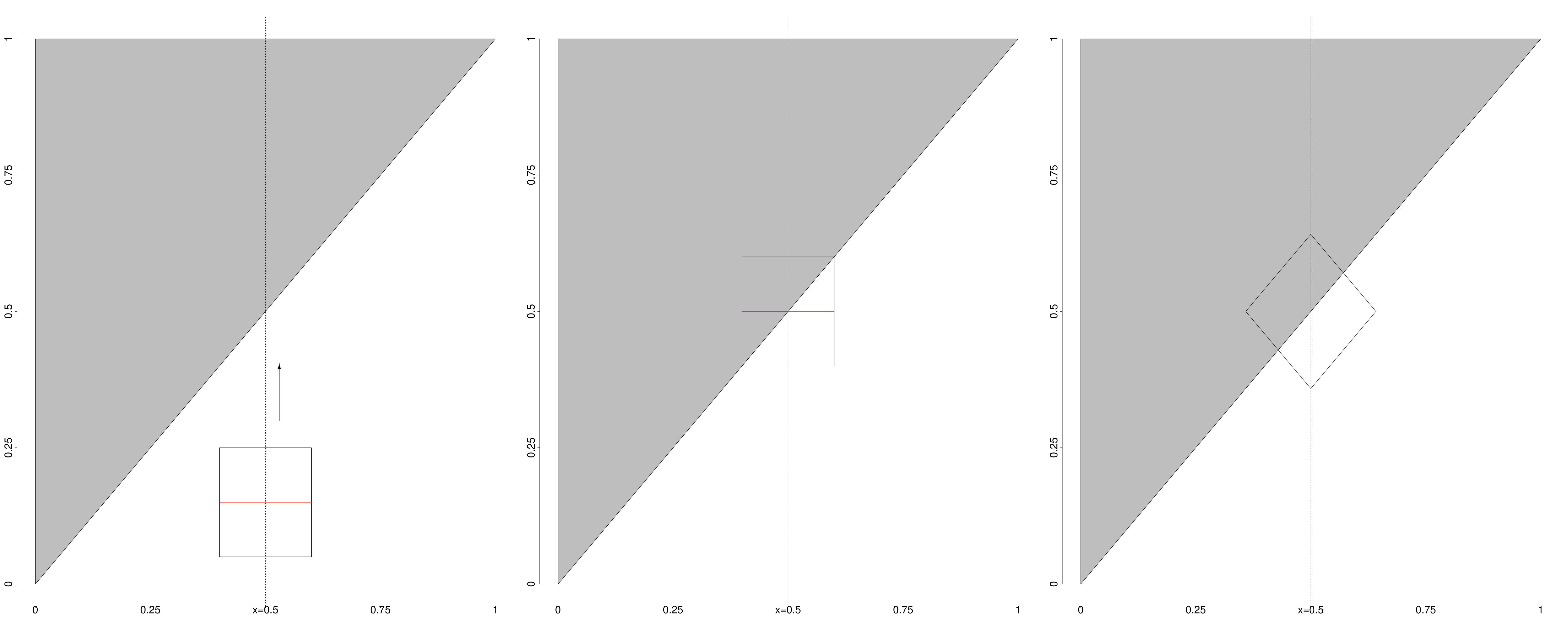}
	\caption{Illustration of rotated difference kernel estimation for a linear jump-location-curve.}
	\label{fig:kernelillustration}
\end{figure}

An illustration of the  mechanism in (\ref{definition_estimator_phi}) is given in Figure~\ref{fig:kernelillustration}. The line within the kernel window indicates that the odd kernel $K_2$ corresponds to the difference of two one-sided kernels. In order to maximize the contrast $\hat{\text{M}}_n(\pvec{z};\psi,h)$, that is the difference of two one-sided estimators, for fixed $x\in(0,1)$ the kernel window is rotated so that the red line is tangential to the edge. 

\begin{bem} \label{remark_relation_to_others}
	The general idea of this estimator was stated in \cite{qiu1997nonparametric}, and in \cite{garlipp2007robust} the estimator was made precise by changing the order of scaling and rotating.
	As initially proposed by \cite{qiu1997nonparametric} one could also use a Gasser--M\"uller type version of the rotational difference kernel estimator defined by
	\begin{align*}
	\hat{\text{M}}_n  (\pvec{z};\psi,h) =  \sum\limits_{i_1,i_2=1}^n \yionetwo 	\int_{\triangle_{i_1,i_2}}  K(\pvec{z}-\pvec{u};\psi,h) \du,
	\end{align*}
	where $\triangle_{i_1,i_2} =   [(i_1-1)/n, i_1/n   ) \times [(i_2-1)/n, i_2/n  )$ denotes a design square.
	However, due to the rectangular fixed design it is easier to analyze the contrast function in (\ref{defi:contrast_emp}) instead.
\end{bem}

\begin{bem} [Several jump curves]\label{rem:severaljumpcurves}
{Also of interest for applications is an} extension to several jump-location-curves, viz.
$$
	m_\phi(x,y)  = m(x,y) + \sum_{k=1}^J j_{\tau_k,\phi_k}(x,y),\\
$$
and, for all $k \in \{ 1,\ldots,J\}$,
$$
	j_{\tau_k,\phi_k}(x,y) =  \tau_k(x) \mathbf{1}_{[0,\phi_k(x)]}(y),
$$	
	where $J \in \N$, $m:[0,1]^2 \to \R$ is the smooth part of the image,  $\tau_k:[0,1] \to \R_+$ are the jump-height curves and $\phi_k:[0,1] \to (0,1)$  are the jump-location-curves.	 In order to guarantee identification, one has to assume that the images of the jump curves are well separated, {i.e.,} that $ B_\rho(\phi_i[0,1]) \cap B_{\rho}( \phi_j[0,1]) = \emptyset$ for $i,j \in \{1,\ldots,J\}$ with $i\neq j$ and some $\rho>0$,  
	where $B_\rho(A)$ denotes the $\rho$-neighborhood of a set $A\subseteq\R$. 
	The theoretical concepts below should be transferable with additional notational effort (see, {e.g.,} Remark 2.6 in \cite{qiu1997nonparametric}) and confidence bands could be constructed based on a Bonferroni-correction.
	Section~\ref{sec:simulation} illustrates this idea by means of a simulated data example.
\end{bem}

\begin{bem}[Closed jump curves] \label{rem:closedjumpcurves}
	The rotated difference kernel method can also be applied to closed curves as follows. 		
	Let $G \subseteq \{ x\in \R^2 : \ \|x\|_2\leq 1 \} $ be a nonempty, compact, convex set with $0\in \text{int}(G)$ and assume that the image function is given by 
$$
	m_G \{ x_\gamma(t)\} = m\{x_\gamma(t)\} + \tau(\gamma) \mathbf{1}_G \{ x_\gamma(t)\},
$$
	where $m:[-1,1]^2 \to \R$ is the smooth part, $\tau:[0,2\pi) \to \R_+$ is the jump height and the location $x_\gamma$ is given by  $x_\gamma:[0,1]\to[-1,1]^2: t \mapsto t(\cos(\gamma),\sin(\gamma))^\top$ with $\gamma \in [0, 2 \pi)$. 		
	Now for $\gamma \in[0,2\pi),\psi \in[-\pi/2,\pi/2]$, consider the estimator 
	\begin{align*} 
	 (\estphi(\gamma), \estpsi(\gamma) )
	\in  \argmax{   t \in [h,1-h], \psi \in[-\pi/2,\pi/2] } \sum\limits_{i_1,i_2=1}^n \yionetwo 	K(x_\gamma(t)-\xionetwo;\psi,h).
	\end{align*}
	While details need to be worked out, our asymptotic analysis below should be transferable to this model. $\hphantom{BLA}$	
\end{bem}
 
 \section{Asymptotic theory} \label{sec:main_results}

 For our asymptotic analysis, we require the following set of assumptions. 
 \begin{assump}[{Errors}] \label{assumption:errors}
 	The errors $\eionetwo$ are square-integrable, centered, independent and identically distributed random variables with common variance $\sigma^2.$
 	%
 	%	Moreover, there exists $\delta\in(0,1]$ such that $\E|\epsilon_1|^k<\infty$ for some $k>4/(2-\delta).$
 	Moreover, $\E|\epsilon_{1,1}|^5<\infty.$
 \end{assump}
 \begin{assump}[{Smoothness}] \label{assumption:smoothness}
 	We have that
 	$	\phi \in C^2[0,1],  \tau \in C^2[0,1]$ and $m\in C^2[0,1]^2.$
 \end{assump}
 \begin{assump}[{Bandwidth}] \label{assumption:bandwidth}
 	The bandwidth $h=h_n$ lies in the range 
 	$	[	 {C_1\ln (n)^\eta}/{n^{1/2}},  {C_2}/{n^{1/3}}]			$
 	for some finite constants $C_1,C_2>0$ and some fixed $\eta>1.$ 
 \end{assump}  
 
 \begin{assump}[{Kernel}] \label{assumption:kernel}
 	The kernel functions $K_1 $ and $K_2$ are three-times continuously differentiable with support in $[-1,1]$ and satisfy the following conditions. 
 	
 	\begin{enumerate}
 		\item [(a)]
		$K_1$ is symmetric, {i.e.,} $K_1(x) = K_1(-x)$, and $K_1(x)>0$ for $x\in (-1,1)$. Further, $K_1$ satisfies, for all $j \in \{ 0, 1, 2\}$,
 		\begin{equation*}
 		%			\label{eq:kernelK1}
 		\int_{[-1,1]} K_1(x) \dx =1 \quad \text{and} \quad K_1^{(j)}(-1)=K_1^{(j)}(1) = 0.
 		\end{equation*}
 		\item [(b)]
		$K_2$ is an odd function, {i.e.,} $K_2(x) = - K_2(-x)$, in particular $K_2(0)=0$, and satisfies, for all $j \in \{ 0, 1, 2\}$,
 		\begin{align*}
 		%\label{eq:kernelK2}
 		\int_0^1 K_2(x) \dx & = - \int_{-1}^0 K_2(x)\dx  = 1, \ 
 		K_2^{(1)}(0) > 0\quad \text{and} \quad  K_2^{(j)}(-1)=K_2^{(j)}(1)=0.
 		\end{align*}	
 	\end{enumerate}
 \end{assump}   
 \begin{assump}[{Kernel moment}] \label{assumption:extra_kernel}
 	$K_2$ satisfies $\int_0^1 xK_2(x) \dx =0.$
 \end{assump}
 \begin{bem} [Assumptions] \label{remark_assumptions}
 	\begin{itemize}
 		\item [(i)]
		Assumption~\ref{assumption:errors} could be relaxed to $\E|\varepsilon_{1,1}|^k <\infty$ for some $k>4/(2-\delta)$ and $\delta \in (2/3,1] ,$  as can be seen from the Gaussian approximation provided in the {Online Supplement; see (\ref{eq_approx_rio})} therein.
 		However, this would introduce a new parameter $\delta,$ which in turn has to be taken into account for the range of admissible bandwidths and thus for sake of convenience we require the slightly stronger assumptions of existing fifth moments which is sufficient for the bandwidth range in Assumption~\ref{assumption:bandwidth}.  
 		\item [(ii)]
		Assumptions~\ref{assumption:smoothness}--\ref{assumption:kernel}  are required for the asymptotic distribution of the jump-location estimate $\estphi$ in (\ref{definition_estimator_phi}), whereas  Assumption~\ref{assumption:extra_kernel} is only needed for the faster rates of convergence of the estimates $\estpsi$ and $\esttau.$
 		Under Assumption~\ref{assumption:smoothness} the moment properties of $K_1$ resp.\ $K_2$ in Assumption~\ref{assumption:kernel} (resp.\ \ref{assumption:extra_kernel}) eliminate lower order terms in Taylor expansions, similarly as for higher order kernels in standard kernel smoothing.
 		We also point out that Assumption~\ref{assumption:bandwidth} corresponds to an undersmoothing for the continuous part of the image function. 
 		In particular, we require $nh_n^2 \to \infty$ and $n h_n^3 \to \text{const.},$ a standard assumption in bivariate nonparametric estimation for Lipschitz functions. Finally, note that $\ln (n)$ and $\ln (h\inv)$ are of the same order.  
 		\item[(iii)]
		 Our new conditions in Assumption~\ref{assumption:kernel} are necessary  for asymptotic normality of the estimates based on the rotated difference kernel method. 
 	\end{itemize}
 \end{bem}
 
 In the following, we fix some compact subinterval $\xset \subset (0,1)$. Let us start with uniform consistency of the estimators on $I$. 
 
 \begin{satz} \label{theorem:consistency}
 	In model \eqref{model}, under Assumptions~\ref{assumption:errors}--\ref{assumption:kernel} we have that
 	\begin{align*}
 	\sup_{x \in \xset} \  \big|\estphi(x) - \phi(x)\big|  =      O_P\{ {\ln (n)^{1/2}}/{n}\}. 
 	\end{align*}
 	If in addition Assumption~\ref{assumption:extra_kernel} holds, then
 	\begin{align*}
 	\sup_{x \in \xset} \  \big|\estpsi(x) - \psi(x)\big|  = O_P\{  {\ln (n)^{1/2}}/{(nh) }\}, \quad 
 	&\sup_{x \in \xset} \  \big|\esttau(x) - \tau(x)\big| = O_P\{ {\ln (n)^{1/2}}/{(nh)}\}.
 	\end{align*}
 \end{satz}
 \begin{bem} [Rates of convergence and optimality] \label{remark_rates_of_conv}
 	Theorem~\ref{theorem:consistency} shows that with our kernel assumptions, {i.e.,} Assumption~\ref{assumption:bandwidth}, the estimator for the jump location achieves the optimal  rate of convergence up to a logarithmic factor, see Theorem~3.3.1 in \cite{korostelev1993minimax} for the corresponding lower bound. Note that this is faster than the rates in Remark~2.4 in \cite{qiu1997nonparametric}, who lacks the kernel conditions in his asymptotic analysis. 
 	
 	Although there seem to be no rigorous results on lower bounds available, it is reasonable to believe that the rates for the jump-slope curve and the jump-height curve are optimal up to a logarithmic term as well, since these correspond to the uniform rate for estimating Lipschitz continuous functions. 	
 \end{bem}
 
 Now we turn to asymptotic normality of $\estphi(x)$ for fixed $x$. To state a version which does not involve a bias correction, we require the maximizer of the deterministic version of the contrast function
 \begin{align*} 
 \phi_n(x)	
 \in  \argmax{   y \in [h,1-h] } \max_{\psi \in[-\pi/2,\pi/2]}  \E\{ \hat{\text{M}}_n((x,y)^\top;\psi,h)\}.	
 \end{align*}
 
 \begin{satz} \label{theorem:main_theorem}
 	In model \eqref{model}, under Assumptions~\ref{assumption:errors}--\ref{assumption:kernel} we have, for any $x \in I$, that as $n \to \infty$,
 	\[  n\{ \estphi(x) -\phi_n(x) \}
 	\rightsquigarrow        
 	\mathcal{N}  [0,  \sigma^2 \ \varphisco(x)/\{\varphihess(x)\}^{2}    ]   , \]
 	where
 		\begin{align} \label{defi_components_of_cov_matrix}
 		\begin{split}
 		\varphihess(x)& = \tau(x) \cos^2\{\psix\}
 		K_2^{(1)}(0), \\ 
 		\varphisco(x) & =  \sin^2\{\psix\} \int_{ [-1,1]^2}  \{ K_1^{(1)}(z_1)K_2(z_2) \} ^2  \mathrm{d} z_1  \mathrm{d} z_2   + \cos^2\{\psix\} \int_{ [-1,1]^2} \{ K_1(z_1)K_2^{(1)}(z_2)\}^2   \,\mathrm{d} z_1 \mathrm{d} z_2.
 		\end{split}	
 		\end{align}       
 \end{satz}
 We even obtain joint asymptotic normality and asymptotic independence of $(\estphi,\estpsi,\esttau)$; see the proofs section, Theorem~\ref{theorem:asymp_norm_of_all}. In (\ref{defi_components_of_cov_matrix}), $\varphihess(x)$ arises as a limit of the Hessian matrices, while $\varphisco(x)$ comes from the asymptotic variance of the score. 
 
  	\medskip	
 \begin{bem}[{Confidence intervals}]
 	In order to construct asymptotic confidence intervals, we choose a consistent estimate $\hat \sigma^2_n$ of the error variance $\sigma^2$; {see, e.g.,} \citet{munk2005difference} and the simulations in Section~\ref{sec:simulation}.  	
 	Given $\alpha \in (0,1)$, we obtain an asymptotic level $1-\alpha$ confidence interval for $ \phi_n(x)$ by {computing}
 		\begin{align} \label{eq:confidence_interval}		
 		\left[   \estphi(x) - \frac{\hat \sigma_n   \sqrt{\varphiscoest  (x)} q_{1-\alpha/2}}{n\varphihessest (x)  } , \  \estphi(x) + \frac{\hat \sigma_n   \sqrt{\varphiscoest (x) }q_{1-\alpha/2}}{n\varphihessest (x)  }\right], 	
 		\end{align}
where $q_\beta = q_\beta( \mathcal{N}_1(0,1))$ is the $\beta$-quantile of $\mathcal{N}(0,1)$ and where we substitute estimators for the unknown parameters in (\ref{defi_components_of_cov_matrix}),  
 		\begin{equation}
 		\begin{aligned} \label{estimates_asymptotic_var_phi}
 		\varphihessest(x) &=  \esttau(x) \cos^2 \{\estpsi(x)\} K_2^{(1)}(0), \\
 		\quad
 		\varphiscoest(x) &= \sin\big(\estpsi(x)\big)^2 \int_{ [-1,1]^2}  \{K_1^{(1)} (z_1) K_2 (z_2) \}^2  \dzone  \dztwo  + \ \cos\{\estpsi(x)\}^2   \int_{ [-1,1]^2}  \{ K_1 (z_1) K_2^{(1)} (z_2) \}^2   \dzone  \dztwo.
 		\end{aligned}
 		\end{equation}	
 \end{bem}
 	
 \begin{bem} [Bias correction] \label{remark_bias_correction}							
 	The proofs show that there exists a bounded sequence $b_n$ such that, as $n \to \infty$,	
 	\[  n\{ \estphi(x) -\phi(x) \} -b_n
 	\rightsquigarrow        
 	\mathcal{N} [0,  \sigma^2 \ \varphisco(x)/\{\varphihess(x)\}^{2} ]  . 
	\]				
 	An explicit bound on $b_n$ based on the Lipschitz constants of $K,m,\phi$ and $\psi$ as well as on  $ \norminf{\tau}$ would be available by an explicit estimation of the second-order term and the Riemann-sum error term, see Lemma~\ref{lemma:asympbias} and the expansion (\ref{defi_repr_taylor}), and hence bias-corrected conservative confidence intervals could be constructed in a similar fashion as in \citet{eubank1993confidence}. Indeed, from \citet{korostelev1993minimax} it is known that the discretization bias for jump-curve estimation in a deterministic design is not negligible theoretically. However, our simulations show that	the confidence intervals in (\ref{eq:confidence_interval})  can typically be used for the actual parameter, and a bias correction is not necessary, indeed, it makes the intervals quite conservative.  	 \end{bem}

 Now let us turn to the construction of uniform confidence sets. 
 For independent standard normal random variables $\xi_{1,1}, \ldots, \xi_{n,n},$ independent of $\yionetwo,$  consider the process
 \begin{align} \label{defi:score_processes_bootstrap}
 \begin{split}
 \zbootphi( x ) & =  \frac{1}{ nh \{ \varphiscoest(x)\}^{1/2} } \\
 &\quad  \times \sumieiz \xi_{\ionetwo} 
 \scalarp{   (\nabla K) ( h\inv \rotmat_{-\estpsi(x)} \{ (x,\estphi(x))^\top -\xionetwo \} ),(\sin \estpsi(x), \cos \estpsi(x))^\top }, 
 \end{split}
 \end{align}
 where $\nabla K(z_1,z_2) = (K_1^{(1)}(z_1) K_2(z_2) , K_1(z_1) K_2^{(1)}(z_2)   )^\top $  
 and by definition of the rotation matrix in (\ref{eq:rotation_matrix}),
 $$	\rotmat_{-\estpsi(x)} \{ (x,\estphi(x))^\top -\xionetwo\} = 
 \left( 	
 \begin{array} {c}
 \cos \{ \estpsi(x) \} (x-x_{i_1}) + \sin \{\estpsi(x)\}\{\estphi(x) -x_{i_2}\} \\
 \cos \{ \estpsi(x) \} \{\estphi(x) -x_{i_2}\} -\sin \{ \estpsi(x)\}(x-x_{i_1}) 
 \end{array}		
 \right).
 $$

 This process corresponds to the normalized score process evaluated at the estimates $(\estphi,\estpsi)$ with independent noise-variables $\xi_{\ionetwo}$. Furthermore, consider the maximum
 \begin{align} \label{eq:bootprocess_phi}
 \bootprocessphi  =  \sup_{x \in I} \big|\zbootphi( x ) \big|. 	
 \end{align}
 The quantiles of this process can be determined by bootstrap simulations.
 The following result is the basis for constructing uniform confidence sets. 
 
 \begin{satz} \label{theorem:main_theorem_uniform}
 	Consider model \eqref{model} under Assumptions~\ref{assumption:errors}--\ref{assumption:kernel}, and assume that $\hat \sigma_n$ is an estimator for $\sigma$ which satisfies $\Pr (| {\hat\sigma_n}/{\sigma}  -1 |\geq s_n) = o(1)$ for some sequence {$s_n = o \{ 1/ \ln (n)\}$}.
 	Then for $\alpha\in(0,1)$, one has $q_{1-\alpha}(\bootprocessphi) = O_P\{(\ln  n )^{1/2}\}$ and for any sequence $t_n=o(1)$ such that $t_n \sqrt{\ln (n)} \to \infty$,  we have that
 	\begin{equation}\label{eq:asympconvband}
 	\liminf_n \ \Pr \left[	\sup_{x\in \xset} \Big| { \varphihessest(x) \{\estphi(x) - \phi(x) \}  }/{\sqrt{\varphiscoest  (x)} }   \Big| \leq  {(1+t_n) \hat \sigma_n q_{1-\alpha}(\bootprocessphi)}/{n}		\right] \geq 1-\alpha.
 	\end{equation}
 \end{satz}
 
 \begin{bem}[Asymptotic confidence bands]
 	Given $\alpha\in(0,1)$ a confidence band for $\phi$ which is asymptotically conservative at level $1-\alpha$ is given by 
 	\begin{align} \label{eq:confidence_bands_phi}
 	\{ [c_{\phi,u}^-(x),c_{\phi,u}^+(x)] : \ x\in I   \}, 
 	\quad c_{\phi,u}^\pm(x) = \estphi(x) \pm {(1+t_n) \hat \sigma_n \sqrt{ \varphiscoest(x)}  q_{1-\alpha}(\bootprocessphi) }/\{n \varphihessest(x)\}.
 	\end{align}
 Note that this is a version of the asymptotic {point-wise} confidence interval in (\ref{eq:confidence_interval}), corrected by the logarithmic factor $q_{1-\alpha}(\bootprocessphi)$ for uniform coverage. Further, the uniform confidence bands also directly apply to the actual parameter $\phi$, the price to pay being that they are asymptotically conservative. Similarly as in Theorem~\ref{theorem:main_theorem_uniform}, one can construct uniform confidence bands for $\psi$ and $\tau$. We sketch this in Section~\ref{sec:sketch_uniform_cb_psi_tau}.	
 \end{bem}
 
 \begin{bem} [Choice of $t_n$] \label{remark_choice_of_tn}
 The factors $t_n$ can be interpreted as a bias correction, and are essential for the validity of (\ref{eq:asympconvband}), even though they do not affect the rate of convergence. We will conduct a sensitivity analysis for the choice of $t_n$ in Section~\ref{sec:simulation}.  $t_n$ ensures that the confidence bands are asymptotically valid, {i.e.,} the significance level $1-\alpha$ is kept asymptotically, the assumption on $t_n$ in Theorem~\ref{theorem:main_theorem_uniform} prevents a straightforward analysis of asymptotically correct coverage of the confidence bands as in Corollary~3.1 of~\cite{chernozhukov2014anti}. 		
 \end{bem}

\section{Simulations}  \label{sec:simulation}
 
 In this section we investigate the finite sample properties of the proposed asymptotic confidence sets for the location $\phi(x)$ of the edge as well as of the estimator 
 \begin{equation}\label{eq:stdest}
 \{ \hat \sigma_n^2  \times  \varphiscoest(x) / \varphihessest(x)^2\}^{1/2}
 \end{equation} 
 for the asymptotic standard deviation  of $n\estphi(x)$ in Theorem~\ref{theorem:main_theorem} using (\ref{estimates_asymptotic_var_phi}). Further, we also investigate the bias in the estimation of the edge when using a deterministic rectangular grid. 
 
 \subsection{Simulation setup}
 
\begin{figure} [t!]
 % FIGURE 2
 	\centering
 	{\includegraphics[width=1\linewidth]{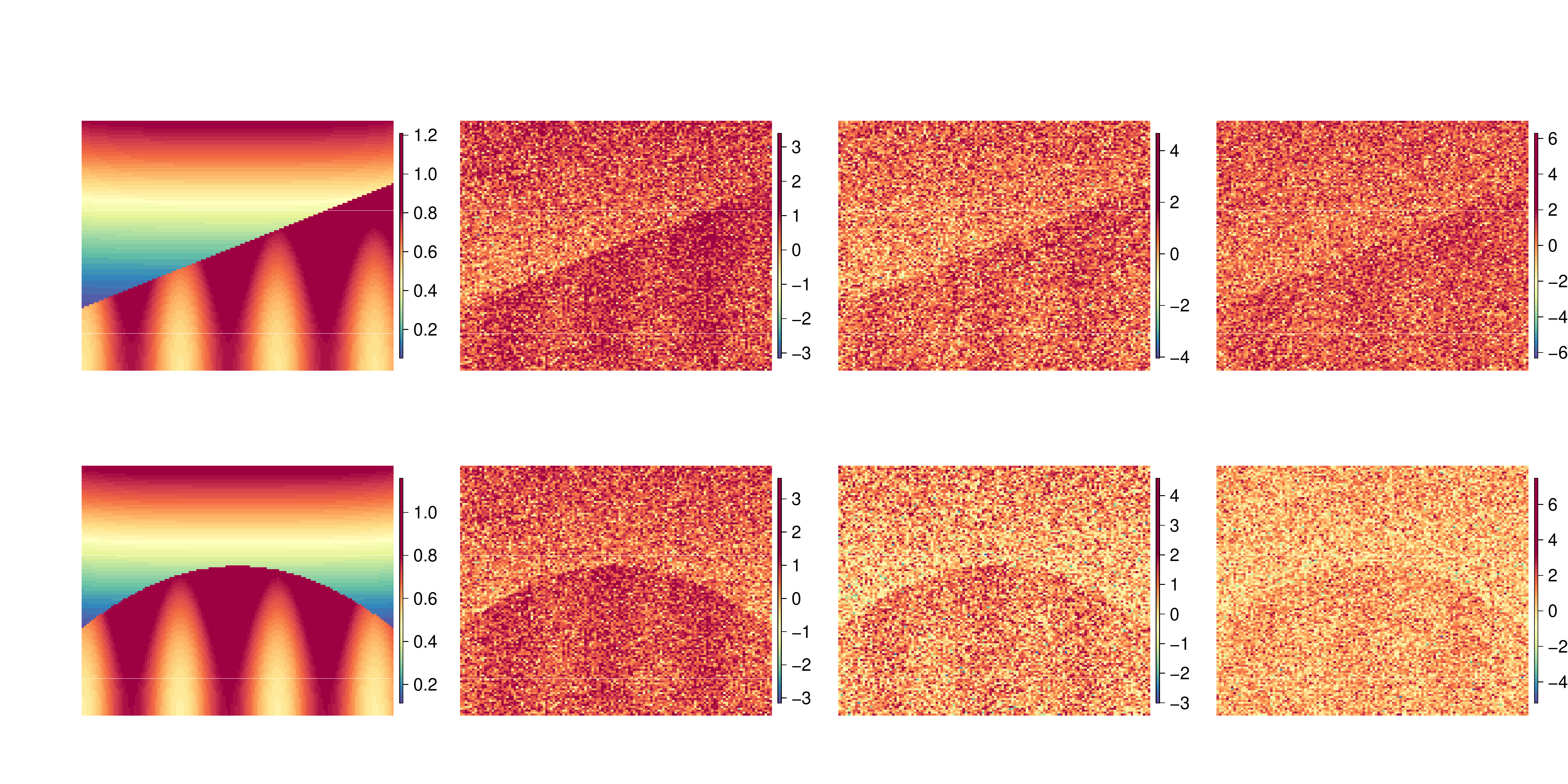}}
 	\caption{Upper panel: Image function $m_{\phi_1}$ with  $\tilde{\sigma}=0$, $0.5$, $0.7$, $0.9$ (from left to right) and $n=128$, respectively. Lower panel:  Image function $m_{\phi_2}$ with $\tilde{\sigma}=0$, $0.5$, $0.7$, $0.9$ (from left to right) and $n=128$, respectively.} 
 	\label{fig:heatplot_quad}
 \end{figure}

 We choose
 $$ 
 m(x,y) = \sin(y^2) \cos\{(x - 1/2 )^2\},\quad  \tau(x) = 3 \sin^2(10 x)/10 + 1/2,
 $$
 as background image and jump height, and consider the following two edge functions
 $$
 \phi_1(x)= 1/4 +  x/2, \quad \phi_2(x)=-\left(x-1/2\right)^2+3/5,
 $$
from which we form the regression functions $m_{\phi_1}$, $m_{\phi_2}$, according to model (\ref{assumption_regression_function}). Further, we choose $\epsilon \sim t_{10}(0,\tilde{\sigma})$, {i.e.,} a {Student $t$} distribution with location parameter zero, scale parameter $\tilde{\sigma}$ and {$10$} degrees of freedom. Thus, the noise-level is given by $\sigma = \tilde{\sigma}\sqrt{10/8}$. For illustration purposes, we display observations from the two models in Figure~\ref{fig:heatplot_quad} for grid size $n^2=128^2$.  We use the kernels
$$
 K_1(x)= C_1 \exp\{-(1-x^2)\inv\} \mathbf{1}_{[-1,1]}(x) ,\quad K_2(x) = C_2 \exp\{-(1-x^2)\inv\} (x^3-x) \mathbf{1}_{[-1,1]}(x), 
$$
where $C_1$ and $C_2$ are normalizing constants such that $\int K_1 =1$ and $\int_0^1 K_2 =1$, which satisfy Assumption~\ref{assumption:kernel}.

For the asymptotic results the bandwidth $h$ need not be chosen according to unknown smoothness parameters, and hence bandwidth selection is not a serious issue from a theoretical point of view. One could choose the bandwidth according to a selection rule like cross-validation or Lepski's rule to optimally estimate the background image $m$. \citet{qiu2005image} discusses simpler, more heuristic alternatives, one of which is to choose the window so that it contains approximately 100 design points. Although this is certainly not a universal rule which works for any $n$, we achieve reasonably good results with this approach for our grid sizes of $n^2 \in \{128^2, 196^2, 256^2\}$. In repeated simulations we use $500$ repetitions. 
 
 \subsection{Estimating the asymptotic standard deviation}
 
We start by investigating the numerical performance of the estimator (\ref{eq:stdest}) of the asymptotic standard deviation. To this end, we need to specify $\hat \sigma_n$, for which we choose squares of differences of all neighboring observation pairs properly normalized. The theory in \cite{munk2005difference} does not immediately apply when estimating on the full image. One possibility is to restrict estimation to a smooth part of the image. We also simulated the estimator $\hat \sigma_n$ of the standard deviation $\sigma$ separately, the results (not displayed) were also satisfactory. 	
 
 \begin{figure}[b!]
 % FIGURE 3
 	\centering
 	\subfigure
 	{\includegraphics[width=0.47\linewidth]{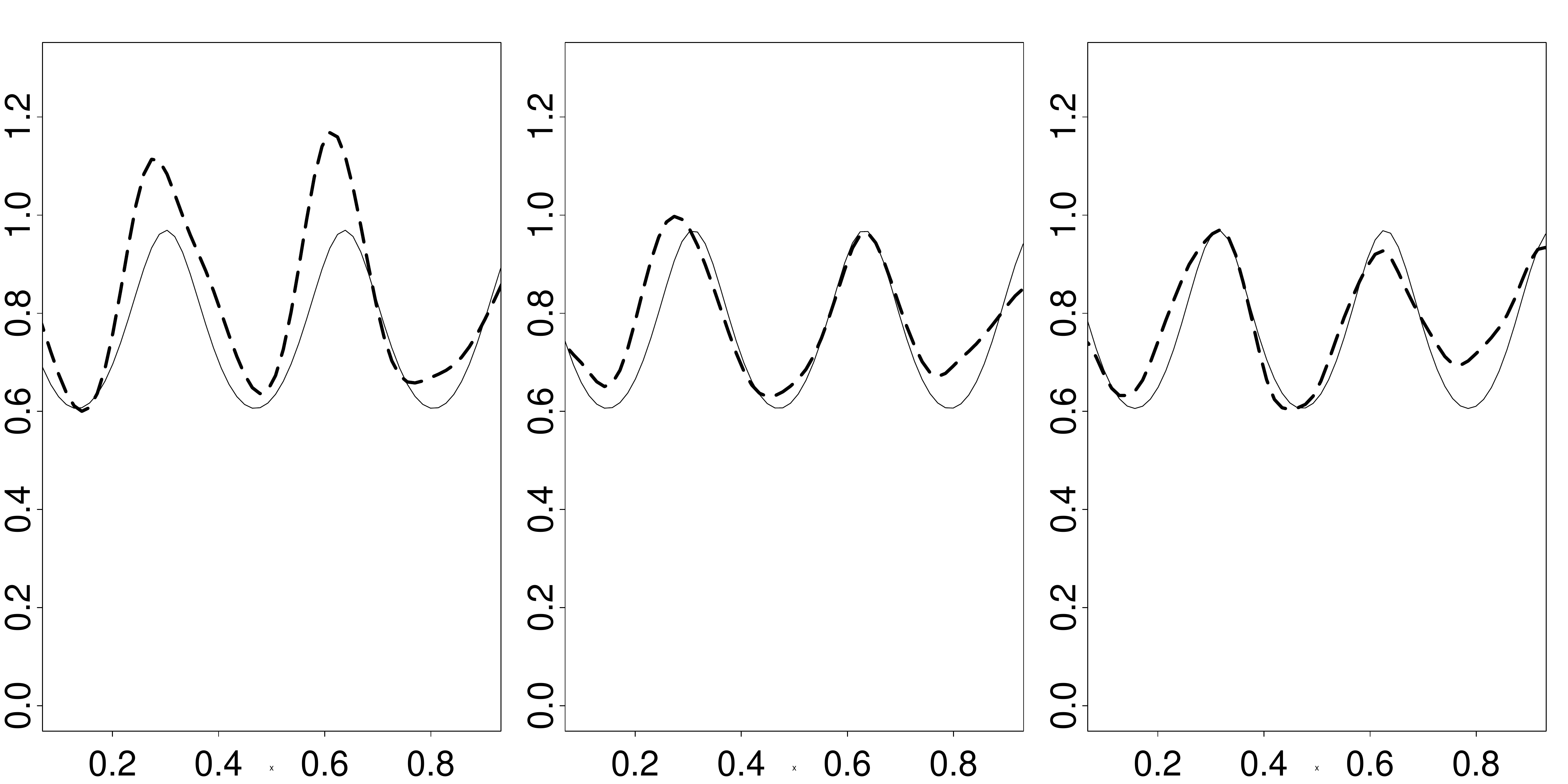}}
 	\subfigure
 	{\includegraphics[width=0.47\linewidth]{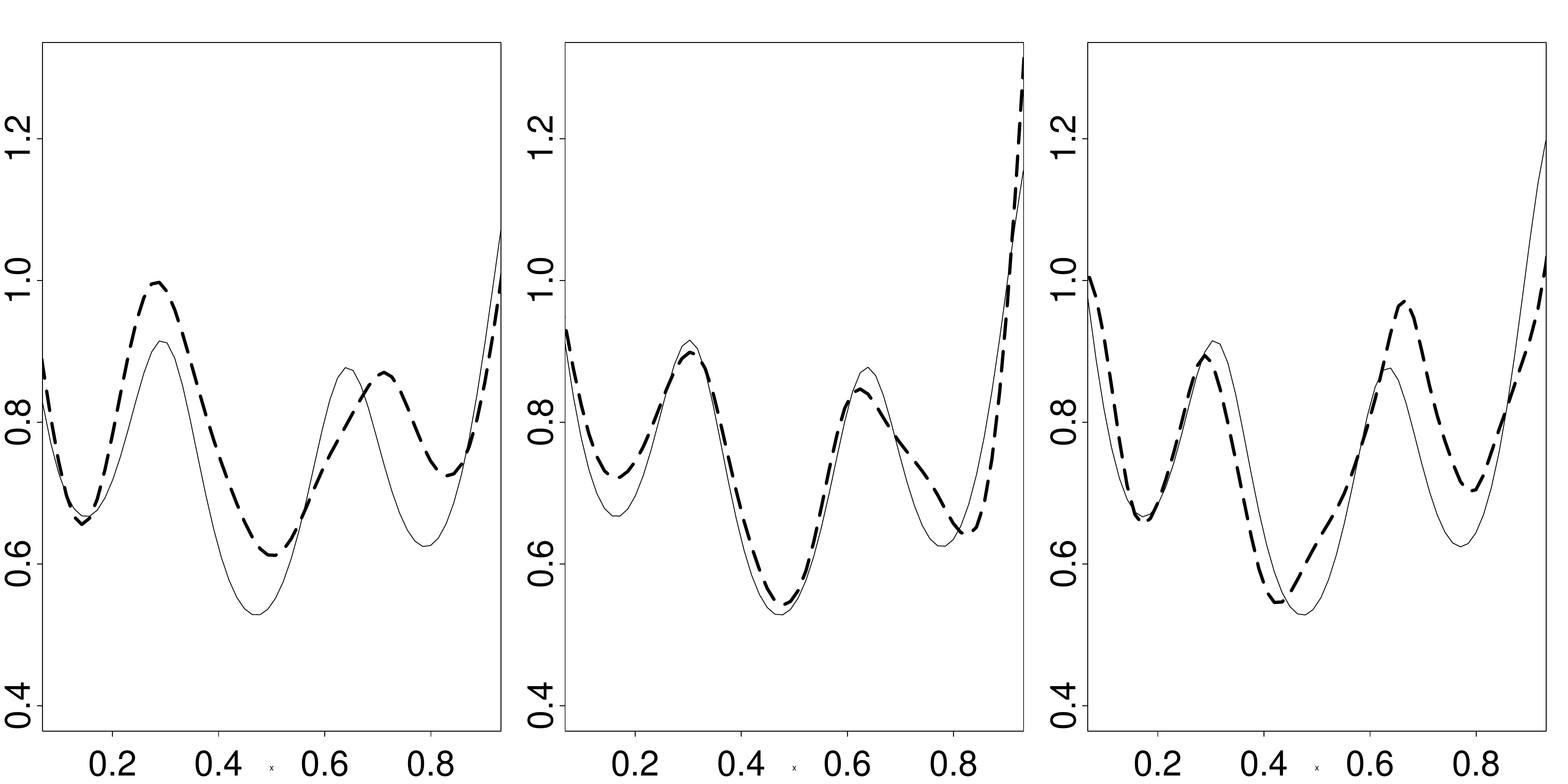}}
 	\caption{
 		Asymptotic standard deviation for $n \estphi$ (solid lines), and its estimates (dotted-dashed lines).
 		Three leftmost pictures: Standard deviation estimation for the image-function $m_{\phi_1}$  for $\tilde{\sigma}=0.5$ and $n=128$, $196$, $256$ (from left to right).
 		Three rightmost pictures: Standard deviation estimation for the image-function $m_{\phi_2}$  for $\tilde{\sigma}=0.5$ and $n=128$, $196$, $256$ (from left to right).}
 	\label{fig:covariance_estimation_phi}
 \end{figure}

 \begin{table}[t!]
 % TABLE 1		
 	\caption{Root of the MSE of the standard deviation estimate in scenario $\phi_1$ for some points $x$ if $\tilde{\sigma}=0.5$ or  $\tilde{\sigma}=0.9$. The last row indicates the mean of the RMSE for 64 points in the corresponding setting.}
	\bigskip
	
 	\centering
 		\begin{tabular}{r|rrr|r|rrr|r }
 			$\phi_1$	& \multicolumn{4}{c|}{  $\tilde{\sigma}=0.5$}  & \multicolumn{4}{c}{  $\tilde{\sigma}=0.9$}    \\ 
 			\hline
 			$x$  & $n=128$ & $n=196$ & $n=256$ & asymp.~sd  &   $n=128$ & $n=196$ & $n=256$ & asymp.~sd  \\ 
 			\hline
 			0.040 & 0.104 & 0.110 & 0.112 & 0.888 & 0.232 & 0.183 & 0.201 & 1.599 \\ 
 			0.142 & 0.101 & 0.080 & 0.071 & 0.611 & 0.209 & 0.154 & 0.119 & 1.100 \\ 
 			0.347 & 0.144 & 0.128 & 0.122 & 0.913 & 0.272 & 0.208 & 0.200 & 1.644 \\ 
 			0.449 & 0.107 & 0.084 & 0.075 & 0.617 & 0.183 & 0.163 & 0.136 & 1.111 \\ 
 			0.653 & 0.186 & 0.160 & 0.136 & 0.935 & 0.325 & 0.222 & 0.244 & 1.683 \\ 
 			0.858 & 0.160 & 0.132 & 0.123 & 0.725 & 0.275 & 0.206 & 0.204 & 1.305 \\ 
 			\hline
 			& 0.148 & 0.124 & 0.111 &   & 0.264 & 0.195 & 0.176 &   \\ 
 	\end{tabular}
 	\label{table:MSE_covariance_estimation_linear}
 \end{table}

 \begin{table}[b!]
 % TABLE 2
 	\caption{Root of the MSE of the standard deviation estimation for some points $x$ if $\tilde{\sigma}=0.5$ resp.\  $\tilde{\sigma}=0.9$ and scenario $\phi_2$. The last row indicates the mean of the RMSE for 64 points in the corresponding setting.}
\bigskip
 	\centering
 		\begin{tabular}{r|rrr|r|rrr|r }
 			$\phi_2$	& \multicolumn{4}{c|}{  $\tilde{\sigma}=0.5$}  & \multicolumn{4}{c}{  $\tilde{\sigma}=0.9$}    \\ 
 			\hline
 			$x$  & $n=128$ & $n=196$ & $n=256$ & asymp.~sd  &   $n=128$ & $n=196$ & $n=256$ & asymp.~sd  \\ 
 			\hline
 			0.040 & 0.137 & 0.148 & 0.153 & 1.149 & 0.274 & 0.244 & 0.243 & 2.069 \\ 
 			0.142 & 0.130 & 0.087 & 0.087 & 0.691 & 0.226 & 0.170 & 0.158 & 1.244 \\ 
 			0.347 & 0.142 & 0.124 & 0.132 & 0.840 & 0.257 & 0.196 & 0.202 & 1.511 \\ 
 			0.449 & 0.101 & 0.078 & 0.067 & 0.540 & 0.191 & 0.133 & 0.123 & 0.972 \\ 
 			0.653 & 0.193 & 0.156 & 0.143 & 0.859 & 0.317 & 0.238 & 0.182 & 1.547 \\ 
 			0.858 & 0.107 & 0.101 & 0.100 & 0.820 & 0.207 & 0.168 & 0.171 & 1.477 \\ 
 			\hline
 			& 0.143 & 0.119 & 0.111 &  & 0.263 & 0.200 & 0.182 &  \\ 
 	\end{tabular}
 	\label{table:MSE_covariance_estimation_quad}
 \end{table}
 
Next we present the results for (\ref{eq:stdest}). Figure~\ref{fig:covariance_estimation_phi} shows smoothed estimates for specific samples for grid sizes $n^2 \in \{128^2, 196^2, 256^2\}$ for the two edge functions $\phi_i$. Further, in Tables \ref{table:MSE_covariance_estimation_linear} and \ref{table:MSE_covariance_estimation_quad} we plot the square roots of the Mean-Squared-Error (RMSE) of the standard deviation estimates for the three sample sizes and two edge curves at various observation points $x$ based on  repetitions. For comparison purposes, the actual asymptotic standard deviation is given as well. One observes that the RMSE in most settings decreases as the number of grid points increases. Further, the magnitude of the RMSE as compared to the actual value of the asymptotic standard deviation is quite small for all~cases.

 \subsection{Confidence intervals and confidence bands}
 
 We investigate the coverage behavior and average width of (\ref{eq:confidence_interval}) as well as of (\ref{eq:confidence_bands_phi}) for the true jump-location curves $\phi_i$ in both settings $i \in \{1,2\}$. The results are summarized in Tables \ref{table:Coverage_low_noise} and \ref{table:Coverage_low_noise_uniform} for the noise-levels $\tilde \sigma =0.5$ and $\tilde \sigma =0.9$. The values in the tables of the {point-wise} confidence intervals correspond to the average of the respective quantity over 64 design points $x$. The quantile $ q_{\beta}(\bootprocessphi) $ for $\beta\in(0,1)$ was simulated based on a multiplier bootstrap sample of size $40{,}000$. Furthermore, as there is no explicit representation of the $t_n$-term we  have chosen it as given in Table~\ref{table:tn_low_noise} for the different scenarios. The $t_n$ decrease for increasing sample size $n$ and are of the same magnitude for both scenarios, i.e., for $\phi_1$ and $\phi_2$. By way of comparison, we give the values of $1/\sqrt{\ln (n)}$ for the different sample sizes $n$ as well. Especially, in the high-noise case the magnitude of our choice for $t_n$ is much smaller as this benchmark.

 \begin{table}[!htbp]
 	\caption{Choice of $t_n$ in  (\ref{eq:confidence_bands_phi}).}
	
	\bigskip
 	\centering
 		\begin{tabular}{r|rrr|rrr}
 			&\multicolumn{3}{|c|}{ $\tilde{\sigma}=0.5$}	& \multicolumn{3}{|c}{ $\tilde{\sigma}=0.9$} \\
 			&   $n=128$ &  $n=196$  &  $n=256$ &   $n=128$ &  $n=196$  &  $n=256$    \\
 			\hline
 			$\phi_1$	&  0.37   & 0.34 &  0.335   & 0.07 	&  0.001  	& 0 	\\
 			$\phi_2$	&  0.4   & 0.37  &  0.25   & 0.14 	&  0.1  	& 0.06	\\
 			\hline
 			$1/{\sqrt{\ln (n)}}$ &  0.45   & 0.44  &  0.42   &  0.45   & 0.44  &  0.42
 	\end{tabular}
 	\label{table:tn_low_noise}
 \end{table}	
 
Overall, the simulated coverage probabilities for the {point-wise} confidence intervals are reasonably close to their nominal values in all scenarios, and the intervals become narrower with increasing numbers of grid points. As expected from the theoretical developments, the uniform confidence bands are somewhat conservative  particularly in the high noise-level case.	
 
 \begin{table}[t!]
 % TABLE 4
 	\caption{Average coverage and width of the {point-wise} confidence intervals for the jump-location in (\ref{eq:confidence_interval}) for  $\tilde{\sigma}=0.5$ and $\tilde{\sigma}=0.9$ over 64 design points.}
	
	\bigskip
 	\centering
 		\begin{tabular}{r|rr|rr|rr|rr}
 			\multicolumn{4}{c}{ $\tilde{\sigma}=0.5$}	& & \multicolumn{4}{c}{ $\tilde{\sigma}=0.9$}		\\
 			\hline
 			\multicolumn{1}{c}{}	& \multicolumn{2}{c|}{95\% nominal coverage} & \multicolumn{2}{c|}{99\% nominal coverage}  & \multicolumn{2}{c|}{95\% nominal coverage} & \multicolumn{2}{c}{99\% nominal coverage}  \\
 			\multicolumn{1}{c}{} &	 coverage & width & coverage & width & coverage & width & coverage & width\\ 
 			\hline
 			\multicolumn{9}{c}{ $n=128$}			\\
 			\hline 
 			$\phi_1$	& 0.960 & 0.025  & 0.991 & 0.033 &  0.943  & 0.043  & 0.985  & 0.056 \\
 			$\phi_2$	& 0.958 & 0.025  & 0.991 & 0.033 &  0.946  & 0.042  & 0.983  & 0.056 \\ 
 			\hline
 			\multicolumn{9}{c}{ $n=196$}			\\
 			\hline 
 			$\phi_1$	& 0.958 & 0.016 & 0.992 & 0.021 & 0.944 & 0.028 & 0.990 & 0.037 \\ 		
 			$\phi_2$	& 0.951 & 0.016 & 0.988 & 0.022 & 0.943 & 0.028 & 0.988 & 0.037 \\			
 			\hline
 			\multicolumn{9}{c}{ $n=256$}			\\
 			\hline 	
 			$\phi_1$	& 0.958 & 0.012 & 0.992 & 0.016 & 0.941 & 0.021 & 0.997 & 0.028 \\ 
 			$\phi_2$	& 0.949 & 0.013 & 0.988 & 0.017 & 0.949 & 0.021 & 0.993 & 0.029
 			\\ 
 			\hline
 	\end{tabular}
 	\label{table:Coverage_low_noise}
 \end{table}	

 \begin{table}[b!]
 % TABLE 5
 	\caption{Average coverage and width of the uniform confidence bands for the jump-location in (\ref{eq:confidence_bands_phi}) for $\tilde{\sigma}=0.5$ and $\tilde \sigma=0.9$.}
	
	\bigskip
 	\centering
 		\begin{tabular}{r|rr|rr|rr|rr}
 			\multicolumn{4}{c}{ $\tilde{\sigma}=0.5$}	& & \multicolumn{4}{c}{ $\tilde{\sigma}=0.9$}		\\
 			\hline
 			\multicolumn{1}{c}{}	& \multicolumn{2}{c|}{95\% nominal coverage} & \multicolumn{2}{c|}{99\% nominal coverage}  & \multicolumn{2}{c|}{95\% nominal coverage} & \multicolumn{2}{c}{99\% nominal coverage}  \\
 			\multicolumn{1}{c}{} &	 coverage & width & coverage & width & coverage & width & coverage & width\\ 
 			\hline
 			\multicolumn{9}{c}{ $n=128$}			\\
 			\hline 
 			$\phi_1$	& 0.943 & 0.051 & 0.986 & 0.059 & 0.955 & 0.062 & 0.999 & 0.071 \\
 			$\phi_2$	& 0.959 & 0.049 & 0.984 & 0.059 & 0.945 & 0.065 & 0.999 & 0.075 \\ 
 			\hline
 			\multicolumn{9}{c}{ $n=196$}			\\
 			\hline 
 			$\phi_1$	& 0.945 & 0.032 & 0.988 & 0.038 & 0.954 & 0.039 & 0.999 & 0.045 \\ 		
 			$\phi_2$	& 0.957 & 0.032 & 0.986 & 0.039 & 0.948 & 0.042 & 1.000 & 0.050 \\			
 			\hline
 			\multicolumn{9}{c}{ $n=256$}			\\
 			\hline 	
 			$\phi_1$	& 0.949 & 0.025 & 0.987 & 0.028 & 0.951 & 0.030 & 0.993 & 0.032 \\ 
 			$\phi_2$	& 0.953 & 0.025 & 0.994 & 0.033 & 0.955 & 0.032 & 0.999 & 0.039
 			\\ 
 			\hline
 	\end{tabular}
 	\label{table:Coverage_low_noise_uniform}
 \end{table}	

\begin{figure} [t!]
 % FIGURE 4
 	\centering
 	\subfigure
 	{\includegraphics[width=20cm,height=5cm,keepaspectratio]{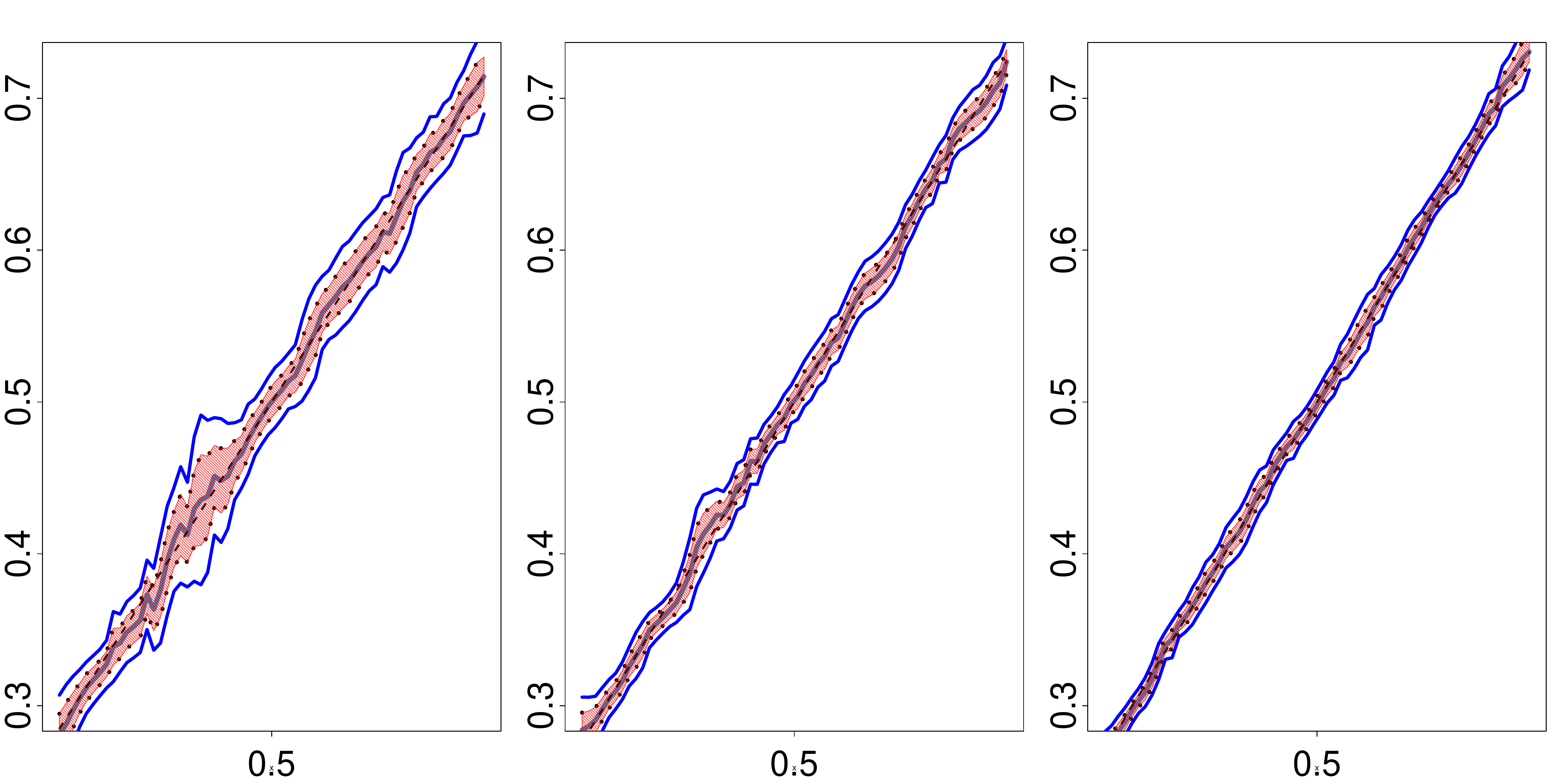}}
 	\subfigure
 	{\includegraphics[width=20cm,height=5cm,keepaspectratio]{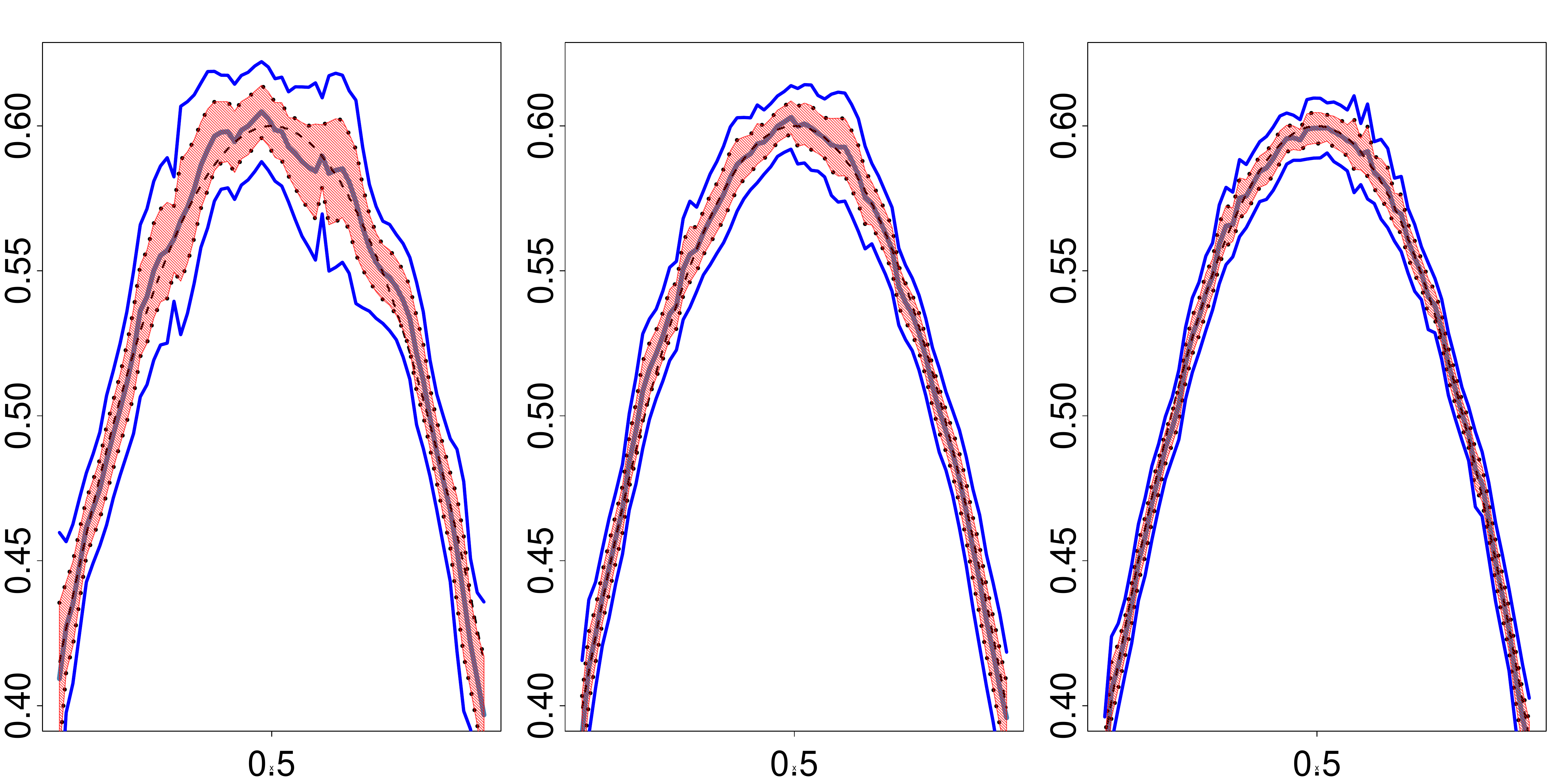}}
 	\caption{Top panel: 95\% Confidence intervals and estimate of the jump-location curve (shaded area and solid line within), uniform confidence bands (solid lines) and true jump-location curve $\phi_1$ (dashed lines inside shaded area) for $n=128,$ $196,$ $256$ and $\tilde{\sigma}=0.5$.
 		Lower panel:  95\% Confidence intervals and estimate of the jump-location curve (shaded area and solid line within), uniform confidence bands (solid lines) and true jump-location curve $\phi_2$ (dashed lines inside shaded area) for $n=128,$ $196,$ $256$ and $\tilde{\sigma}=0.5$.}
 	\label{fig:example_cis}
\end{figure}

Figure~\ref{fig:example_cis} illustrates the estimated curve as well as the confidence intervals and bands for $\phi_1$ and $\phi_2$ in the low noise-level for increasing numbers of grid points  for $\alpha= 0.05$, that is, asymptotic $95\%$ coverage probability. Apparently, the variability of the jump-location estimator decreases and the confidence intervals or bands become narrower. Besides, the confidence bands adapt to the shape of the {point-wise} confidence intervals as the width-terms only differ in the choice of the quantile.	
 
 \subsection{Sensitivity analysis of $t_n$ }
 
\begin{figure}[t!]
 % FIGURE 5
 	\centering
 	\includegraphics[width=0.7\linewidth]{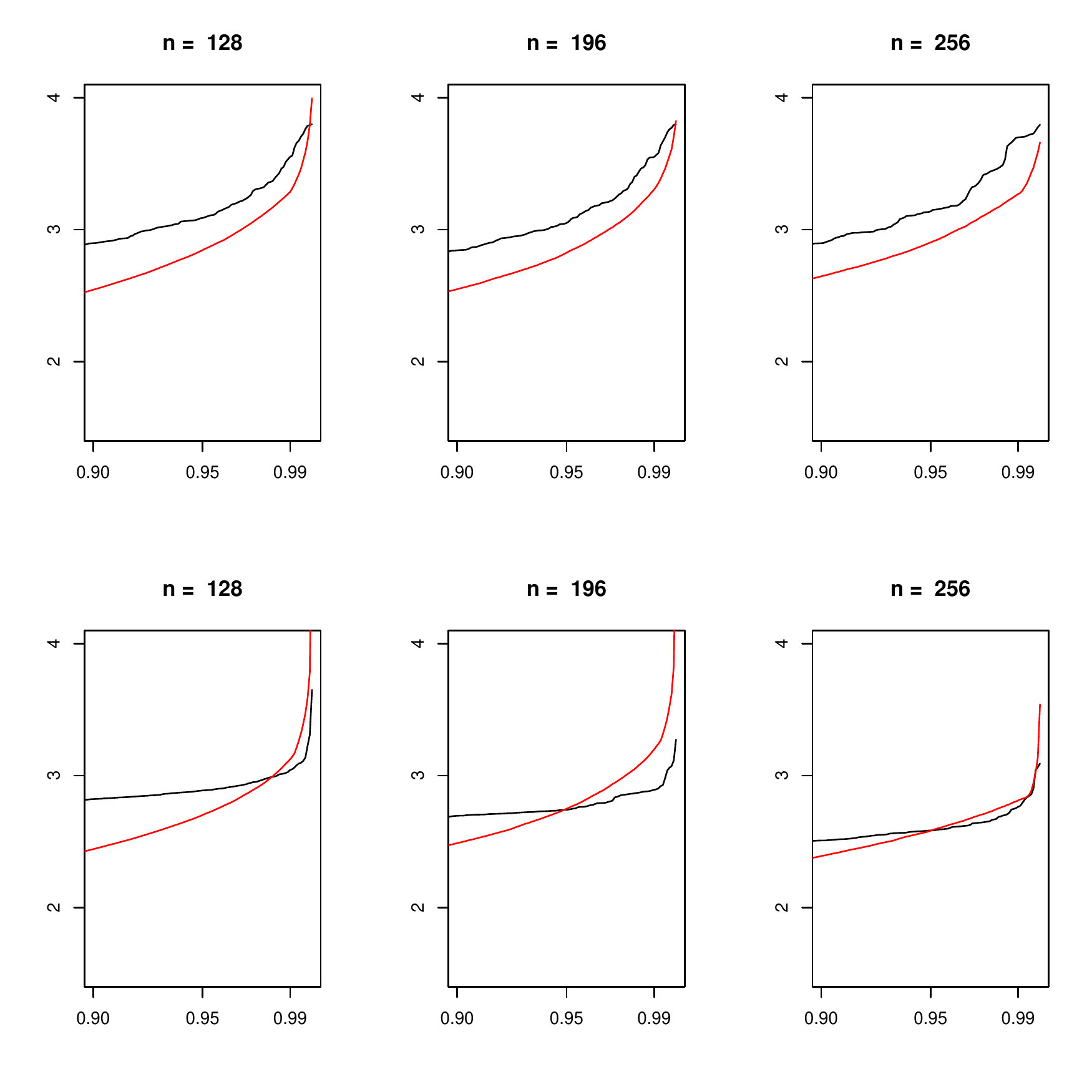}
 	\caption{Empirical quantiles of (\ref{sup_process}) for $\phi_1$ (black solid lines) and the corresponding bootstrapped quantiles (red solid line) for sample sizes $n \in \{128,196,256\}.$  Top panel: low-noise case $\tilde \sigma =0.5.$ Bottom panel: high-noise case $\tilde \sigma =0.9.$ }
 	\label{fig:tnlinear}
 \end{figure}

In order to get some insight on the role of the sequence $t_n,$ we display in Figure~\ref{fig:tnlinear} the empirical quantiles of 
\begin{align} \label{sup_process}
\sup_{x\in \xset} \ n \times  \Big| \frac{ \varphihessest(x) \{\estphi(x) - \phi(x) \}  }{ \hat \sigma_n \sqrt{\varphiscoest  (x)} }   \Big|			
\end{align}	
together with the simulated quantile curve $\{  q_{1-\alpha}(\bootprocessphi) : \alpha \in (0,1)\}$ for $\phi_1$ (the results for $\phi_2$ were similar). The simulated quantile curve is below the empirical quantile curve in particular for the low-noise level, while in the high noise-level the bootstrapped quantile curve is below the empirical quantile curve only for values smaller than a unique intersection point at approximately $0.95,$ and from then on above the empirical quantile curve.  
Thus, to guarantee an appropriate coverage of the confidence bands the bias correction term $t_n$ sequence must be of a higher magnitude in the low-noise level than in the high-noise level. 
Heuristically this is reasonable as the presence of the bias is much more noticeable in the low-noise case than in the high-noise case. The asymptotic choice $1/\sqrt{\ln  (n)}$ would lead to valid but conservative confidence bands.
 
 \subsection{Comparing bias and standard deviation}
 
Next, we investigate the order of the bias numerically and compare it to the standard deviation. Tables \ref{table:bias_var_ratio_linear} and \ref{table:bias_var_ratio_quad} contain the results for the ratio of the bias and the standard deviation for different design points $x$ in the low noise-level and the high noise-level case. The ratios are quite small, showing that the bias indeed is often of smaller magnitude than the standard deviation. 
  
 \begin{table} [b!]
 % TABLE 6
 	\caption{Ratio between computed bias and estimated standard deviation for different points $x$ in scenario $\phi_1$. The last line contains the average ratio over  64 design points for $\tilde{\sigma}=0.5$ and $\tilde{\sigma}=0.9.$}
	
	\bigskip
 	\centering
 		\begin{tabular}{r|rrr|rrr}
 			$\phi_1$ & \multicolumn{3}{c|}{  $\tilde{\sigma}=0.5$} & \multicolumn{3}{c}{  $\tilde{\sigma}=0.9$}   \\ 
 			\hline
 			$x$  & $n=128$ & $n=196$ & $n=256$ & $n=128$ & $n=196$ & $n=256$ \\ 
 			\hline
 			0.040 & 0.077 & 0.036 & 0.009 & 0.026 & 0.004 & 0.048 \\ 
 			0.142 & 0.001 & 0.037 & 0.013 & 0.032 & 0.035 & 0.010 \\ 
 			0.347 & 0.064 & 0.037 & 0.041 & 0.029 & 0.035 & 0.086 \\ 
 			0.449 & 0.043 & 0.005 & 0.013 & 0.014 & 0.093 & 0.085 \\ 
 			0.653 & 0.048 & 0.072 & 0.032 & 0.063 & 0.033 & 0.002 \\ 
 			0.858 & 0.018 & 0.004 & 0.022 & 0.003 & 0.046 & 0.075 \\ 
 			\hline
 			& 0.033 & 0.034 & 0.028 & 0.029 & 0.050 & 0.049 \\ 			
 	\end{tabular}
 	\label{table:bias_var_ratio_linear}
 \end{table}
 
 \begin{table} [t!]
 % TABLE 7
 	\caption{Ratio between empirical bias and estimated standard deviation for different points $x$ in scenario $\phi_1$. The last line contains the average ratio over 64 design points for $\tilde{\sigma}=0.5$ and $\tilde{\sigma}=0.9.$}
	
	\bigskip
 	\centering
 		\begin{tabular}{r|rrr|rrr}
 			$\phi_2$ & \multicolumn{3}{c|}{  $\tilde{\sigma}=0.5$} & \multicolumn{3}{c}{  $\tilde{\sigma}=0.9$}   \\ 
 			\hline
 			$x$  & $n=128$ & $n=196$ & $n=256$ & $n=128$ & $n=196$ & $n=256$ \\ 
 			\hline
 			0.040 & 0.091 & 0.121 & 0.051 & 0.037 & 0.021 & 0.042 \\ 
 			0.142 & 0.227 & 0.267 & 0.172 & 0.059 & 0.116 & 0.096 \\ 
 			0.347 & 0.181 & 0.216 & 0.030 & 0.041 & 0.144 & 0.073 \\ 
 			0.449 & 0.091 & 0.337 & 0.189 & 0.061 & 0.305 & 0.063 \\ 
 			0.653 & 0.120 & 0.095 & 0.079 & 0.020 & 0.081 & 0.117 \\ 
 			0.858 & 0.218 & 0.223 & 0.161 & 0.093 & 0.149 & 0.116 \\ 
 			\hline
 			& 0.160 & 0.217 & 0.138 & 0.075 & 0.130 & 0.101 \\ 
 	\end{tabular}
 	\label{table:bias_var_ratio_quad}
 \end{table}
 
 \subsection{Several jump curves}

 We briefly investigate numerically the extension indicated in Remark \ref{rem:severaljumpcurves} to two jump curves. The data are generated from
 \begin{align*}	
 m_{\phi}(x,y) &=	m(x,y) + j_{\tau_1,\phi_{1}}(x,y)+ j_{\tau_2,\phi_{2}}(x,y),
 \quad \tau_1 = \tau_2 \equiv 3/2, \\
 \phi_1(x) &= -\left(x-1/2\right)^2+21/50, \quad
 \phi_2(x) = \phi_1(x)+21/50,	
 \end{align*}
 and $m$ is as before, while $\sigma=0.1,$ $n=64,$ and $h=0.15.$  

\begin{figure}
 	\centering
 	\includegraphics[width=0.7\linewidth]{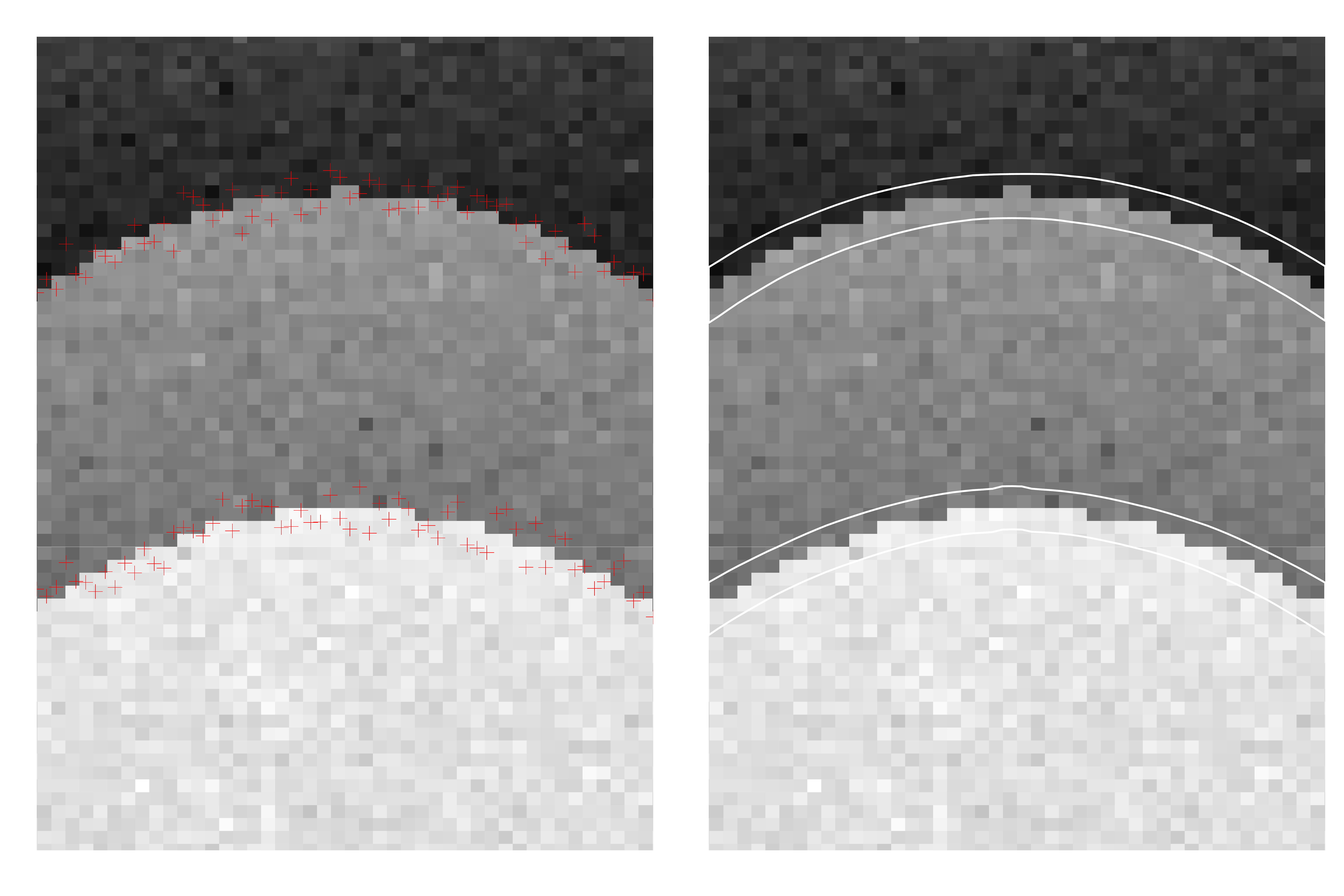}
 	\caption{Left picture: Noisy regression function and red crosses indicate the thinned out candidate set of change points. Right picture: Resulting Bonferroni-corrected 95\% confidence bands.}
 	\label{fig:manyedges}
 \end{figure}

In each strip $\{x\}\times [0,1],$ in a first step we choose as estimates the points with maximal contrast and having at least $h$-distance, see the left picture in Figure~\ref{fig:manyedges}. In the second step, we compute the maximizers in some $\delta$-neighborhood of the $y$-coordinate of the candidate points and construct the simultaneous confidence bands by a Bonferroni-correction. The results are displayed in the right picture of Figure~\ref{fig:manyedges}. Both edges are within the 95\% confidence bands and are of satisfactory width. 
 
 \subsection{Real data illustration}
  
\begin{figure} [ht!]
 	\centering
 	{\includegraphics[width=1\linewidth]{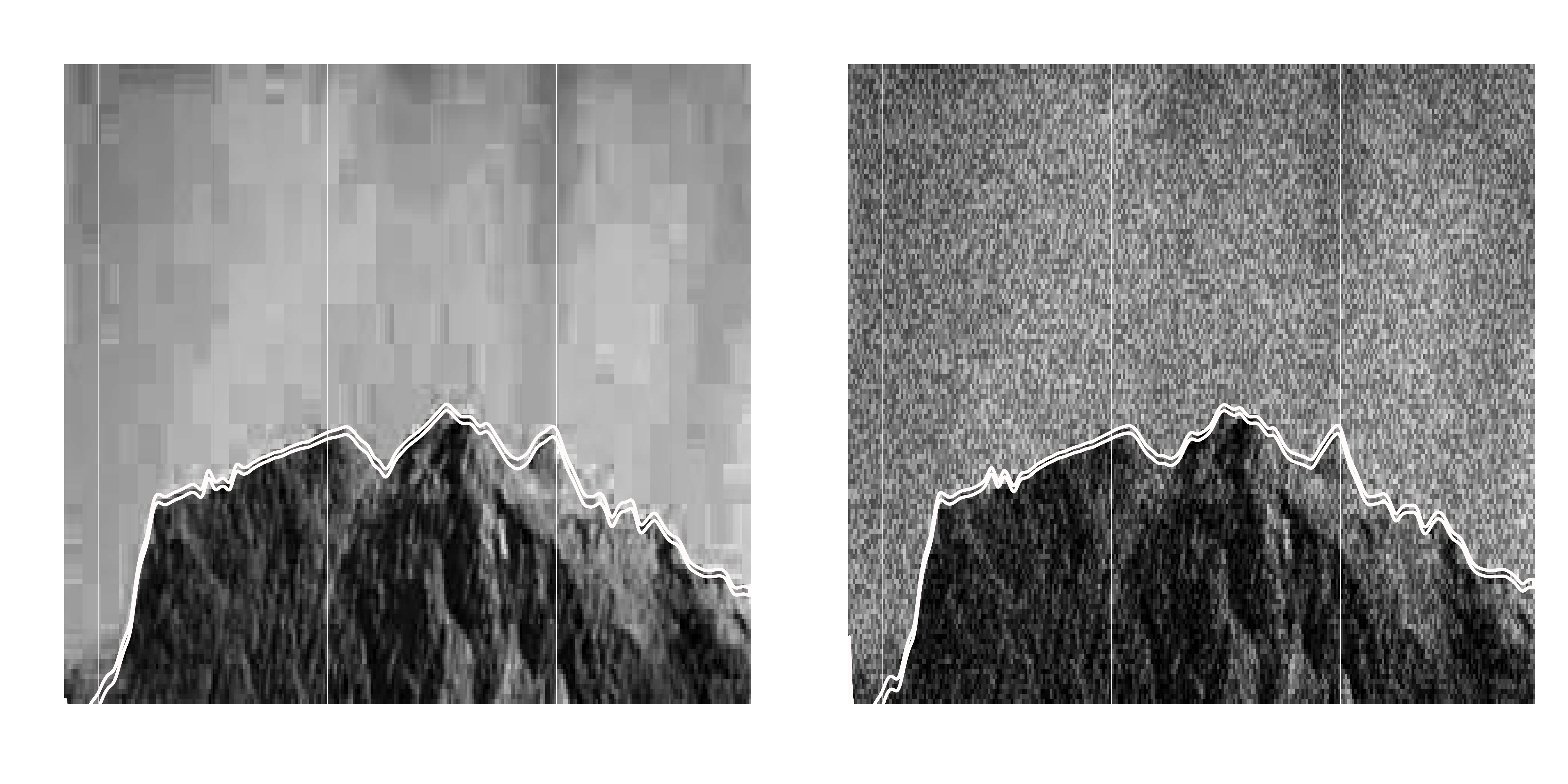}}
 	\caption{Left: $300 \times 128$ Gray-scale picture of a rock taken from a camera with reasonable ISO-configuration and with focus on the rock. Solid lines correspond to the 95\%-level confidence band for the  noisy picture on the right. Right:  $300 \times 128$ Gray-scale picture of a rock  taken from a camera with inappropriate ISO-configuration and with no focus. Solid lines correspond to the 95\%-level confidence band for this noisy picture.} 
 	\label{fig:rock-pictures}
 \end{figure}

 Finally we apply our method to two $300 \times 128$ Gray-scale real-data images, taken by camera by one of the authors. They contain the outline of a rock in front of a gray background. 
 Once, an appropriate ISO--configuration and focus on the rock and once inappropriate ISO--configuration of the camera and no focus at all are employed. 
 In both cases we apply our method to estimate the boundary of the rock and construct 0.95-level uniform confidence sets using $t_n =0$ and a naive estimator for the noise-level.
 \nopagebreak
 Figure~\ref{fig:rock-pictures} contains the results. 

 \nopagebreak
 The jump-location curve lies mostly inside the constructed confidence band and the width is quite satisfying, although the noise-level of the picture is rather low.
This illustration might be representative for more serious applications. For instance, a polar explorer  might intend to monitor the height of a glacier in order to analyze the effect of global warming. 
 Due to external effects such as bad weather conditions the pictures of the glacier could be noisy, so that its ridge needs to be estimated. The confidence bands are useful to assess whether the height of the glacier actually decreased. 	
 \nopagebreak

\section{Discussion}\label{sec:discuss}

\begin{text}
	In this paper we developed methods to construct asymptotic confidence sets for the jump curve in an otherwise smooth two-dimensional regression function, for which to the best of our knowledge no methods were previously available. As a further step, a combination of our results with the issue of edge detection would be desirable. A direct thresholding of the contrast function (\ref{defi:contrast_emp}) without prior localization as in (\ref{def:2D-deltafunktion.bivariat}) seems not to result in confidence sets which decrease at a near optimal rate, so that further methods would be required.  An extension of our results to the multivariate setting, especially to three dimensions would also be relevant. Further, images are often observed with blurring, that is, after convolution with a point-spread function. Thus, extensions to a deconvolution setting would also be of practical importance. 	
\end{text}

\section{Outline of proofs} \label{sec:proofs}

For future reference we collect additional notation that is used in the following. 

\medskip

\subsection{Notation}

We shall use the following notation. For a vector $\pvec{z}=(z_1,z_2)^\top \in \R^2$ we denote the coordinate projection onto the {$i$th} coordinate as $(\pvec{z})_i=z_i$ for $i \in \{ 1,2\}$. Furthermore, we write $\pvec{z}^\alpha=z_1^{\alpha_1} z_2^{\alpha_2}$ for $\alpha=(\alpha_1, \alpha_2)^\top \in \N^2$ and 
$$
f^{(\alpha)}(\pvec{z}) = \frac{\partial^{\alpha_1+ \alpha_2}}{\partial z_1^{\alpha_1} \partial z_2^{\alpha_2} }f(\pvec{z})
$$ 
for a function $f:\R^2 \to \R$. In particular, $\nabla  f (\pvec{z}) = (f^{(1,0)} , f^{(0,1)} )^\top(\pvec{z})$ and 
$$ \nabla \nabla^\top f(\pvec{z}) = \lb \begin{array} {cc}
f^{(2,0)} & f^{(1,1)} \\ f^{(1,1)} & f^{(0,2)}
\end{array}  \rb (\pvec{z})  .  
$$

We also write $\partial_{z_i} f(\pvec{z})$ for $f^{\pvec{e_i}^\top}(\pvec{z}),$ where $\pvec{e_i}$ is the {$i$th} canonical unit vector. If $g: X \to Y$ for $X,Y \subset \R,$ then we let $\epi(g)=\{  (x,y)^\top \in X \times Y  :g(x)\leq y \}$ be the epigraph of $g$. Let $A \triangle B$ denote the symmetric difference between two sets $A$ and $B$, {i.e.,} for any $A,B \subset \R^2$,
$$
A \ \triangle \ B  =  (A\backslash B) \cup (B \backslash A) = (A \cup B) \backslash (A \cap B).
$$

Write $\pvec{z} + A  =  \{\pvec{z} + \mathbf{y} : \mathbf{y} \in A\}$ for $\pvec{z} \in \R^2$ and a set $A \subset \R^2$. For sequences $a_n$ and $b_n$ we write $a_n \cong b_n$ if there exist constants $0<C_1<C_2$ and $n_0 \in \N$ such that $ C_1 \leq | a_n / b_n  | \leq C_2  $ for any $n\geq n_0$. Denote by $\norm{\cdot  }  $ a norm on  $\R^2$ as well as on $ \R^{2\times 2}$, where the dimension should be clear from the context. We only assume that the matrix-norm is compatible with the vector-norm, that is $\norm{\pvec{A} \pvec{z}  }\leq \norm{\pvec{A}  } \times  \norm{ \pvec{z}  }$ for a matrix $\pvec{A}$ and a vector $\pvec{z}$ and that the matrix-norm is submultiplicative, {i.e.,} $\norm{\pvec{A} \pvec{B}  }\leq \norm{\pvec{A}  } \times  \norm{ \pvec{B}  }$ for matrices $\pvec{A},\pvec{B}$. For a function $f:I \to \R^{i \times j}$ for $ i,j\in\{1,2\}$ with $i\geq j$ and $I$ is a compact subinterval of $(0,1)$ we define the uniform-norm as 
$$	\norminf{f} = \begin{cases}
\sup_{x\in \xset} \ \ \abs{f(x)} & \mbox{if }  i=j=1, \\
\sup_{x\in \xset} \norm{f(x)} & \text{otherwise.} 
\end{cases}			
$$

\subsection{Rescaling of the contrast function}

From now on, we shall always assume that $h$ is so small ($n$ is sufficiently large) that $I \subset [h,1-h]$. Define, for $x\in I$, the set of the rescaled parameters as 
\begin{align*}
\begin{split}
\maxset  = \{ w \in \R: \phi(x) + wh \in [ h, 1 - h ]  \} \times [-\pi/2,\pi/2],
\end{split}
\end{align*} 
and for $(x,\psi,w)^\top \in I \times \R^2$ the rescaled contrast function
\begin{align} \label{def:2D-deltafunktion.bivariat}
\empcrit(w, \psi;x) = \hat{\text{M}}_n( (x,\phi(x)+hw)^\top ;\psi,h)  = \hat{\text{M}}_n( \px + w h \pvec{e_2} ;\psi,h),
\end{align}
where  $\hat{\text{M}}_n$ is given in (\ref{defi:contrast_emp}) and $ \px = (x, \phi(x)  )^\top$  as well as $\pvec{e_2} = (0, 1)^\top$. With this, we have for the maximizers 
\begin{align} \label{defi_rescaled_maximizers}
\left(\estw(x),\estpsi(x) \right) \in \argmax{ (w,\psi)\in \maxset } \empcrit  (w,\psi;x)
\end{align}
that $ \estw(x) = \{ \estphi(x) - \phi(x)\} / h$. The parameter set  is denoted by $\paramset = \bigcup_{x\in I} \{x\} \times \maxset$. We also introduce the deterministic contrast function
$$	
\detcrit(w,\psi;x)= \E\{ \empcrit(w, \psi;x) \}.		
$$	

\subsection{Uniform consistency}

In this section we show uniform consistency over $\xset$ for the maximizers in (\ref{defi_rescaled_maximizers}).

\begin{prop} \label{prop:consistency}
Under Assumptions~\ref{assumption:errors}--\ref{assumption:kernel} we have, {as $n \to \infty$,} 
$$
\norminf{ \big(\estw(\cdot), (\estpsi -\psi)(\cdot) \big)^\top } \Pto 0. 
$$
\end{prop}

For the proof, we make use of the following adaptation of Theorem~5.7 in \cite{van2000asymptotic}, which is proved in {Section~\ref{sec:stochunifconv} of the Online Supplement}.

\begin{prop} \label{prop:vdv_consistency:uniform}
Assume that $\maxset \subset \R^d$ are compact sets and set $ \Theta = \cup_n \paramset,$ where $\paramset = \cup_{x\in I} \{x\} \times \maxset$. Let $\hat f_n: \R^d \times I \to \R$ be random functions and 	let $f: \R^d \times I \to \R$ be a deterministic function. Suppose that 
\begin{align}\label{eq:unifconsproof1}
	&\sup_{ (x,\theta) \in \paramset } |   \hat f_n(\theta,x) - f(\theta,x)  | \Pto 0, 
\end{align}
and that there exists a  map $\theta_0: I \to \Theta$ such that, for any $\epsilon>0$,
\begin{align}\label{eq:unifconsproof2}
	&\inf_{x \in I } \left[ f \{ \theta_0(x),x\} - \sup_{	\theta \in \R^d : |\theta-\theta_0(x)|\geq \epsilon } f(\theta,x) 	\right] >0.
\end{align}
Then for any estimator $\hat  \theta_n(x) \in  \maxset ,$ where $x\in I,$ which satisfies 
$$
\inf_{x \in I } \big[ \hat f_n \{ \hat \theta_n(x),x\} - \hat f_n \{ \theta_0(x),x\} \big]\geq  -o_P(1) 
$$
one has $	 \norminf{ \hat \theta_n   -  \theta_0   }		\Pto 0$ as $n \to \infty$.
\end{prop} 

In order to obtain Proposition~\ref{prop:consistency}, in the next two lemmas  we check the requirements of Proposition~\ref{prop:vdv_consistency:uniform}.  
The next lemma determines the limit function for $\empcrit(w,\psi;x)$ and shows that they fulfill the assumption (\ref{eq:unifconsproof1}). 

\begin{lemma} \label{lemma:pointwise_convergence} 
	Under Assumptions~\ref{assumption:errors}--\ref{assumption:kernel}, we have that
	$$  \sup_{(x,w,\psi)^\top\in \paramset } | \empcrit(w,\psi;x)  - \asympcrit\left(w,\psi;x\right)| = O_P(h) + O_P\{ \ln (n)^{1/2}  (nh)\inv\}, 
	$$	
	where
	\begin{align}
	\label{def:det-limit-crit-fct}
	\asympcrit(w, \psi;x) = \tau (x)\int_{ \Hr \{ \psix - \psi\} + w (\sin \psi, \cos \psi )^\top }  K(\pvec{z})
	\dz,
	\end{align}
	and, for any $\psi \in [-\pi, \pi]$, $\Hr(\psi) = \rotmat_\psi \{\R \times [0, \infty) \}$.
\end{lemma}

\begin{proof}[{\sl Proof of Lemma~\ref{lemma:pointwise_convergence}}]
	Split $\empcrit(w,\psi;x)$ into a smooth-image-related part, a jump-related part and an error-related part as follows
	\begin{align*}
	\detcrit(w,\psi;x)	
	&= \left(nh\right)^{-2}   \sum\limits_{i_1,i_2=1}^n m (\xionetwo ) K[h\inv \rotmat_{-\psi}\{ \px +hw \pvec{e_2} -\xionetwo\}]   \\
	&\quad + \left(nh\right)^{-2}   \sum\limits_{i_1,i_2=1}^n j_\tau\left(\xionetwo\right) K[ h\inv \rotmat_{-\psi}\{\px +hw \pvec{e_2} -\xionetwo\}]   \\
	&\quad + \left(nh\right)^{-2}   \sum\limits_{i_1,i_2=1}^n \eionetwo K[h\inv \rotmat_{-\psi}\{\px +hw \pvec{e_2} -\xionetwo\}]  \\
	& = S_n(w,\psi;x) + J_n(w,\psi;x) + E_n(w,\psi;x).
	\end{align*}
	
	Take
	$$ g_1(\pvec{z})=m(\pvec{z}), 
	\quad f(\pvec{z})= 4^{-1} \left( \begin{array}{cc}
	K(\pvec{z}) & K(\pvec{z})  \\ K(\pvec{z}) & K(\pvec{z}) 
	\end{array} \right), \quad g_2(\pvec{z})= g_3(\pvec{z}) =(1,1)^\top. 
	$$ 
{Lemma~\ref{lemma:S_n_term} in the Online Supplement} with $r_1=k=2$, $r_2=r_3=0$ and $j=1$ states that
	\begin{align*}
	S_n (w,\psi;x)	
	= m \{ \px + hw \pvec{e}_2\} \int_{[-1,1]^2} K\left(\pvec{z}\right) \dz 
	+ O (h ) 
	+ O\{ (nh )\inv\}
	= O (h ) + O \{(nh )\inv\}
	\end{align*}
	uniformly over $\paramset.$
	With the same functions $f,g_2,g_3$ as above, {one deduces} from {Lemma~\ref{lemma:J_n_term}} that 
	\begin{align*}
	J_n(w,\psi;x)
	= & \tau (x)  \int_{ \Hr \{ \psix-\psi \}  + w\left(\sin\left(\psi\right),\cos\left(\psi\right)\right)^\top}  K\left(\pvec{z}\right) \dz + O\left(h \right)  + O\{ (nh )\inv \} \\
	= & \asympcrit\left(w,\psi;x\right) + O\left(h  \right)  + O\{ (nh )\inv \},
	\end{align*} 	
	uniformly over $\paramset.$ 
	Finally, {Lemma~\ref{lemma:E_n_term}} with the same $f,g_2,g_3$ implies $E_n(w,\psi;x) = O_p\{ \ln (n)^{1/2} (nh)\inv\}$ uniform over $\paramset.$
	This concludes the proof of the lemma.
\end{proof}

In the next lemma, we rewrite the asymptotic form of the contrast and show that it has a unique well-separated maximum, that is the requirement (\ref{eq:unifconsproof2}).

\begin{lemma} \label{lemma:derivatives_asymptotic}
	Under Assumptions~\ref{assumption:smoothness} and  \ref{assumption:kernel}, one has that
	if $\psix-\psi = \pm \pi/2$ then $\asympcrit\left(w, \psi;x\right) = 0$ for all $w$, while otherwise  
	\begin{equation}\label{eq:asympcritrepres}
	\asympcrit\left(w, \psi;x\right)= -\tau (x)  \int_{-1}^1 K_1\left(y\right) \bar K_2\{a_x(\psi)y + b_x(w,\psi)\} \dy,
	\end{equation}
	where
	\begin{align*}
	\bar K_2\left(y\right)  = \int_{-1}^y K_2\left(t\right) \dt ,\quad  	
	a_x\left(\psi\right) =\tan\{\psix-\psi\}, \quad 
	b_x\left(w,\psi\right)= w \, \frac{\cos\{\psix\}}{\cos\{\psix-\psi\}}.
	\end{align*}
	This implies that for any $\epsilon>0$,
	\begin{equation}\label{eq:contrast}
	\inf_{x \in [0,1]}  \Big\{\asympcrit(0,\psix;x) -  \sup_{  (w,\psi)\in \R^2: \ \epsilon \leq \max( |w|, |\psix - \psi|  ) }  \asympcrit(w,\psi;x) \Big\} >0. 
	\end{equation}
	Moreover, $\asympcrit(0,\psix;x)= \tau(x)$.
\end{lemma}

\begin{figure}[t!]
% FIGURE 8
	\centering
	\includegraphics[width=0.8\linewidth]{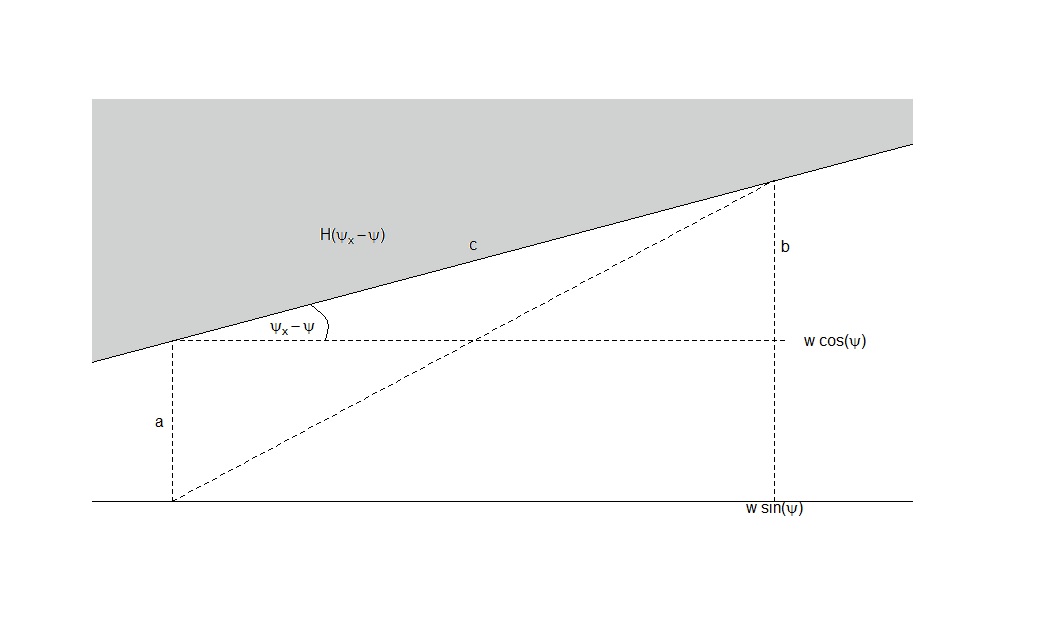}
	\caption{Gray area corresponds to $\Hr\{\psix-\psi\}$.}
	\label{fig:help_graphic_asymp_limit}
\end{figure}

\begin{proof}[{\sl Proof of Lemma  \ref{lemma:derivatives_asymptotic}}]
	If $\psix-\psi =  \pi/2$ then $\Hr\{\psix - \psi\} = (-\infty, 0] \times \R $. Since $\int_\R K_2 = \int_{-1}^1 K_2 = 0$, the statement $\asympcrit\left(w, \psi;x\right) = 0$ follows from (\ref{def:det-limit-crit-fct}), and similarly for $\psix-\psi =  -\pi/2$. 
	If $\psix-\psi \not= \pm  \pi/2$, we start by showing that 
	\[ 
	\Hr\{\psix - \psi\} + w\left(\sin \psi, \cos \psi\right)^\top = \Hr\{\psix - \psi\} + w \frac{\cos\{\psix\}}{\cos\{\psix-\psi\}} \, \pvec{e_2},
	\]
	for an illustration see Figure~\ref{fig:help_graphic_asymp_limit}.
	If $w=0$ the assertion is trivial. 	If $w\not=0$, we need to determine $a$ such that the vector $ a \, \pvec{e_2}$ is on the boundary of the set $\Hr\{ \psix - \psi\} + w\left(\sin \psi, \cos \psi\right)^\top$; see Figure~\ref{fig:help_graphic_asymp_limit}. 
	Setting $c = w\sin\left(\psi\right)/\cos\{\psix-\psi\}$ and $b= c \sin\{\psix -\psi\}$ we get that (see Figure~\ref{fig:help_graphic_asymp_limit})  
	$$
	a= w\cos\left(\psi\right) - b  = \frac{ \cos\left(\psi\right)\cos\{\psix-\psi\} - \sin\left(\psi\right)\sin\{\psix-\psi\}  }{\cos\left(\psix-\psi\right)} = \frac{\cos\{\psix\}}{\cos\{\psix-\psi\}}.  
	$$
	From (\ref{def:det-limit-crit-fct}) we obtain that
	\begin{align*}
	\asympcrit\left(w, \psi;x\right) 
	&= \tau  (x) \int_{ \Hr\{ \psix - \psi \} + w \left(\sin \psi, \cos \psi \right)^\top }  K\left(\pvec{z}\right)	\dz   \\
	&= \tau (x)  \int_{ \Hr\{ \psix - \psi \} + b_x\left(w,\psi\right) \pvec{e}_2}  K_1\left(z_1\right) K_2\left(z_2\right)	\dzone \dztwo  \\
	&= \tau (x)  \int_{-1}^1 K_1\left(z_1\right)  \int_{a_{x}(\psi)z_1 + b_{x}(w,\psi)}^{1} K_2\left(z_2\right)	\dztwo  \dzone \\
	&= \tau (x)  \int_{-1}^1 K_1\left(z_1\right)  \left[ \bar K_2 (1 ) - \bar K_2 \{ a_{x}(\psi)z_1 + b_{x}(w,\psi)\} \right]  \dzone \\
	&= -\tau (x)  \int_{-1}^1 K_1\left(z_1\right)   \bar K_2 \{ a_{x}(\psi)z_1 + b_{x}(w,\psi)\}  \dzone.	
	\end{align*}
	Since $ |a_x(\psi)y + b_x(w,\psi)| \geq 1$ if $|w| \geq 1/{|\cos \{ \psi(x)\}|}+1$
	for all $y \in [-1,1]$ and $\psi \in [-\pi/2,\pi/2]$, $\asympcrit\left(w, \psi;x\right)$ vanishes outside a compact set of values of $w$ and $x$.  
	
Turning to (\ref{eq:contrast}), if we show that it holds for individual $x$, then since the supremum in (\ref{eq:contrast}) can be taken over a single compact set, we have that the left-hand side of (\ref{eq:contrast}) is a continuous function in $x$ which is positive for any  $x \in [0,1]$. Hence, the infimum over $x\in[0,1]$ is still positive.
	
To show that (\ref{eq:contrast}) holds for individual $x$, we observe that $-\bar K_2(0)=1$, $-\bar K_2(y)<1$ if $y \not = 0$, and $\bar K_2(y) = 0$ if $|y| \geq 1$ by Assumption~\ref{assumption:kernel}. Since $\asympcrit\left(w, \psi;x\right)$ is continuous in $w$ and $\psi$, as can also be seen from (\ref{def:det-limit-crit-fct}), it is enough to show that $(0,\psi(x))$ is the unique maximizer of $\asympcrit (w, \psi;x )$. But this is immediate from (\ref{eq:asympcritrepres}) and the above properties of $\bar K_2$ and the positivity of $K_1$, since if $(w,\psi) \not= (0, \psi(x))$, there is at most a single value of $z_1$ for which $a_x(\psi)z_1 + b_x(w,\psi)=0$ and hence $-\bar K_2 \{ a_x(\psi)z_1 + b_x(w,\psi) \} =1$.  
\end{proof}

\begin{proof}[{\sl Proof of Proposition~\ref{prop:consistency}}]
The uniform consistency of the estimates $\estw$ and $\estpsi$ is immediate from Proposition~\ref{prop:vdv_consistency:uniform}, as Lemma~\ref{lemma:pointwise_convergence} and Lemma~\ref{lemma:derivatives_asymptotic} provide the necessary assumptions.
\end{proof}

\subsection{Rate of convergence: Proof of Theorem~\ref{theorem:consistency}}

To prove the theorems in the main part we start with a simple linearization. By the Mean Value Theorem, we have that 
\begin{align*}
\nabla \empcrit\left(0,\psix;x \right) & = \nabla \empcrit\left(0,\psix;x \right) - \nabla \empcrit \big(\estw(x),\estpsi(x);x \big) \\
&= - \int_0^1 \nabla \nabla^\top \empcrit\big(t\estw(x),\psix+ t\{ \estpsi(x)-\psix\} ;x\big) \dt\  (\estw (x) ,\estpsi (x) -\psix )^\top ,
\end{align*}
since $\nabla \empcrit\left(\estw(x),\estpsi(x);x \right) =0$. This implies
\begin{align}\label{exp:estim} 
\left(\estw (x) ,\estpsi (x) -\psix\right)^\top =   -   \hat{\pvec{H}}\inv_n (x) \ \nabla \empcrit\left(0,\psix;x \right), 
\end{align}
where
$$
\hat{\pvec{H}}_n (x) =  \int_0^1 \nabla \nabla^\top \empcrit\left(t\estw(x),\psix+ t\{ \estpsi(x)-\psix\} ;x\right) \dt,
$$
and the existence of the inverse of the Hessian matrix $\hat{\pvec{H}}_n (x) $ uniformly in $x \in I$ for large $n$ and with high probability follows from Lemma~\ref{lemma:convergence_hessian} below. 

Similarly as in (\ref{exp:estim}), for the rescaled maximizers of the deterministic contrast function, i.e.,
\begin{align*} 
(w_n(x), \psi_n(x) ) \in  \argmax{   (w,\psi)^\top \in \maxset }  \detcrit(w,\psi,x),
\end{align*}
we have that 
\begin{align}\label{exp:detmax} 
\left(w_n (x),\psi_n (x) -\psix\right)^\top =   -   {\pvec{H}}\inv_n (x)\ \nabla \detcrit\left(0,\psix;x \right), 
\end{align}
where $w_n (x) = \{ \phi_n(x) - \phi(x)\}/h$ and the deterministic Hessian matrix is given by
$$
{\pvec{H}}_n (x) =  \int_0^1 \nabla \nabla^\top \detcrit\left(t\,w_n(x),\psix+ t\{ \psi_n(x)-\psix\};x\right) \dt. 
$$

By means of (\ref{exp:detmax}), we can decompose the representation in (\ref{exp:estim})	into a stochastic part (or score part), a bias part and a remainder term as follows:
\begin{align} \label{defi_repr_taylor}
\left(\estw (x) ,\estpsi (x) -\psix\right)^\top 
&=   \hat{\pvec{H}}\inv_n (x) \scorepartphipsi(0,\psi(x);x)  -  {\pvec{H}}\inv_n (x) \ \biaspartphipsi(x) + \remainderpartphipsi(x),
\end{align}
where 
\begin{align} 
\begin{split}
\scorepartphipsi(w,\psi;x)&=  (\scorepartphi(w,\psi;x),\scorepartpsi(w,\psi;x)  )^\top =  - \{\nabla \empcrit\left(w,\psi;x \right)-\nabla \detcrit\left(w,\psi;x \right) \},\\	
\biaspartphipsi(x) &=  (\biaspartphi(x),\biaspartpsi(x) )^\top =  \nabla \detcrit (0,\psix;x ), \nonumber \\ 
\remainderpartphipsi(x) &= \{ \pvec{H}\inv_n(x) -   \hat{\pvec{H}}\inv_n (x)   \} \biaspartphipsi(x) .
\end{split}
\end{align}

Note that by (\ref{eq:asympcritrepres}) we have $\tau(x)= \asympcrit(0,\psix;x)$ and $\esttau(x)=\empcrit(\estw(x),\estpsi(x);x)$ such that by the Mean Value Theorem as well as adding and subtracting $\detcrit(0,\psix;x)$, one has
\begin{align} \label{eq_taylor_esttau}
\begin{split}
\esttau(x) -\tau(x) & = \empcrit(0,\psix;x) - \asympcrit(0,\psix;x)  \\	
&\quad  + (\estw(x),\estpsi(x)-\psix) \int_0^1 \nabla  \empcrit\left(t\estw(x),\psix + t\{ \estpsi(x)-\psix\} ;x\right) \dt \\	
& = \scoreparttau(0,\psix;x) + \biasparttau(x) + \remainderparttau(x),
\end{split}
\end{align}
where 
\begin{align}
\begin{split}		
\scoreparttau(w,\psi;x) &=\empcrit(w,\psi;x) - \detcrit(w,\psi;x), \\		
\biasparttau(x)	&= \detcrit(0,\psix;x) - \asympcrit(0,\psix;x), \nonumber \\		
\remainderparttau(x) &= (\estw(x),\estpsi(x)-\psix) \int_0^1 \nabla  \empcrit\left(t\estw(x),\psix+ t\{ \estpsi(x)-\psix\};x\right) \dt.		
\end{split}
\end{align}

\subsection{Asymptotic bias}

\begin{lemma} \label{lemma:asympbias}
One has $\nabla \detcrit( w,\psi;x ) = \E\{ \nabla \empcrit(w,\psi;x )\}$ for any $(x,w,\psi)\in \paramset$. Moreover, under Assumptions~\ref{assumption:smoothness} and  \ref{assumption:kernel} we have that 
\[ 
\sup_{x\in \xset} \norm{  \biaspartphi } =O(h^2) + O\{ (nh)\inv\}.
\]
If in addition Assumption~\ref{assumption:extra_kernel} is fulfilled, then  		
\[ 
\sup_{x\in \xset} \norm{  \biaspartpsi } =O(h^2) + O \{ (nh)\inv\}, \quad  \sup_{x\in \xset} \norm{  \biasparttau } =O(h^2) + O\{ (nh)\inv\}.  
\]	
The constants in the remainder terms depend only on the Lipschitz constants of $K,m,\phi$ and $\psi$ as well as on $\norminf{\tau}$. In particular, if Assumption~\ref{assumption:bandwidth} is satisfied, then
\[ 
\sup_{x\in \xset} \ \max\{  \norm{  \biaspartphi(x)}, \norm{ \biaspartpsi(x)  } , \norm{\biasparttau(x)}   \} = O \{ (nh)\inv\}.
\]
\end{lemma}	

The proof of this lemma is provided in {Section~\ref{sec:asymbias} of the Online Supplement}. The bias consists of a smoothing bias of order $O(h^2),$ and a discretization bias of order $O \{ (nh)\inv\}$, which dominates the order of the bias in case of undersmoothing,  that is under Assumption~\ref{assumption:bandwidth}.	

\subsection{Convergence of the Hessian matrix}

\begin{lemma} \label{lemma:convergence_hessian}
Under Assumptions~\ref{assumption:errors}--\ref{assumption:kernel}, we have, as $n \to \infty$,
$$
\norminf{	\hat{\pvec{H}}_n  -\pvec{H} } \Pto 0 \quad \text{and} \quad  \norminf{	\pvec{H}_n -\pvec{H} } \to 0,
$$
where $ \pvec{H}(x) = \diag \{\varphihess(x) , \varpsihess(x) \}$ and  $\varphihess(x)$  as in Theorem~\ref{theorem:main_theorem}, while
$$
\varpsihess(x) =  \tau(x)	K_2^{(1)}(0) \int_{-1}^{1} y^2 K_1(y)\dy.
$$ 
\end{lemma}

The proof of Lemma~\ref{lemma:convergence_hessian} is provided in {Section~\ref{sec:convhessian} of the Online Supplement.} Note that the limit matrix $\pvec{H}$ corresponds to the Hessian matrix of the asymptotic criterion function $\asympcrit$ at the parameters $w = 0$ and $\psi = \psix$. Furthermore, we only need the uniform consistency in Lemma~\ref{lemma:convergence_hessian} to derive the proof of Theorem~\ref{theorem:consistency} and \ref{theorem:main_theorem}. However, for the proof of Theorem~\ref{theorem:main_theorem_uniform} we need the uniform convergence even with an explicit rate of convergence.

\begin{prop} \label{prop:estimates_uniform_O_term}
Under Assumptions~\ref{assumption:errors}--\ref{assumption:kernel}, we have that 
\begin{align*}
	\sup_{x \in \xset}	\norm{ \scorepartphipsi(0,\psix;x) } &= 	O_P\{(\ln  n)^{1/2} /( nh)\}, \\
	\sup_{x \in I } \norm{ \scoreparttau(0,\psix;x) } &= O_P\{(\ln  n)^{1/2} /( nh)\}.
	\end{align*}
\end{prop}

The proof of this proposition is provided in {Section~\ref{sec:asympscore} of the Online Supplement}.\hfill $\Box$

\begin{proof} [{\sl Proof of Theorem~\ref{theorem:consistency}}]
	
	Lemma~\ref{lemma:convergence_hessian} implies that $\hat{\pvec{H}}\inv_n$ is almost surely a stochastically bounded matrix-valued sequence uniformly in $x$.  
	Hence, in combination with Lemma~\ref{lemma:asympbias} deduce
	\begin{align} \label{eq_first_help_conv}
	\sup_{x \in \xset} 	\norm{ \hat{\pvec{H}}\inv_n(x)\  \biaspartphipsi(x)  }	= O_P(h^2)  + O_P\{ (nh )\inv \}
	\end{align}
	as well as 
	\begin{align} \label{eq_sec_help_conv}
	\sup_{x \in \xset} 	\norm{	\remainderpartphipsi(x)	} =o_P(h^2)  + o_P\{ (nh )\inv \} .
	\end{align}
	Next, Slutzky's Lemma and Proposition~\ref{prop:estimates_uniform_O_term} imply
	\begin{align*}
	\sup_{x \in \xset} \norm{ \hat{\pvec{H}}\inv_n(x) \ \scorepartphipsi(0,\psix;x ) } 
	= 	O_P\{ (\ln  n)^{1/2} /(nh)\}.
	\end{align*}
	\normalsize
	From (\ref{defi_repr_taylor})  it immediately follows by combining the three {latter} displays that
	\begin{align*}
	\norminf{ ( \estw, \estpsi- \psi )^\top  } = O_P\{ (\ln  n)^{1/2}/{(nh)} \} + O_P(h^2) = O_P\{ (\ln  n)^{1/2}/(nh) \},
	\end{align*}
	where the last equation is due to Assumption~\ref{assumption:bandwidth}.
	Note that Assumption~\ref{assumption:extra_kernel} is only needed for the bias of the slope in (\ref{eq_first_help_conv}). 

	The uniform rate of convergence for $\esttau$ now follows with a similar argument.
	Indeed, due to representation (\ref{eq_taylor_esttau}) and Lemma~\ref{lemma:asympbias} and Proposition~\ref{prop:estimates_uniform_O_term} it suffices to show that 
	\begin{align} \label{eq_asymp_remainder_tau}
	\norminf{\remainderparttau} = o_P\{ (nh)\inv \}.
	\end{align}
	Since we  verified above that
	$\norminf{\estw,(\estpsi-\psi)}=O_P \{(\ln  n)^{1/2}/(nh) \}$ holds, we only need to show that 
	$$
	 \sup_{x \in I } \norm{f(\estw(x),\estpsi(x);x)}=o_P\{(\ln  n)^{-1/2}\},
$$ 
where we abbreviate 
$$
f(w,\psi;x) =  \int_0^1 \nabla  \empcrit\left(tw,\psix+ t \{ \psi-\psix\} ;x\right) \dt.
$$
By the Mean Value Theorem (differentiating under the integral), Lemma~\ref{lemma:convergence_hessian} and since $\hat{\pvec{H}}_n$ is uniformly bounded,
\begin{align} \label{eq_help_tau_remainder_1}
	\begin{split}
	\sup_{x \in I } \norm{ f(\estw(x),\estpsi(x);x)  - f(0,\psix;x)	} &\leq  \norm{ \hat{\pvec{H}}_n } \sup_{x \in I } |	\big( \estw(x),\estpsi(x) -\psix \big)^\top		|	\\
	&=O_P\{ (\ln  n)^{1/2}/(nh) \} =o_P\{ (\ln  n)^{-1/2}\},	
\end{split}
\end{align}
where the last equation is {justified by} Assumption~\ref{assumption:bandwidth}.
Note that 
$$
f(0,\psix;x)   = \nabla \empcrit(0,\psix;x)   = \scorepartphipsi(0,\psix;x) - \biaspartphipsi(x).  
$$
Hence, due to Lemma~\ref{lemma:asympbias}, Proposition~\ref{prop:estimates_uniform_O_term} and the triangle inequality
\begin{align} \label{eq_help_tau_remainder_2}
\sup_{x \in I }\norm{  f(0,\psix;x)	} 	=O_P \{(\ln  n)^{1/2}/(nh) \} =o_P\{ (\ln  n)^{-1/2}\},		
\end{align}	 
so that $ \sup_{x \in I } \norm{f(\estw(x),\estpsi(x);x)}=o_P\{(\ln  n)^{-1/2}\}$ follows by (\ref{eq_help_tau_remainder_1}) and (\ref{eq_help_tau_remainder_2}) by using the triangle inequality.
\end{proof}

\subsection{Asymptotic normality: Proof of Theorem~\ref{theorem:main_theorem}}

{We first establish the asymptotic normality of the score.}

\begin{lemma} \label{lemma:asymptitc_normality}
Under Assumptions~\ref{assumption:errors}, \ref{assumption:bandwidth} and \ref{assumption:kernel} we have that, for any $x\in I$, as $n \to \infty$,
	\begin{align*}
	nh \left(  \begin{array} {c}
	\scorepartphi(0,\psix;x) \\ \scorepartpsi(0,\psix;x) \\ \scoreparttau\ \ (0,\psix;x)
	\end{array}   \right) \rightsquigarrow \mathcal{N}_3 [ \pvec{0},\sigma^2 \diag ( \varphisco(x), \varpsisco ,  \vartausco )  ],
	\end{align*}	
	where $\varphisco$ is as in Theorem~\ref{theorem:main_theorem} and
	\begin{align*}	
	\varpsisco & = \int_{-1}^1 \Big [ \int_{-1}^{1} \{  K_1(z_1)K_2^{(1)}(z_2)   z_1 -  K_1^{(1)}(z_1)K_2(z_2)   z_2 \}^2\mathrm{d} z_1  \Big] \mathrm{d} z_2, \quad 
	\vartausco  = \int_{ [-1,1]^2}  \{ K(\pvec{z}) \}^2\dz.	
	\end{align*}	
\end{lemma}

{Section~\ref{sec:asympscore} of the Online Supplement} contains the proof of this lemma. 	
Theorem~\ref{theorem:main_theorem} now follows immediately from the following more general theorem, since $ \estw(x) = \{ \estphi(x) - \phi(x)\} / h$ as well as $ w_n(x) = \{ \phi_n(x) - \phi(x)\} / h$.

\begin{satz} \label{theorem:asymp_norm_of_all}
Under Assumptions~\ref{assumption:errors}--\ref{assumption:kernel}, one has, as $n \to \infty$,
$$
nh\,\left(
\begin{array} {c}
\estw(x)-  w_n(x) \\ \estpsi(x)-\psi_n(x) \\ \esttau(x) - \tau_n(x)
\end{array}
\right) \rightsquigarrow \mathcal{N}_3 [ \pvec{0},\sigma^2 \,\Sigma(x)  ],	
$$
where $\Sigma(x)  = \diag  [  \varphisco(x)/\{ \varphihess(x)\}^{2}, \varpsisco /  \{ \varpsihess(x) \}^{2}  ,  \vartausco ]$.
\end{satz}

\begin{proof}[{\sl Proof of Theorem~\ref{theorem:asymp_norm_of_all}}]
From (\ref{exp:estim}) and (\ref{exp:detmax}), we obtain that	
$$
nh\,\left(\estw (x) -  w_n (x) ,\estpsi (x) -\psi_n (x) \right)^\top
=      \hat{\pvec{H}}\inv_n (x) \ nh\scorepartphipsi(0,\psix;x) 
	+ nh\remainderpartphipsi(x).
$$
	Now, by the Mean Value Theorem one has that
	\begin{align*}
	\tau_n(x) -\tau(x) 
	&= 	\detcrit(w_n(x),\psi_n(x);x) - \asympcrit(0,\psix;x) \\
	&=	\biasparttau(x) 
	+ (w_n(x),\psi_n(x)- \psix) \int_0^1 \nabla  \detcrit\left(tw,\psix+ t\{\psi-\psix\};x\right) \dt \\
	&= \biasparttau(x) + R_{n,2}^{\tau}(x).
	\end{align*}		
	With this, (\ref{eq_taylor_esttau}) and the first display of this proof, we have
	\begin{align} \label{eq_asymp_decomp_of_all}
	nh\,\left(
	\begin{array} {c}
	\estw(x)-  w_n(x) \\ \estpsi(x)-\psi_n(x) \\ \esttau(x) - \tau_n(x)
	\end{array}
	\right)
	= &  nh\left( \begin{array} {c}		
	\hat{\pvec{H}}\inv_n (x) \ \scorepartphipsi(0,\psix;x) \\			
	\scoreparttau(0,\psix;x)			
	\end{array}  \right)
	+ nh\left( \begin{array} {c}		
	\remainderpartphipsi(x) \\		
	\remainderparttau(x) - R_{n,2}^{\tau}(x)		
	\end{array}  \right) .		
	\end{align}

On the one hand, the first term is asymptotically normally distributed with covariance matrix $\Sigma$ as in the assumption, since by 
Lemmas~\ref{lemma:convergence_hessian} and \ref{lemma:asymptitc_normality} and Slutzky's Lemma one has, for all $x\in (0,1)$,
$$
\big( \Sigma(x) \big)_{1:2,1:2} =	\sigma^2 \,	\hat{\pvec{H}}^{-2}(x) \diag  \{\varphisco(x),  \varpsisco \},
$$
where $ ( \Sigma(x) )_{1:2,1:2}$ is the matrix given consisting of the first two columns or first two rows of $\Sigma(x)$.
	On the other hand, the second term on the {right-hand} side (\ref{eq_asymp_decomp_of_all}) is $o_P(1)$ by (\ref{eq_sec_help_conv}) and (\ref{eq_asymp_remainder_tau}), whereas one can show $R_{n,2}^{\tau}(x) = o_p \{ (nh)\inv\}$ similar as (\ref{eq_asymp_remainder_tau}). 
	This concludes the proof of Theorem~\ref{theorem:asymp_norm_of_all}.
\end{proof}	

\subsection{Uniform confidence bands: Proof of Theorem~\ref{theorem:main_theorem_uniform}}

For the proof of Theorem~\ref{theorem:main_theorem_uniform} we require rates of convergence of the normalized estimators of the Hessian matrix.
\begin{lemma} \label{lemma:convergence_hessian_uniform_aux_result}
	Under Assumptions~\ref{assumption:errors}--\ref{assumption:kernel} we have that
	\begin{align*}		
	&\norminf{\varphihessest  - \varphihess} = O_P(h) + O_P \{ \sqrt{\ln (n)}/{(nh)}\}, \quad 	
	\norminf{\varphiscoest  - \varphisco} =  O_P\{ \sqrt{\ln (n)}/{(nh)}\}, \\		
	&\norminf{\varpsihessest  - \varpsihess} = O_P(h) + O_P \{ \sqrt{\ln (n)}/{(nh)}\},
	\quad 		
	\norminf{ \pvec{\hat H}_n\inv  - \pvec{H}\inv} = O_P(h) + O_P\{ \sqrt{\ln (n)}/{(nh)} \}.	
	\end{align*}
	Here, $\varpsihessest =   \esttau(x)	K_2^{(1)}(0) \int_{-1}^{1} y^2 K_1(y)\dy.$
\end{lemma}

The proof is provided in {Section~\ref{sec:rate_of_conhessian} of the Online Supplement}. 

Next we extend our notation by incorporating the following definition. The L\'{e}vy-concentration function of a random variable $X$ is given,
for all $\zeta \geq 0$, by
\begin{align*} 
\levyanti(X,\zeta) = \sup_{x \in \R} \Pr (|X-x|\leq \zeta).
\end{align*}

We introduce the normalized score process
\begin{align} \label{eq:score_process_phi}
\begin{split}
\zscorephi(x) 
&=   \sigma\inv    \{\varphisco(x)\}^{-1/2}  \ nh \scorepartphi(0,\psix;x). \\
\end{split}
\end{align}
It easily follows {from Lemma~\ref{lemma:first_deriv_empcrit} in the Online Supplement} that
\begin{align*}
\scorepartphi(0,\psix;x) &=  \partial_w \{ \empcrit (0,\psix;x) - \detcrit (0,\psix;x) \} \\	
&=  (n^2 h^2)^{-1} \sum\limits_{i_1,i_2=1}^n \eionetwo \scalarp{   (\nabla K) \{ h\inv \rotmat_{-\psix} (\px -\xionetwo)\} ,(\sin \psix, \cos \psix)^\top },
\end{align*}
such that 
\begin{align*}
\zscorephi(x) 	
&=  \frac{1}{ nh  \sqrt{\varphisco(x)} \sigma } \sumieiz \eionetwo 
\scalarp{   (\nabla K) \{ h\inv \rotmat_{-\psix}  (x,\phi(x))^\top -\xionetwo)\} , (\sin \psix, \cos \psix )^\top }.
\end{align*}

The score process in (\ref{defi:score_processes_bootstrap}) can be obtained from the latter display by replacing the actual parameters $(\phi,\psi)$ by their estimates $(\estphi,\estpsi)$ and the noise of the observations $\eionetwo/\sigma$ by the sequence $\xi_{\ionetwo}$. Similarly to (\ref{eq:bootprocess_phi}) we set
\begin{align*} 
\orgprocessphi = \norminf{	\zscorephi}.	
\end{align*}

\begin{lemma} \label{lemma:levy-anti-concentration_phi}
Under Assumptions~\ref{assumption:errors}--\ref{assumption:kernel} the following statements are valid.
\begin{itemize}
\item [(i)] 
By possibly enriching the probability space, one has
\begin{equation*} 		
|\orgprocessphi  -	\bootprocessphi|	= O_P\Big\{ \frac{(\ln  n)^{1/2}}{n^{1/2}h} \Big\}.
\end{equation*}
\item [(ii)] 
One has $	\bootprocessphi = O_P\{ (\ln  n)^{1/2}\}$ and for any sequence $\delta_n=o \{(\ln  n)^{-1/2}\}$, one has $\levyanti( \bootprocessphi ,\delta_n ) =o(1)$.
\item [(iii)] 
Moreover,  $q_{1-\alpha}(\bootprocessphi) \cong (\ln \,n)^{1/2}.$
\end{itemize}	
\end{lemma}

The proof is given in {Section~\ref{sec:gauss_approx_phi_psi} of the Online Supplement}. 

\begin{proof}[{\sl Proof of Theorem~\ref{theorem:main_theorem_uniform}}]
Set $ \hat d_n = \hat\sigma_n q_{1-\alpha}(\bootprocessphi)$. Recall that for $\pvec{z}=(z_1,z_2)^\top \in \R^2$ we write $(\pvec{z})_i = z_i$ for $i \in \{ 1,2\}$. From (\ref{defi_repr_taylor}) we obtain that 
\begin{align*}
\Pr \Big[	\sup_{x\in \xset} \Big| \frac{ \varphihessest(x) \{\estphi(x) - \phi(x) \}  }{\{ \varphiscoest(x)\}^{1/2}}  &  \Big| \geq {\hat d_n (1+t_n)}/{n}		\Big] \\	
	& = 	\Pr \Big[	\sup_{x\in \xset} \Big| \frac{ \varphihessest(x)
	}{\{ \varphiscoest(x)\}^{1/2}}\, nh \Big[ \hat{\pvec{H}}\inv_n(x) \{ \scorepartphipsi(0,\psix;x )  +   \biaspartphipsi(x)\} \Big]_1 \Big| \geq   \hat d_n (1+t_n)		\Big] \\
	&\leq \Pr \Big[	\sup_{x\in \xset} \Big|\frac{ \varphihessest(x)  }{\{ \varphiscoest(x)\}^{1/2}} \, nh \Big[ \hat{\pvec{H}}\inv_n(x) \, \scorepartphipsi(0,\psix;x )  \Big]_1\Big| \geq \hat d_n		\Big] \\
	&\quad + \Pr\Big[ \sup_{x\in \xset} \Big|  \frac{ \varphihessest(x)  }{\{ \varphiscoest(x)\}^{1/2}} \Big[ \hat{\pvec{H}}\inv_n(x)\ nh\biaspartphipsi(x)  \Big]_1  \Big| \geq t_n \hat d_n		\Big].
	\end{align*}

	From the last statement of Lemma~\ref{lemma:levy-anti-concentration_phi} and the assumption on $t_n$ we deduce that $t_n \hat d_n	\to \infty$ as $n \to \infty$. 
	Hence, by using Lemmas \ref{lemma:asympbias} and \ref{lemma:convergence_hessian_uniform_aux_result} the second term in the preceding display can be bounded by 
	\begin{align*}	
	\Pr \Big[ \sup_{x\in \xset} \Big|  \frac{ \varphihessest(x)  }{\{ \varphiscoest(x)\}^{1/2}} \Big[ \hat{\pvec{H}}\inv_n(x)\ nh\biaspartphipsi(x)  \Big]_1  \Big| \geq t_n \hat d_n	\Big]
	=o(1).
	\end{align*}

Recalling the notation in (\ref{eq:score_process_phi}), the first term can be estimated by 
\begin{align} \label{eq:help_uniform_cond_band1}
\begin{split}	
\Pr \Big[ 	\sup_{x\in \xset} & \Big|\frac{ \varphihessest(x)  }{\{ \varphiscoest(x)\}^{1/2}}nh \Big[ \hat{\pvec{H}}\inv_n(x) \scorepartphipsi(0,\psix;x )  \Big]_1\Big| \geq \hat d_n	\Big]  \\
	& \leq \Pr \Big[  \sup_{x\in \xset} \abs{ \Big[ \Big[\frac{\varphihessest(x)}{\{ \varphiscoest(x)\}^{1/2}} \hat{\pvec{H}}\inv_n(x) - \frac{ \varphihess(x)}{ \{ \varphisco(x)\}^{1/2}}\pvec{H}\inv(x) \Big]  \ nh\,\scorepartphipsi(0,\psix;x ) \Big]_1   } \geq  \ln (n)^{-1} \ \Big] \\
	& \quad +   \Pr  \Big[	\orgprocessphi \geq  \{ \hat d_n	 - \ln (n)^{-1} \} /\sigma	\Big].
	\end{split}	
	\end{align}
	
Lemma~\ref{lemma:convergence_hessian_uniform_aux_result} implies that
\begin{align*}
& \norminf{ ({\varphihessest}/{\sqrt{\varphiscoest}} ) \hat{\pvec{H}}\inv_n - (\varphihess/\sqrt{ \varphisco})  \pvec{H}\inv} = O_P(h) + 
O_P\{ (\ln  n)^{1/2}/{(nh)} \}, 
\end{align*}
and together with Proposition~\ref{prop:estimates_uniform_O_term} the first term in (\ref{eq:help_uniform_cond_band1}) can be bounded by 
\[
O_P \{ (\ln  n)^{3/2}h\} + O_P\{ (\ln  n)^{2}/ nh\} = o_P(1),
\]
due to Assumption~\ref{assumption:bandwidth}. As for the second term, plugging in $ \hat d_n = \hat\sigma_n q_{1-\alpha}(\bootprocessphi)$ gives
	\begin{align*}
\Pr \Big[	\orgprocessphi \geq & \{ \hat\sigma_n   q_{1-\alpha}(\bootprocessphi) 	 - \ln (n)^{-1}\}/\sigma	\Big] \\
	& \leq \Pr \{	\orgprocessphi  -	\bootprocessphi  \geq  \ln (n)^{-1}/\sigma	\} + \Pr [\bootprocessphi \geq \{ \hat\sigma_n   q_{1-\alpha}(\bootprocessphi)	 - 2\,\ln (n)^{-1}\}/\sigma ]\\ 		
	\leq & \,o(1) +  \Pr \Big\{ \bootprocessphi \geq q_{1-\alpha}(\bootprocessphi)\hat\sigma_n / \sigma	\Big\}  
	+ \levyanti\big(	\bootprocessphi ,2\,\ln (n)^{-1}/\sigma\big)
	\end{align*}
by using the first part of Lemma~\ref{lemma:levy-anti-concentration_phi} together with the choice of $h$, and the definition of the L\'{e}vy-concentration function. The last term in this display is $o(1)$ by using the second part of Lemma~\ref{lemma:levy-anti-concentration_phi}. Finally, observe that 
\begin{align*}
\Pr \Big\{ \bootprocessphi \geq q_{1-\alpha}(\bootprocessphi)\hat\sigma_n / \sigma	\Big\}  
	&\leq \Pr \Big\{		\bootprocessphi \geq   q_{1-\alpha}(\bootprocessphi) 	\Big\} + \Pr \Big\{	\big|	\bootprocessphi - q_{1-\alpha}(\bootprocessphi)\big| \leq   \big|\hat\sigma_n / \sigma- 1\big|q_{1-\alpha}(\bootprocessphi) 	\Big\} \\	
	&\leq \Pr \big\{		\bootprocessphi \geq   q_{1-\alpha}(\bootprocessphi) 	\big\} + \levyanti\big(	\bootprocessphi ,s_n q_{1-\alpha}(\bootprocessphi)\,\big) + \Pr  ( |\hat\sigma_n / \sigma- 1 | \geq s_n  )\\
	&= \Pr \Big\{\bootprocessphi \geq   q_{1-\alpha}(\bootprocessphi) 	\Big\} +o(1),
\end{align*}
where $s_n$ is as in the assumption of the theorem and where we used the second statement of Lemma~\ref{lemma:levy-anti-concentration_phi} together with the third statement of Lemma~\ref{lemma:levy-anti-concentration_phi} to obtain $s_n q_{1-\alpha}(\bootprocessphi)  =o \{ (\ln  n)^{-1/2}\}$.  

Summarizing, we obtain that 
\[
\Pr \Bigg [	\sup_{x\in \xset} \Big| \frac{ \varphihessest(x) \{ \estphi(x) - \phi(x) \}  }{\{ \varphiscoest(x)\}^{1/2}}   \Big| \geq {\hat d_n (1+t_n)}/{n} \Bigg ]
\leq \Pr \big\{		\bootprocessphi \geq   q_{1-\alpha}(\bootprocessphi) 	\big\} + o(1) = \alpha + o(1),
\]
which corresponds to (\ref{eq:asympconvband}). Eventually, the third statement of Lemma~\ref{lemma:levy-anti-concentration_phi} concludes the proof.
\end{proof}

\subsection{Uniform confidence bands for jump-slope and jump-height curve} \label{sec:sketch_uniform_cb_psi_tau}

For independent standard normally distributed random variables $\xi_{1,1}, \ldots, \xi_{n,n},$ independent of $\yionetwo,$  consider in the spirit of (\ref{defi:score_processes_bootstrap})  the processes
\begin{align*}
\zbootpsi( x ) &= \frac{1}{ nh^2 \big(\varpsisco\big)^{1/2} } \\
&\quad  \times 	\sumieiz \xi_{\ionetwo} 
\scalarp{   (\nabla K) [ h\inv \rotmat_{-\estpsi(x)} \{(x,\estphi(x))^\top -\xionetwo \} ] , [\rotmat_{3\pi /2 -\estpsi(x)} \{(x,\estphi(x))-\xionetwo \} ]^\top}, \\	
\zboottau( x ) &= \frac{1}{ nh \big(\vartausco\big)^{1/2} }  	\sumieiz \xi_{\ionetwo} 
K [ h\inv \rotmat_{-\estpsi(x)} \{ (x,\estphi(x))^\top -\xionetwo \} ],
\end{align*}
where $\varpsisco$ and $\vartausco$ are given in Lemma~\ref{lemma:asymptitc_normality}, as well as the maxima of the processes
$$
\bootprocesspsi = \sup_{x \in I} \big|\zbootpsi( x )\big| \quad \mbox{and} \quad  \bootprocesstau = \sup_{x \in I} \big|\zboottau( x )\big|.
$$

The following theorem can be used to construct uniform confidence bands for $\psi$ or $\tau$.

\begin{satz} \label{theorem_uniform_cb_psi_tau}	
Consider model \eqref{model} under Assumptions~\ref{assumption:errors}--\ref{assumption:extra_kernel}, and assume that $\hat \sigma_n$ is an estimator for $\sigma$ which satisfies $\Pr (|{\hat\sigma_n}/{\sigma}  -1 |\geq s_n) = o(1)$ for some sequence $s_n =o \{ \ln (n)^{-1}\}$. Then for $\alpha\in(0,1)$, one has $q_{1-\alpha}(\bootprocesspsi) = O\{(\ln  n )^{1/2}\}$ and for any sequence $t_n=o(1)$ such that $t_n \sqrt{\ln (n)} \to \infty$, one also has
$$
\liminf_n \Pr \left[	\sup_{x\in \xset}  \Big|  { \varpsihessest(x) \{\estpsi(x) - \psix \}  }/{\sqrt{\varpsisco}}   \Big| \leq {(1+t_n) \hat \sigma_n q_{1-\alpha}(\bootprocesspsi)}/{(nh)}		\right] \geq 1-\alpha, 
$$	
as well as $q_{1-\alpha}(\bootprocesstau) = O \{(\ln  n )^{1/2}\}$ and
$$
\liminf_n \Pr \left(	\sup_{x\in \xset}  \Big|  { \{\esttau(x) - \tau(x)\}  }/{\sqrt{\vartausco}}   \Big| \leq  {(1+t_n) \hat \sigma_n q_{1-\alpha}(\bootprocesstau)}/{(nh)} \right] \geq 1-\alpha. 
$$	
\end{satz}

\begin{proof}[{\sl Sketch of proof of Theorem~\ref{theorem_uniform_cb_psi_tau}}]
	For the jump-slope and the jump-height, we introduce in the spirit of (\ref{eq:score_process_phi}) the following processes
	\begin{align*} 
	\begin{split}
	\zscorepsi(x) 
	&=   \sigma\inv    (\varpsisco)^{-1/2}   \ nh \scorepartpsi(0,\psix;x)  \quad \mbox{and}  \quad
	\zscoretau(x) 
	=   \sigma\inv    (\vartausco)^{-1/2}   \ nh \scoreparttau(0,\psix;x)  , 
	\end{split}
	\end{align*}
	and their suprema $\orgprocesspsi = \norminf{\zscorepsi}$ and $\orgprocesstau = \norminf{\zscoretau }$. Following the lines of proof of Lemma~\ref{lemma:levy-anti-concentration_phi}, one derives the following result.
	
\begin{lemma} \label{lemma:levy-anti-concentration_psi}
Under Assumptions~\ref{assumption:errors}--\ref{assumption:kernel}, and by possibly enriching the probability space, one has
\begin{equation*} 
\norminf{\zscorepsi  -  \zbootpsi} = O_P\Big\{ \frac{(\ln  n)^{1/2}}{n^{1/2}h} \Big\} \quad  \mbox{and} \quad	\norminf{\zscoretau  - \zboottau   }	= O_P\Big\{ \frac{(\ln  n)^{1/2}}{n^{1/2}h} \Big\}.
\end{equation*}
Moreover, $\bootprocesspsi = O_P \{ (\ln  n)^{1/2}\}$ resp.\ $\bootprocesstau = O_P\{(\ln  n)^{1/2}\} $ and for any sequence $\delta_n=o\{(\ln  n)^{-1/2}\}$, we have that
$$
\levyanti(  	\bootprocesspsi ,\delta_n ) =o(1),\quad \levyanti(  	\bootprocesstau ,\delta_n ) =o(1).
$$
Eventually, one has $q_{1-\alpha}(\bootprocesspsi) \cong (\ln \,n)^{1/2}$ resp.\ $q_{1-\alpha}(\bootprocesstau) \cong (\ln \,n)^{1/2}.$
\end{lemma}

With this, it is straightforward to obtain the proof following the lines of the proof of Theorem~\ref{theorem:main_theorem_uniform}.
\end{proof}

\section*{Acknowledgments}

The authors are thankful to Axel Munk for helpful discussions in the early stage of the project, as well as to the Editor-in-Chief of the Journal of Multivariate Analysis, Christian Genest, an Associate Editor and two anonymous reviewers for helpful comments. Financial support of the Deutsche Forschungsgemeinschaft, grant Ho 3260/5-1, is gratefully acknowledged. 

%\section*{References}



\section*{Supplementary material}

# Supplement to "Asymptotic confidence sets for the jump curve in bivariate regression problems"

Viktor Bengs, Matthias Eulert, Hajo Holzmann[*]

*Fb. 12 - Mathematik und Informatik, Philipps-Universität Marburg, Hans-Meerwein-Straße 6, 35032 Marburg, Germany*

**Abstract**

This supplementary material contains detailed proofs of the results in Bengs et al. [1].

*Keywords:* Image Processing, Jump Detection, M-Estimation, Rotated Difference Kernel Estimator.

## A. Supplement: Proofs for uniform consistency and asymptotic normality

*Additional Notation.* In the following we shall denote by $[\nabla K(\mathbf{z}) \mid \nabla K(\mathbf{z})]$ the $2 \times 2$ - matrix with columns equal to $\nabla K(\mathbf{z})$.

### A.1. Uniform consistency of M-estimates: proof of Proposition 2

*Proof of Proposition 2.* From (17) we have that $\sup_{x \in I} |\hat{f}_n(\theta_0(x), x) - f(\theta_0(x), x)| = o_P(1)$. By the property of the estimator $\hat{\theta}_n(x)$ it follows that $\sup_{x \in I} (\hat{f}_n(\theta_0(x), x) - \hat{f}_n(\hat{\theta}_n(x), x)) \leq o_P(1)$. With this, $\sup_{x \in I} (f(\theta_0(x), x) - \hat{f}_n(\hat{\theta}_n(x), x)) \leq o_P(1)$. Next, using (17) again yields

$$\sup_{x \in I} (f(\theta_0(x), x)) - f(\hat{\theta}_n(x), x) = \sup_{x \in I} (f(\theta_0(x), x) - \hat{f}_n(\hat{\theta}_n(x), x) + \hat{f}_n(\hat{\theta}_n(x), x) - f(\hat{\theta}_n(x), x))$$

$$\leq o_P(1) + \sup_{(x, \theta) \in \Theta_n} |\hat{f}_n(\theta, x) - f(\theta, x)| = o_P(1).$$

Given $\varepsilon > 0$ choose $\eta > 0$ as the left side of (18). Then

$$P(||\hat{\theta}_n(x) - \theta_0(x)||_\infty > \varepsilon) = P(\exists\, x \in I : |\hat{\theta}_n(x) - \theta_0(x)| \geq \varepsilon)$$
$$\leq P(\exists\, x \in I : f(\theta_0(x), x) - f(\hat{\theta}_n(x), x) \geq \eta)$$
$$= P\big(\sup_{x \in I} (f(\theta_0(x), x) - f(\hat{\theta}_n(x), x)) \geq \eta\big) = o(1).$$

□

### A.2. Order of the asymptotic bias: proof of Lemma 3

Recall from (15) that the rescaled contrast function is given by

$$\hat{\mathbb{M}}_n(w, \psi; x) = (nh)^{-2} \sum_{i_1, i_2 = 1}^{n} Y_{i_1, i_2} K\left(h^{-1} \mathcal{R}_{-\psi}(\mathbf{p}(x) + wh\mathbf{e_2} - \mathbf{x_{i_1, i_2}})\right), \qquad \mathbf{p}(x) = (x, \phi(x))^T. \tag{1}$$

---

Let
$$T(\mathbf{z};w,\psi) = \mathcal{R}_{3\pi/2-\psi}(\mathbf{p}(x)+wh\mathbf{e_2}-\mathbf{z}),$$
so that
$$T(\mathbf{p}(x)+wh\mathbf{e_2}-h\mathcal{R}_\psi\mathbf{z};w,\psi) = h\mathcal{R}_{3\pi/2}\mathbf{z} = h(z_2,-z_1)^T. \tag{2}$$

**Lemma 1.** *It holds that*
$$\partial_w\hat{\mathbb{M}}_n(w,\psi;x) = (n^2h^2)^{-1}\sum_{i_1,i_2=1}^n Y_{i_1,i_2}\langle (\nabla K)\left(h^{-1}\mathcal{R}_{-\psi}(\mathbf{p}(x)+wh\mathbf{e_2}-\mathbf{x_{i_1,i_2}})\right),(\sin(\psi),\cos(\psi))^T\rangle,$$
$$\partial_\psi\hat{\mathbb{M}}_n(w,\psi;x) = (n^2h^3)^{-1}\sum_{i_1,i_2=1}^n Y_{i_1,i_2}\langle (\nabla K)\left(h^{-1}\mathcal{R}_{-\psi}(\mathbf{p}(x)+wh\mathbf{e_2}-\mathbf{x_{i_1,i_2}})\right),T(\mathbf{x_{i_1,i_2}};w,\psi))\rangle.$$

*Proof of Lemma 1.* We have that
$$\partial_w h^{-1}\mathcal{R}_{-\psi}(\mathbf{p}(x)+wh\mathbf{e_2}-\mathbf{x_{i_1,i_2}}) = (\sin(\psi),\cos(\psi))^T,$$
$$\partial_\psi h^{-1}\mathcal{R}_{-\psi}(\mathbf{p}(x)+wh\mathbf{e_2}-\mathbf{x_{i_1,i_2}}) = h^{-1}\mathcal{R}_{3/2\pi-\psi}(\mathbf{p}(x)+wh\mathbf{e_2}-\mathbf{x_{i_1,i_2}}) = h^{-1}T(\mathbf{x_{i_1,i_2}};w,\psi).$$

The statement follows together with (1) by using the chain rule. □

**Lemma 2.** *Under the Assumption 4 it holds that*

(i) $\int_{[-1,1]^2} K^{(0,1)}(z_1,z_2)z_1\,\mathrm{d}z_1\mathrm{d}z_2 = \int_{[-1,1]^2} K^{(1,0)}(z_1,z_2)z_2\,\mathrm{d}z_1\mathrm{d}z_2 = 0;$

(ii) $\int_{[-1,1]\times[0,1]} K^{(1,0)}(z_1,z_2)\,z_1\,\mathrm{d}z_1\mathrm{d}z_2 = \int_{[-1,1]\times[0,1]} K^{(0,1)}(z_1,z_2)\,z_2\,\mathrm{d}z_1\mathrm{d}z_2 = -1;$

(iii) $\int_{[-1,1]\times[0,1]} K^{(0,1)}(z_1,z_2)\,z_1\,\mathrm{d}z_1\mathrm{d}z_2 = \int_{[-1,1]\times[0,1]} K^{(1,0)}(z_1,z_2)\,z_2\,\mathrm{d}z_1\mathrm{d}z_2 = 0;$

*If in addition Assumption 5 is fulfilled, then*

(iv) $\int_{[-1,1]^2} K^{(1,0)}(z_1,z_2)\,z_1\,z_2\,\mathrm{d}z_1\mathrm{d}z_2 = \int_{[-1,1]^2} K^{(0,1)}(z_1,z_2)\,z_1\,z_2\,\mathrm{d}z_1\mathrm{d}z_2 = 0;$

(v) $\int_{[-1,1]^2} K^{(1,0)}(z_1,z_2)\,z_2^2\,\mathrm{d}z_1\mathrm{d}z_2 = \int_{[-1,1]^2} K^{(0,1)}(z_1,z_2)\,z_1^2\,\mathrm{d}z_1\mathrm{d}z_2 = 0;$

(vi) $\int_{[-1,1]\times[0,1]} K^{(1,0)}(z_1,z_2)\,z_1\,z_2\,\mathrm{d}z_1\mathrm{d}z_2 = \int_{[-1,1]\times[0,1]} K^{(0,1)}(z_1,z_2)\,z_1\,z_2\,\mathrm{d}z_1\mathrm{d}z_2 = 0;$

(vii) $\int_{[-1,1]\times[0,1]} K^{(1,0)}(z_1,z_2)\,z_2^2\,\mathrm{d}z_1\mathrm{d}z_2 = \int_{[-1,1]\times[0,1]} K^{(0,1)}(z_1,z_2)\,z_1^2\,\mathrm{d}z_1\mathrm{d}z_2 = 0.$

*Proof of Lemma 2.* Some equalities follow by symmetry and normalization of $K_1$ and $K_2$, others require the boundary properties. For example, by integration by parts
$$\int_{-1}^1 K_1^{(1)}(x)\,x\,\mathrm{d}x = -\int_{-1}^1 K(x)\,\mathrm{d}x = -1$$
since $K_1(1) = K_1(-1) = 0$ and $\int K = 1$ by assumption. In addition, $\int_0^1 K_2(x) = 1$ which yield the first statement in (ii). □

*Proof of Lemma 3.* By means of Lemma 1 and the definition of $\mathbb{M}_n$, it is apparent for any $(x,w,\psi) \in \Theta_n$ that
$$\nabla\mathbb{M}_n(w,\psi;x) = \mathrm{E}\left(\nabla\hat{\mathbb{M}}_n(w,\psi;x)\right).$$

For sake of brevity we write $\mathbf{v}_{\psi(x)} = (\sin\psi(x),\cos\psi(x))^T$. Hence, we discuss the expected values of the partial derivatives of $\hat{\mathbb{M}}_n(0,\psi(x);x)$, that is $\mathcal{B}_n^\phi$ and $\mathcal{B}_n^\psi$ respectively.

*The bias of $\hat{\phi}_n$*

Using Lemma 1 yields

$$\mathbf{B}_n^\phi(x) = \partial_w \mathbb{M}_n(0, \psi(x); x) = (nh)^{-2} \sum_{i_1,i_2=1}^n m(\mathbf{x}_{\mathbf{i_1},\mathbf{i_2}}) \langle (\nabla K) \left(h^{-1} \mathcal{R}_{-\psi(x)}(\mathbf{p}(x) - \mathbf{x}_{\mathbf{i_1},\mathbf{i_2}})\right), \mathbf{v}_{\psi(x)} \rangle$$

$$+ (nh)^{-2} \sum_{i_1,i_2=1}^n j_\tau(\mathbf{x}_{\mathbf{i_1},\mathbf{i_2}}) \langle (\nabla K) \left(h^{-1} \mathcal{R}_{-\psi(x)}(\mathbf{p}(x) - \mathbf{x}_{\mathbf{i_1},\mathbf{i_2}})\right), \mathbf{v}_{\psi(x)} \rangle$$

$$=: S_n + J_n.$$

Note that we omitted the dependency on $x$ for $S_n$ resp. $J_n$. From Lemma 10 one obtains

$$S_n = \int_{[-1,1]^2} \left(m(\mathbf{p}(x) - h\mathcal{R}_{\psi(x)}\mathbf{z}) - m(\mathbf{p}(x))\right) \langle \nabla K(\mathbf{z}), \mathbf{v}_{\psi(x)} \rangle \, d\mathbf{z}$$

$$+ \int_{[-1,1]^2} m(\mathbf{p}(x)) \langle \nabla K(\mathbf{z}), \mathbf{v}_{\psi(x)} \rangle \, d\mathbf{z} + O((nh)^{-1})$$

by taking $g_1(\mathbf{z}) = m(\mathbf{z})$, $f(\mathbf{z}) = [\nabla K(\mathbf{z}) \mid \nabla K(\mathbf{z})]$, $g_2 \equiv \mathbf{e}_1$, $g_3 \equiv \mathbf{v}_{\psi(x)}$, such that $j = 1$, $r_1 = 2$, $r_2 = r_3 = 0$, where we recall that $[\nabla K(\mathbf{z}) \mid \nabla K(\mathbf{z})]$ denotes the $2 \times 2$ - matrix with columns equal to $\nabla K(\mathbf{z})$. As the components $\nabla K(\mathbf{z})$ are odd, the second term is zero. Concerning the first term, obtain by Taylor expansion

$$\int_{[-1,1]^2} \left(m(\mathbf{p}(x) - h\mathcal{R}_{\psi(x)}\mathbf{z}) - m(\mathbf{p}(x))\right) \langle \nabla K(\mathbf{z}), \mathbf{v}_{\psi(x)} \rangle \, d\mathbf{z}$$

$$= h m^{(1,0)}(\mathbf{p}(x)) \int_{[-1,1]^2} (\mathcal{R}_{\psi(x)}\mathbf{z})_1 \langle \nabla K(\mathbf{z}), \mathbf{v}_{\psi(x)} \rangle \, d\mathbf{z} + h m^{(0,1)}(\mathbf{p}(x)) \int_{[-1,1]^2} (\mathcal{R}_{\psi(x)}\mathbf{z})_2 \langle \nabla K(\mathbf{z}), \mathbf{v}_{\psi(x)} \rangle \, d\mathbf{z}$$

$$+ O(h^2).$$

Both integrals on the right-hand side are zero due to Lemma 5, (vii) and Lemma 2, (i). Thus, $S_n = O(h^2) + O((nh)^{-1})$. Similarly, with the same functions $f, g_2, g_3$ as above for Lemma 11 obtain

$$J_n = \int_{h^{-1}\mathcal{R}_{-\psi(x)}(\mathbf{p}(x) - [0,1]^2 \backslash \text{epi}(\phi))} \left(\tau\big(x - (h\mathcal{R}_{\psi(x)}\mathbf{z})_1\big) - \tau(x)\right) \langle \nabla K(\mathbf{z}), \mathbf{v}_{\psi(x)} \rangle \, d\mathbf{z}$$

$$+ \int_{h^{-1}\mathcal{R}_{-\psi(x)}(\mathbf{p}(x) - [0,1]^2 \backslash \text{epi}(\phi))} \tau(x) \langle \nabla K(\mathbf{z}), \mathbf{v}_{\psi(x)} \rangle \, d\mathbf{z} + O((nh)^{-1}).$$

Without loss of generality let $h$ be so small ($n$ large enough) such that

$$[-1,1] \times [0,1] \subset h^{-1}\mathcal{R}_{-\psi(x)}\left(\mathbf{p}(x) - [0,1]^2 \backslash \text{epi}(\phi)\right). \tag{3}$$

Note that $\int_{[-1,1] \times [0,1]} \langle \nabla K(\mathbf{z}), \mathbf{v}_{\psi(x)} \rangle \, d\mathbf{z} = 0$, (e.g. Lemma 5 (xi)) such that the second term is zero in the latter display. For the first term one has by Taylor expansion

$$\int_{[-1,1] \times [0,1]} \left(\tau\big(x - (h\mathcal{R}_{\psi(x)}\mathbf{z})_1\big) - \tau(x)\right) \langle \nabla K(\mathbf{z}), \mathbf{v}_{\psi(x)} \rangle \, d\mathbf{z}$$

$$= h \tau^{(1)}(x) \int_{[-1,1] \times [0,1]} (\mathcal{R}_{\psi(x)}\mathbf{z})_1 \langle \nabla K(\mathbf{z}), \mathbf{v}_{\psi(x)} \rangle \, d\mathbf{z} + O(h^2),$$

where the integral can be further simplified by Lemma 2, (ii) and (iii) to

$$\int_{[-1,1]\times[0,1]} (\mathcal{R}_{\psi(x)}\mathbf{z})_1 \langle \nabla K(\mathbf{z}), \mathbf{v}_{\psi(x)} \rangle \, d\mathbf{z}$$
$$= \int_{[-1,1]\times[0,1]} \big(\cos(\psi(x))z_1 + \sin(-\psi(x))z_2\big)\big(\sin(\psi(x))K^{(1,0)}(z_1,z_2) + \cos(\psi(x))K^{(0,1)}(z_1,z_2)\big) \, dz_1 dz_2$$
$$= -\cos(\psi(x))\sin(\psi(x)) - \cos(\psi(x))\sin(-\psi(x)) = 0.$$

In summary, $J_n = O(h^2) + O((nh)^{-1})$.

*The bias of $\hat{\psi}_n$*

Using Lemma 1 yields

$$\mathcal{B}_n^{\psi}(x) = \partial_{\psi}\mathbb{M}_n(0,\psi(x);x) = (nh^2)^{-3} \sum_{i_1,i_2=1}^n m(\mathbf{x}_{\mathbf{i}_1,\mathbf{i}_2}) \langle (\nabla K)(h^{-1}\mathcal{R}_{-\psi}(\mathbf{p}(x) - \mathbf{x}_{\mathbf{i}_1,\mathbf{i}_2})), T(\mathbf{x}_{\mathbf{i}_1,\mathbf{i}_2}; 0, \psi(x))) \rangle$$
$$+ (nh^2)^{-3} \sum_{i_1,i_2=1}^n j_{\tau}(\mathbf{x}_{\mathbf{i}_1,\mathbf{i}_2}) \langle (\nabla K)(h^{-1}\mathcal{R}_{-\psi}(\mathbf{p}(x) - \mathbf{x}_{\mathbf{i}_1,\mathbf{i}_2})), T(\mathbf{x}_{\mathbf{i}_1,\mathbf{i}_2}; 0, \psi(x))) \rangle$$
$$=: S_n + J_n.$$

Lemma 10 gives us together with (2) that uniformly over $x$

$$S_n = \int_{[-1,1]^2} \big(m(\mathbf{p}(x) - h\mathcal{R}_{\psi(x)}(z_1,z_2)^T) - m(\mathbf{p}(x))\big) \langle \nabla K(z_1,z_2), (z_2,-z_1)^T \rangle \, dz_1 dz_2$$
$$+ \int_{[-1,1]^2} m(\mathbf{p}(x)) \langle \nabla K(z_1,z_2), (z_2,-z_1)^T \rangle \, dz_1 dz_2 + O((nh)^{-1}),$$

where we used $g_1, g_2$ and $f$ as before and $g_3(\mathbf{z}) = T(\mathbf{z}; 0, \psi(x)))$, so that $j = 1$, $r_1 = 3$, $r_2 = 0$ and $r_3 = 1$. The second term on the right-hand side of the latter display is zero since the functions $x \mapsto K_1^{(1)}(x)$ and $x \mapsto K_1(x)x$ are odd. By Taylor expansion in the first term

$$\int_{[-1,1]^2} \big(m(\mathbf{p}(x) - h\mathcal{R}_{\psi(x)}(z_1,z_2)^T) - m(\mathbf{p}(x))\big) \langle \nabla K(z_1,z_2), (z_2,-z_1)^T \rangle \, dz_1 dz_2$$
$$= hm^{(1,0)}(\mathbf{p}(x)) \int_{[-1,1]^2} \big(\mathcal{R}_{\psi(x)}(z_1,z_2)^T\big)_1 \langle \nabla K(z_1,z_2), (z_2,-z_1)^T \rangle \, dz_1 dz_2$$
$$+ hm^{(0,1)}(\mathbf{p}(x)) \int_{[-1,1]^2} \big(\mathcal{R}_{\psi(x)}(z_1,z_2)^T\big)_2 \langle \nabla K(z_1,z_2), (z_2,-z_1)^T \rangle \, dz_1 dz_2 + O(h^2).$$

Lemma 2, (iv) and (v) imply that both integrals on the right-hand side are zero. Therefore, one has $S_n = O(h^2) + O((nh)^{-1})$. We still assume that $h$ is so small respectively $n$ is so large that (3) holds. Thus, application of Lemma 11 with the same functions $g_2, g_2$ and $f$ as just yields

$$J_n = \int_{[-1,1]\times[0,1]} \big(\tau(x - (h\mathcal{R}_{\psi(x)}\mathbf{z})_1) - \tau(x)\big) \langle \nabla K(z_1,z_2), (z_2,-z_1)^T \rangle \, dz_1 dz_2$$
$$+ \int_{[-1,1]\times[0,1]} \tau(x) \langle \nabla K(z_1,z_2), (z_2,-z_1)^T \rangle \, dz_1 dz_2 + O((nh)^{-1}).$$

4

Since $\int_{[-1,1]\times[0,1]} \langle \nabla K(z_1,z_2), (z_2,-z_1)\rangle \,\mathrm{d}z_1\mathrm{d}z_2 = 0$, by Lemma 2, (iii), we only discuss the first term on the right-hand side of the latter display. Similarly as before by a Taylor expansion,

$$\int_{[-1,1]\times[0,1]} \left(\tau\big(x - (h\mathcal{R}_{\psi(x)}\mathbf{z})_1\big) - \tau(x)\right) \langle \nabla K(z_1,z_2), (z_2,-z_1)^T\rangle \,\mathrm{d}z_1\mathrm{d}z_2$$
$$= h\tau^{(1)}(x) \int_{[-1,1]\times[0,1]} (\mathcal{R}_{\psi(x)}\mathbf{z})_1 \langle \nabla K(z_1,z_2), (z_2,-z_1)^T\rangle \,\mathrm{d}z_1\mathrm{d}z_2 + O(h^2),$$

where by Lemma 2, (vi) and (vii), the integral is equivalent to

$$\int_{[-1,1]\times[0,1]} (\mathcal{R}_{\psi(x)}\mathbf{z})_1 \langle \nabla K(z_1,z_2), (z_2,-z_1)^T\rangle \,\mathrm{d}z_1\mathrm{d}z_2$$
$$= \int_{[-1,1]\times[0,1]} \big(\cos(\psi(x))z_1 + \sin(-\psi(x))z_2\big) \big(K^{(1,0)}(z_1,z_2)z_2 - K^{(0,1)}(z_1,z_2)z_1\big) \,\mathrm{d}z_1\mathrm{d}z_2 = 0.$$

Consequently, $J_n = O(h^2) + O((nh)^{-1})$.

*The bias of $\hat{\tau}_n$*

We start with a decomposition of $\mathbb{M}_n(0, \psi(x); x)$

$$\mathbb{M}_n(0, \psi(x); x) = (nh)^{-2} \sum_{i_1,i_2=1}^{n} m\big(\mathbf{x}_{\mathbf{i_1},\mathbf{i_2}}\big) K\big(h^{-1}\mathcal{R}_{-\psi}\big(\mathbf{p}(x) - \mathbf{x}_{\mathbf{i_1},\mathbf{i_2}}\big)\big)$$
$$+ (nh)^{-2} \sum_{i_1,i_2=1}^{n} j_\tau\big(\mathbf{x}_{\mathbf{i_1},\mathbf{i_2}}\big) K\big(h^{-1}\mathcal{R}_{-\psi}\big(\mathbf{p}(x) - \mathbf{x}_{\mathbf{i_1},\mathbf{i_2}}\big)\big) = S_n + J_n.$$

Note that we need a stricter result as Lemma 1, since we now consider special values of $w$ and $\psi$. Using Lemma 10 with

$$g_1(\mathbf{z}) = m(\mathbf{z}), \quad f(\mathbf{z}) = 4^{-1} \begin{pmatrix} K(\mathbf{z}) & K(\mathbf{z}) \\ K(\mathbf{z}) & K(\mathbf{z}) \end{pmatrix}, \quad g_2(\mathbf{z}) = g_3(\mathbf{z}) = (1,1)^T,$$

such that $j = 1$, $r_1 = 2$, $r_2 = r_3 = 0$, we obtain that

$$S_n = \int_{[-1,1]^2} \big(m(\mathbf{p}(x) - h\mathcal{R}_{\psi(x)}\mathbf{z}) - m(\mathbf{p}(x))\big) K(\mathbf{z})\,\mathrm{d}\mathbf{z} + \big(m(\mathbf{p}(x))\big) \int_{[-1,1]^2} K(\mathbf{z})\,\mathrm{d}\mathbf{z} + O\big((nh)^{-1}\big)$$
$$= \int_{[-1,1]^2} \big(m(\mathbf{p}(x) - h\mathcal{R}_{\psi(x)}\mathbf{z}) - m(\mathbf{p}(x))\big) K(\mathbf{z})\,\mathrm{d}\mathbf{z} + O\big((nh)^{-1}\big),$$

where the last equality holds as $K_2$ is odd. By Taylor expansion and since $K_2$ is odd resp. due to Assumption 5 we deduce for the first term on the right hand side of the latter display that it is $O(h^2)$ and therefore $S_n = O(h^2) + O\big((nh)^{-1}\big)$. Similarly, using Lemma 11 with $g_2, g_3$ and $f$ as just,

$$J_n = \int_{[-1,1]\times[0,1]} \big(\tau\big(x - (h\mathcal{R}_{\psi(x)}\mathbf{z})_1\big) - \tau(x)\big) K(z_1,z_2)\,\mathrm{d}z_1\mathrm{d}z_2 + \tau(x) \int_{[-1,1]\times[0,1]} K(z_1,z_2)\,\mathrm{d}z_1\mathrm{d}z_2 + O((nh)^{-1})$$
$$= \int_{[-1,1]\times[0,1]} \big(\tau\big(x - (h\mathcal{R}_{\psi(x)}\mathbf{z})_1\big) - \tau(x)\big) K(z_1,z_2)\,\mathrm{d}z_1\mathrm{d}z_2 + \tau(x) + O((nh)^{-1}),$$

where the last equality is due to Assumption 4, i.e. $\int K_1 = \int_0^1 K_2 = 1$. For the first term on the right hand side of the latter display, we obtain that it is $O(h^2)$, which can be seen again by Taylor expansion of $\tau$ and since $x \mapsto K_1(x)x$ is odd. Hence, $J_n = \tau(x) + O(h^2) + O((nh)^{-1})$, and recalling that $\tau(x) = \mathbb{M}(0, \psi(x); x)$ derive

$$\mathcal{B}_n^\tau(x) = \mathbb{M}_n(0, \psi(x); x) - \tau(x) = O(h^2) + O((nh)^{-1}),$$

which completes the proof. □

*A.3. Convergence of the Hessian matrix: proof of Lemma 4*

Similarly to Lemma 1 it is straightforward to calculate the second derivatives.

**Lemma 3.** *It holds for any $x \in [0,1]$,*

$$\partial_w^2 \hat{\mathbb{M}}_n(w, \psi; x)$$
$$= (n^2 h^2)^{-1} \sum_{i_1, i_2=1}^n Y_{i_1, i_2} \langle (\nabla \nabla^T K) \left( h^{-1} \mathcal{R}_{-\psi}(\mathbf{p}(x) + w h \mathbf{e_2} - \mathbf{x_{i_1, i_2}}) \right) (\sin(\psi), \cos(\psi))^T, (\sin(\psi), \cos(\psi))^T \rangle,$$

$$\partial_\psi \partial_w \hat{\mathbb{M}}_n(w, \psi; x)$$
$$= (n^2 h^3)^{-1} \sum_{i_1, i_2=1}^n Y_{i_1, i_2} \langle (\nabla \nabla^T K) \left( h^{-1} \mathcal{R}_{-\psi}(\mathbf{p}(x) + w h \mathbf{e_2} - \mathbf{x_{i_1, i_2}}) \right) (\sin(\psi), \cos(\psi))^T, T(\mathbf{x_{i_1, i_2}}; w, \psi) \rangle$$
$$+ (n^2 h^2)^{-1} \sum_{i_1, i_2=1}^n Y_{i_1, i_2} \langle (\nabla K) \left( h^{-1} \mathcal{R}_{-\psi}(\mathbf{p}(x) + w h \mathbf{e_2} - \mathbf{x_{i_1, i_2}}) \right)^T, (\cos(\psi), -\sin(\psi)) \rangle,$$

$$\partial_\psi^2 \hat{\mathbb{M}}_n(w, \psi; x)$$
$$= (n^2 h^4)^{-1} \sum_{i_1, i_2=1}^n Y_{i_1, i_2} \langle (\nabla \nabla^T K) \left( h^{-1} \mathcal{R}_{-\psi}(\mathbf{p}(x) + w h \mathbf{e_2} - \mathbf{x_{i_1, i_2}}) \right) T(\mathbf{x_{i_1, i_2}}; w, \psi(x)), T(\mathbf{x_{i_1, i_2}}; w, \psi) \rangle$$
$$- (n^2 h^3)^{-1} \sum_{i_1, i_2=1}^n Y_{i_1, i_2} \langle (\nabla K) \left( h^{-1} \mathcal{R}_{-\psi}(\mathbf{p}(x) + w h \mathbf{e_2} - \mathbf{x_{i_1, i_2}}) \right)^T, \mathcal{R}_{-3\pi/2} T(\mathbf{x_{i_1, i_2}}; w, \psi) \rangle.$$

We need the following two auxiliary results to verify the stochastic convergence of the Hessian matrix.

**Lemma 4.** *Let $d \in \mathbb{N}$ and $(\hat{f}_n)_n : \mathbb{R}^d \times I \to \mathbb{R}$ be random functions, which are uniformly Lipschitz, i.e.*

$$\left| \hat{f}_n(\mathbf{z}_1; x) - \hat{f}_n(\mathbf{z}_2; x) \right| \leq L_n |\mathbf{z}_1 - \mathbf{z}_2|,$$

*where $L_n = O_P(1)$ uniformly over $I$. Furthermore, assume there exists $\eta_0 : I \to \mathbb{R}^d$ such that the functions $\hat{f}_n$ are weakly uniformly consistent with the function $f : \mathbb{R}^d \times I \to \mathbb{R}$ at $\eta_0$:*

$$\sup_{x \in I} |\hat{f}_n(\eta_0(x); x) - f(\eta_0(x); x)| \xrightarrow{P} 0.$$

*Then for each sequence of random functions $(\hat{\eta}_n)_{n \in \mathbb{N}}$ with $\hat{\eta}_n : I \to \mathbb{R}^d$ such that $\|\hat{\eta}_n - \eta_0\|_\infty \xrightarrow{P} 0$ it holds that*

$$\sup_{x \in I} |\hat{f}_n(\hat{\eta}_n(x); x) - f(\eta_0(x); x)| \xrightarrow{P} 0.$$

*Proof of Lemma 4.*

$$\sup_{x \in I} \left| \hat{f}_n(\hat{\eta}_n(x); x) - f(\eta_0(x); x) \right| \leq \sup_{x \in I} \left| \hat{f}_n(\hat{\eta}_n(x); x) - \hat{f}_n(\eta_0(x); x) \right| + \sup_{x \in I} \left| \hat{f}_n(\eta_0(x); x) - f(\eta_0(x); x) \right|$$
$$\leq L_n \sup_{x \in I} |\hat{\eta}_n(x) - \eta_0(x)| + o_P(1) = o_P(1).$$

□

**Lemma 5.** *Under the Assumption 4 it holds that*

(i) $\int_{[-1,1]^2} K^{(2,0)}(\mathbf{z}) \, d\mathbf{z} = \int_{[-1,1]^2} K^{(1,1)}(\mathbf{z}) \, d\mathbf{z} = \int_{[-1,1]^2} K^{(0,2)}(\mathbf{z}) \, d\mathbf{z} = 0$;

(ii) $\int_{[-1,1]^2} K^{(2,0)}(z_1, z_2) \, z_2 \, dz_1 dz_2 = \int_{[-1,1]^2} K^{(1,1)}(z_1, z_2) \, z_1 \, dz_1 dz_2 = 0$;

(iii) $\int_{[-1,1]^2} K^{(1,1)}(z_1, z_2) \, z_2 \, dz_1 dz_2 = \int_{[-1,1]^2} K^{(0,2)}(z_1, z_2) \, z_1 \, dz_1 dz_2 = 0$;

(iv) $\int_{[-1,1]^2} K^{(1,0)}(\mathbf{z}) \, d\mathbf{z} = \int_{[-1,1]^2} K^{(0,1)}(\mathbf{z}) \, d\mathbf{z} = 0$;

(v) $\int_{[-1,1]^2} K^{(2,0)}(z_1, z_2) \, z_2^2 \, dz_1 dz_2 = \int_{[-1,1]^2} K^{(1,1)}(z_1, z_2) \, z_1 z_2 \, dz_1 dz_2 = 0$;

(vi) $\int_{[-1,1]^2} K^{(0,2)}(z_1,z_2) \, z_1^2 \, dz_1 dz_2 = 0$;

(vii) $\int_{[-1,1]^2} K^{(1,0)}(z_1,z_2) \, z_1 \, dz_1 dz_2 = \int_{[-1,1]^2} K^{(0,1)}(z_1,z_2) \, z_2 \, dz_1 dz_2 = 0$.

(viii) $\int_{[-1,1]\times[0,1]} K^{(2,0)}(\mathbf{z}) \, d\mathbf{z} = \int_{[-1,1]\times[0,1]} K^{(1,1)}(\mathbf{z}) \, d\mathbf{z} = 0$, $\quad \int_{[-1,1]\times[0,1]} K^{(0,2)}(\mathbf{z}) \, d\mathbf{z} = -K_2^{(1)}(0)$;

(ix) $\int_{[-1,1]\times[0,1]} K^{(2,0)}(z_1,z_2) \, z_2 \, dz_1 dz_2 = \int_{[-1,1]\times[0,1]} K^{(1,1)}(z_1,z_2) \, z_1 \, dz_1 dz_2 = 0$;

(x) $\int_{[-1,1]\times[0,1]} K^{(1,1)}(z_1,z_2) \, z_2 \, dz_1 dz_2 = \int_{[-1,1]\times[0,1]} K^{(0,2)}(z_1,z_2) \, z_1 \, dz_1 dz_2 = 0$;

(xi) $\int_{[-1,1]\times[0,1]} K^{(1,0)}(\mathbf{z}) \, d\mathbf{z} = \int_{[-1,1]\times[0,1]} K^{(0,1)}(\mathbf{z}) \, d\mathbf{z} = 0$;

(xii) $\int_{[-1,1]\times[0,1]} K^{(1,1)}(z_1,z_2) \, z_1 z_2 \, dz_1 dz_2 = 1$, $\quad \int_{[-1,1]\times[0,1]} K^{(2,0)}(z_1,z_2) \, z_2^2 \, dz_1 dz_2 = 0$;

(xiii) $\int_{[-1,1]\times[0,1]} K^{(0,2)}(z_1,z_2) \, z_1^2 \, dz_1 dz_2 = -K_2^{(1)}(0) \int_{-1}^{1} K_1(y) y^2 \, dy$;

(xiv) $\int_{[-1,1]\times[0,1]} K^{(1,0)}(z_1,z_2) \, z_1 \, dz_1 dz_2 = \int_{[-1,1]\times[0,1]} K^{(0,1)}(z_1,z_2) \, z_2 \, dz_1 dz_2 = -1$.

*Proof of Lemma 5.* Some equalities follow by symmetry and normalization of $K_1$ and $K_2$, others require the boundary properties. For example,

$$\int_0^1 K_2^{(2)}(z) \, d\mathbf{z} = K_2^{(1)}(1) - K_2^{(1)}(0) = -K_2^{(1)}(0)$$

since $K_2^{(1)}(1) = 0$ by assumption, which yield the last statement in (viii). $\square$

*Proof of Lemma 4.* We show that

$$\sup_{x\in I} \|\nabla\nabla^T \hat{\mathbb{M}}_n(0,\psi(x);x) - \mathbf{H}(x)\| = O_P(h) + O_P(\sqrt{\log(n)}(nh)^{-1})$$

as well as the stochastic Lipschitz continuity of

$$\hat{f}_n(w,\psi;x) := \int_0^1 \nabla\nabla^T \hat{\mathbb{M}}_n(tw,\psi(x) + t(\psi - \psi(x));x) \, dt$$

on $\Theta_n$. Note that $\hat{f}_n(0,\psi(x);x) = \nabla\nabla^T \hat{\mathbb{M}}_n(0,\psi(x);x)$. The lemma then follows from Lemma 4 together with Proposition 1 by taking $\eta_0(x) = (0,\psi(x))$.

*Uniform stochastic convergence*

As we have to show stochastic convergence of a symmetric matrix, we break it down to showing stochastic convergence of the components.

For sake of brevity we write $\mathbf{v}_{\psi(x)} = (\sin\psi(x), \cos\psi(x))^T$. Also note that all the $O$-terms in Lemma 10, Lemma 11 and the $o_P$-term in Lemma 12 are uniform in $x$, so that all occurring $O$-Terms in the following Steps 1 – 3 hold uniformly in $x$.

*Step 1:*

We show that

$$\sup_{x\in I} |\partial_w^2 \hat{\mathbb{M}}_n(0,\psi(x);x) + \tau(x)\cos^2(\psi(x))K_2^{(1)}(0)| = O_P(h) + O_P(\sqrt{\log(n)}(nh)^{-1}).$$

Use Lemma 3 and split $\partial_w^2 \hat{\mathbb{M}}_n(0,\psi(x);x)$ into three terms

$$\partial_w^2 \hat{\mathbb{M}}_n(0,\psi(x);x) = (n^2h^2)^{-1} \sum_{i_1,i_2=1}^n Y_{i_1,i_2} \langle (\nabla\nabla^T K)(h^{-1}\mathfrak{R}_{-\psi(x)}(\mathbf{p}(x) - \mathbf{x}_{\mathbf{i_1},\mathbf{i_2}})) \mathbf{v}_{\psi(x)}, \mathbf{v}_{\psi(x)} \rangle,$$

$$=: S_n + J_n + E_n,$$

7

where

$$S_n = (nh)^{-2} \sum_{i_1,i_2=1}^{n} m\left(\mathbf{x}_{\mathbf{i_1},\mathbf{i_2}}\right) \langle (\nabla\nabla^T K)\left(h^{-1}\mathcal{R}_{-\psi(x)}\left(\mathbf{p}(x) - \mathbf{x}_{\mathbf{i_1},\mathbf{i_2}}\right)\right) \mathbf{v}_{\psi(x)}, \mathbf{v}_{\psi(x)} \rangle,$$

$$J_n = (nh)^{-2} \sum_{i_1,i_2=1}^{n} j_\tau\left(\mathbf{x}_{\mathbf{i_1},\mathbf{i_2}}\right) \langle (\nabla\nabla^T K)\left(h^{-1}\mathcal{R}_{-\psi(x)}\left(\mathbf{p}(x) - \mathbf{x}_{\mathbf{i_1},\mathbf{i_2}}\right)\right) \mathbf{v}_{\psi(x)}, \mathbf{v}_{\psi(x)} \rangle,$$

$$E_n = (nh)^{-2} \sum_{i_1,i_2=1}^{n} \varepsilon_{i_1,i_2} \langle (\nabla\nabla^T K)\left(h^{-1}\mathcal{R}_{-\psi(x)}\left(\mathbf{p}(x) - \mathbf{x}_{\mathbf{i_1},\mathbf{i_2}}\right)\right) \mathbf{v}_{\psi(x)}, \mathbf{v}_{\psi(x)} \rangle.$$

Note we dropped the inputs $(x, w, \psi)$ for the latter terms for sake of convenience. By Lemma 10 with $g_1(\mathbf{z}) = m(\mathbf{z})$, $g_2(\mathbf{z}) = g_3(\mathbf{z}) \equiv \mathbf{v}_{\psi(x)}$, $f(\mathbf{z}) = \nabla\nabla^T K(\mathbf{z})$ such that $r_1 = 2, r_2 = r_3 = 0, j = 1$, obtaining

$$S_n = m(\mathbf{p}(x)) \int_{[-1,1]^2} \langle \nabla\nabla^T K(\mathbf{z}) \mathbf{v}_{\psi(x)}, \mathbf{v}_{\psi(x)} \rangle d\mathbf{z} + O(h) + O\left((nh)^{-1}\right)$$
$$= O(h) + O\left((nh)^{-1}\right),$$

where the last equation is due to Lemma 5 (i). Applying Lemma 11 with $f, g_2$ and $g_3$ as above and Lemma 5 (viii) yield

$$J_n = \tau(x) \int_{[-1,1]\times[0,1]} \langle \nabla\nabla^T K(\mathbf{z}) \mathbf{v}_{\psi(x)}, \mathbf{v}_{\psi(x)} \rangle d\mathbf{z} + O(h) + O\left((nh)^{-1}\right)$$
$$= -\tau(x) \cos(\psi(x))^2 K_2^{(1)}(0) + O(h) + O\left((nh)^{-1}\right).$$

Using Lemma 12 with the same functions, we get immediately $E_n = O_P(\sqrt{\log(n)}(nh)^{-1})$. Hence,

$$\partial_w^2 \hat{\mathbb{M}}_n(0, \psi(x); x) = S_n + J_n + E_n \xrightarrow{P} -\tau(x)\cos(\psi(x))^2 K_2^{(1)}(0).$$

As all the $O$-terms are uniform in $x$ we have actually stochastic uniform convergence in the preceding display.

*Step 2:*

We show that
$$\sup_{x \in I} |\partial_w \partial_\psi \hat{\mathbb{M}}_n(0, \psi(x); x)| = O_P(h) + O_P(\sqrt{\log(n)}(nh)^{-1}).$$

Use Lemma 3 and split $\partial_w \partial_\psi \hat{\mathbb{M}}_n(0, \psi(x); x)$ into six terms

$$\partial_w \partial_\psi \hat{\mathbb{M}}_n(0, \psi(x); x)$$
$$= (n^2 h^3)^{-1} \sum_{i_1,i_2=1}^{n} Y_{i_1,i_2} \langle (\nabla\nabla^T K)\left(h^{-1}\mathcal{R}_{-\psi(x)}\left(\mathbf{p}(x) - \mathbf{x}_{\mathbf{i_1},\mathbf{i_2}}\right)\right) \mathbf{v}_{\psi(x)}, T\left(\mathbf{x}_{\mathbf{i_1},\mathbf{i_2}}; 0, \psi(x)\right) \rangle$$
$$+ (n^2 h^2)^{-1} \sum_{i_1,i_2=1}^{n} Y_{i_1,i_2} \langle (\nabla K)\left(h^{-1}\mathcal{R}_{-\psi(x)}\left(\mathbf{p}(x) - \mathbf{x}_{\mathbf{i_1},\mathbf{i_2}}\right)\right)^T, (\cos(\psi(x)), -\sin(\psi(x))) \rangle,$$
$$=: S_{n,1} + J_{n,1} + E_{n,1} + S_{n,2} + J_{n,2} + E_{n,2},$$

8

where

$$S_{n,1} = (n^2h^3)^{-1} \sum_{i_1,i_2=1}^{n} m\left(\mathbf{x_{i_1,i_2}}\right) \langle (\nabla\nabla^T K)\left(h^{-1}\mathcal{R}_{-\psi(x)}\left(\mathbf{p}(x) - \mathbf{x_{i_1,i_2}}\right)\right) \mathbf{v}_{\psi(x)}, T\left(\mathbf{x_{i_1,i_2}}; 0, \psi(x)\right)\rangle,$$

$$J_{n,1} = (n^2h^3)^{-1} \sum_{i_1,i_2=1}^{n} j_\tau\left(\mathbf{x_{i_1,i_2}}\right) \langle (\nabla\nabla^T K)\left(h^{-1}\mathcal{R}_{-\psi(x)}\left(\mathbf{p}(x) - \mathbf{x_{i_1,i_2}}\right)\right) \mathbf{v}_{\psi(x)}, T\left(\mathbf{x_{i_1,i_2}}; 0, \psi(x)\right)\rangle,$$

$$E_{n,1} = (n^2h^3)^{-1} \sum_{i_1,i_2=1}^{n} \varepsilon_{i_1,i_2} \langle (\nabla\nabla^T K)\left(h^{-1}\mathcal{R}_{-\psi(x)}\left(\mathbf{p}(x) - \mathbf{x_{i_1,i_2}}\right)\right) \mathbf{v}_{\psi(x)}, T\left(\mathbf{x_{i_1,i_2}}; 0, \psi(x)\right)\rangle,$$

$$S_{n,2} = (n^2h^2)^{-1} \sum_{i_1,i_2=1}^{n} m\left(\mathbf{x_{i_1,i_2}}\right) \langle (\nabla K)\left(h^{-1}\mathcal{R}_{-\psi(x)}\left(\mathbf{p}(x) - \mathbf{x_{i_1,i_2}}\right)\right)^T, (\cos(\psi(x)), -\sin(\psi(x)))\rangle,$$

$$J_{n,2} = (n^2h^2)^{-1} \sum_{i_1,i_2=1}^{n} j_\tau\left(\mathbf{x_{i_1,i_2}}\right) \langle (\nabla K)\left(h^{-1}\mathcal{R}_{-\psi(x)}\left(\mathbf{p}(x) - \mathbf{x_{i_1,i_2}}\right)\right)^T, (\cos(\psi(x)), -\sin(\psi(x)))\rangle,$$

$$E_{n,2} = (n^2h^2)^{-1} \sum_{i_1,i_2=1}^{n} \varepsilon_{i_1,i_2} \langle (\nabla K)\left(h^{-1}\mathcal{R}_{-\psi(x)}\left(\mathbf{p}(x) - \mathbf{x_{i_1,i_2}}\right)\right)^T, (\cos(\psi(x)), -\sin(\psi(x)))\rangle.$$

By Lemma 10 with $g_1(\mathbf{z}) = m(\mathbf{z})$, $f(\mathbf{z}) = \nabla\nabla^T K(\mathbf{z})$, $g_2 \equiv \mathbf{v}_{\psi(x)}$, $g_3(\mathbf{z}) = T(\mathbf{z}; 0, \psi(x))$, where $r_1 = 3$, $r_2 = 0$, $r_3 = j = 1$, together with (2) lead to

$$S_{n,1} = m(\mathbf{p}(x)) \int_{[-1,1]^2} \langle \nabla\nabla^T K(\mathbf{z})\mathbf{v}_{\psi(x)}, (z_2, -z_1)^T\rangle \, dz_1 dz_2 + O(h) + O\left((nh)^{-1}\right)$$
$$= O(h) + O\left((nh)^{-1}\right),$$

where the last line is due to Lemma 5 (ii) and (iii). Similarly, applying Lemma 11 we get

$$J_{n,1} = \tau(x) \int_{[-1,1]\times[0,1]} \langle \nabla\nabla^T K(z_1,z_2)\mathbf{v}_{\psi(x)}, (z_2, -z_1)^T\rangle \, dz_1 dz_2 + O(h) + O\left((nh)^{-1}\right)$$
$$= O(h) + O\left((nh)^{-1}\right),$$

in which the second equality follows by Lemma 5 (ix) and (x). Using Lemma 12 we get immediately $E_{n,1} = O_P(\sqrt{\log(n)}(nh)^{-1})$. Further, we use Lemma 10 with $g_1(\mathbf{z}) = m(\mathbf{z})$, $f(\mathbf{z}) = [\nabla K(\mathbf{z}) \mid \nabla K(\mathbf{z})]$, $g_2 \equiv \mathbf{e}_1$ and $g_3(\mathbf{z}) = (\cos(\psi(x)), -\sin(\psi(x)))^T$, so that $r_1 = 2$, $r_2 = r_3 = 0$, $j = 1$, which together with Lemma 5 (xi) implies

$$S_{n,2} = m(\mathbf{p}(x)) \int_{[-1,1]^2} \langle \nabla K(\mathbf{z}), (\cos(\psi(x)), -\sin(\psi(x)))\rangle \, d\mathbf{z} + O(h) + O\left((nh)^{-1}\right)$$
$$= O(h) + O\left((nh)^{-1}\right).$$

Using for Lemma 11 the same functions as just and Lemma 5 (xi) one gets

$$J_{n,2} = \tau(x) \int_{[-1,1]\times[0,1]} \langle \nabla K(\mathbf{z}), (\cos(\psi(x)), -\sin(\psi(x)))\rangle \, d\mathbf{z} + O(h) + O\left((nh)^{-1}\right)$$
$$= O(h) + O\left((nh)^{-1}\right).$$

By means of Lemma 12, obtain $E_{n,2} = O_P(\sqrt{\log(n)}(nh)^{-1})$. Summarizing, we deduce that

$$\partial_w \partial_\psi \hat{\mathbb{M}}_n(0, \psi(x); x) = S_{n,1} + J_{n,1} + E_{n,1} + S_{n,2} + J_{n,2} + E_{n,2} \xrightarrow{P} 0.$$

As the $O$-terms are all uniform in $x$ the latter display holds uniformly in $x$.

*Step 3:*
It remains to prove

$$\sup_{x \in I} |\partial_\psi^2 \hat{\mathbb{M}}_n(0, \psi(x); x) + \tau(x) K_2^{(1)}(0) \int_{-1}^1 K_1(y) y^2 \, dy| = O_P(h) + O_P(\sqrt{\log(n)}(nh)^{-1}).$$

Use Lemma 3 and split $\partial_\psi^2 \hat{\mathbb{M}}_n(0, \psi(x); x)$ into six terms

$$\partial_\psi^2 \hat{\mathbb{M}}_n(0, \psi(x); x)$$
$$= (n^2 h^4)^{-1} \sum_{i_1, i_2=1}^n Y_{i_1, i_2} \langle (\nabla \nabla^T K) (h^{-1} \mathcal{R}_{-\psi(x)} (\mathbf{p}(x) - \mathbf{x}_{\mathbf{i}_1, \mathbf{i}_2})) T(\mathbf{x}_{\mathbf{i}_1, \mathbf{i}_2}; 0, \psi(x)), T(\mathbf{x}_{\mathbf{i}_1, \mathbf{i}_2}; 0, \psi(x)) \rangle$$
$$- (n^2 h^3)^{-1} \sum_{i_1, i_2=1}^n Y_{i_1, i_2} \langle (\nabla K) (h^{-1} \mathcal{R}_{-\psi(x)} (\mathbf{p}(x) - \mathbf{x}_{\mathbf{i}_1, \mathbf{i}_2})), \mathcal{R}_{-3\pi/2} T(\mathbf{x}_{\mathbf{i}_1, \mathbf{i}_2}; 0, \psi(x)) \rangle$$
$$=: S_{n,1} + J_{n,1} + E_{n,1} - S_{n,2} - J_{n,2} - E_{n,2},$$

in which

$$S_{n,1} = (n^2 h^4)^{-1} \sum_{i_1, i_2=1}^n m(\mathbf{x}_{\mathbf{i}_1, \mathbf{i}_2}) \langle (\nabla \nabla^T K) (h^{-1} \mathcal{R}_{-\psi(x)} (\mathbf{p}(x) - \mathbf{x}_{\mathbf{i}_1, \mathbf{i}_2})) T(\mathbf{x}_{\mathbf{i}_1, \mathbf{i}_2}; 0, \psi(x)), T(\mathbf{x}_{\mathbf{i}_1, \mathbf{i}_2}; 0, \psi(x)) \rangle,$$

$$J_{n,1} = (n^2 h^4)^{-1} \sum_{i_1, i_2=1}^n j_\tau(\mathbf{x}_{\mathbf{i}_1, \mathbf{i}_2}) \langle (\nabla \nabla^T K) (h^{-1} \mathcal{R}_{-\psi(x)} (\mathbf{p}(x) - \mathbf{x}_{\mathbf{i}_1, \mathbf{i}_2})) T(\mathbf{x}_{\mathbf{i}_1, \mathbf{i}_2}; 0, \psi(x)), T(\mathbf{x}_{\mathbf{i}_1, \mathbf{i}_2}; 0, \psi(x)) \rangle,$$

$$E_{n,1} = (n^2 h^4)^{-1} \sum_{i_1, i_2=1}^n \varepsilon_{i_1, i_2} \langle (\nabla \nabla^T K) (h^{-1} \mathcal{R}_{-\psi(x)} (\mathbf{p}(x) - \mathbf{x}_{\mathbf{i}_1, \mathbf{i}_2})) T(\mathbf{x}_{\mathbf{i}_1, \mathbf{i}_2}; 0, \psi(x)), T(\mathbf{x}_{\mathbf{i}_1, \mathbf{i}_2}; 0, \psi(x)) \rangle,$$

$$S_{n,2} = (n^2 h^3)^{-1} \sum_{i_1, i_2=1}^n m(\mathbf{x}_{\mathbf{i}_1, \mathbf{i}_2}) \langle (\nabla K) (h^{-1} \mathcal{R}_{-\psi(x)} (\mathbf{p}(x) - \mathbf{x}_{\mathbf{i}_1, \mathbf{i}_2})), \mathcal{R}_{-3\pi/2} T(\mathbf{x}_{\mathbf{i}_1, \mathbf{i}_2}; 0, \psi(x)) \rangle,$$

$$J_{n,2} = (n^2 h^3)^{-1} \sum_{i_1, i_2=1}^n j_\tau(\mathbf{x}_{\mathbf{i}_1, \mathbf{i}_2}) \langle (\nabla K) (h^{-1} \mathcal{R}_{-\psi(x)} (\mathbf{p}(x) - \mathbf{x}_{\mathbf{i}_1, \mathbf{i}_2})), \mathcal{R}_{-3\pi/2} T(\mathbf{x}_{\mathbf{i}_1, \mathbf{i}_2}; 0, \psi(x)) \rangle,$$

$$E_{n,2} = (n^2 h^3)^{-1} \sum_{i_1, i_2=1}^n \varepsilon_{i_1, i_2} \langle (\nabla K) (h^{-1} \mathcal{R}_{-\psi(x)} (\mathbf{p}(x) - \mathbf{x}_{\mathbf{i}_1, \mathbf{i}_2})), \mathcal{R}_{-3\pi/2} T(\mathbf{x}_{\mathbf{i}_1, \mathbf{i}_2}; 0, \psi(x)) \rangle.$$

Application of Lemma 10 with the functions $g_1(\mathbf{z}) = m(\mathbf{z})$, $f(\mathbf{z}) = \nabla \nabla^T K(\mathbf{z})$, $g_2(\mathbf{z}) = g_3(\mathbf{z}) = T(\mathbf{z}; 0, \psi(x))$ ( $r_1 = 4$, $r_2 = r_3 = 1$, $j = 1$) yields with (2)

$$S_{n,1} = m(\mathbf{p}(x)) \int_{[-1,1]^2} \langle \nabla \nabla^T K(z_1, z_2) (z_2, -z_1)^T, (z_2, -z_1)^T \rangle \, dz_1 dz_2 + O(h) + O((nh)^{-1})$$
$$= O(h) + O((nh)^{-1}),$$

where the second equality is by Lemma 5 (v) and (vi). With $g_2$ and $f$ as just and Lemma 11 combined with Lemma 5 (xii) and (xiii),

$$J_{n,1} = \tau(x) \int_{[-1,1] \times [0,1]} \langle \nabla \nabla^T K(z_1, z_2) (z_2, -z_1)^T, (z_2, -z_1)^T \rangle \, dz_1 dz_2 + O(h) + O((nh)^{-1})$$
$$= -2\tau(x) - \tau(x) K_2^{(1)}(0) \int_{-1}^1 K_1(y) y^2 \, dy + O(h) + O((nh)^{-1}).$$

By Lemma 12 obtain that $E_{n,1} = O_P(\sqrt{\log(n)}(nh)^{-1})$. Next, using $g_1(\mathbf{z}) = m(\mathbf{z})$, $f(\mathbf{z}) = [\nabla K(\mathbf{z}) \mid \nabla K(\mathbf{z})]$, $g_2 \equiv \mathbf{e}_1$

and $g_3(\mathbf{z}) = \mathcal{R}_{-3\pi/2} T(\mathbf{z}; 0, \psi(x))$ in Lemma 10 ($r_1 = 3$, $r_2 = 0$, $r_3 = j = 1$) and Lemma 5 (vii) yield

$$S_{n,2} = m(\mathbf{p}(x)) \int_{[-1,1]^2} \langle \nabla K(z_1, z_2), (z_1, z_2)^T \rangle \, dz_1 dz_2 + O(h) + O((nh)^{-1})$$
$$= O(h) + O((nh)^{-1}).$$

With $g_2, g_3$ and $f$ as just, we get by applying Lemma 11 and Lemma 5 (xiv),

$$J_{n,2} = \tau(x) \int_{[-1,1]\times[0,1]} \langle \nabla K\left((z_1, z_2)^T\right), (z_1, z_2)^T \rangle \, dz_1 dz_2 + O(h) + O((nh)^{-1})$$
$$= -2\tau(x) + O(h) + O((nh)^{-1}).$$

and $E_{n,2} = O_P(\sqrt{\log(n)}(nh)^{-1})$ from Lemma 12. Finally, one obtains uniformly in $x$ that

$$\partial_\psi^2 \hat{\mathbb{M}}_n(0, \psi(x); x) = S_{n,1} + J_{n,1} + E_{n,1} - S_{n,2} - J_{n,2} - E_{n,2} \xrightarrow{P} -\tau(x) K_2^{(1)}(0) \int_{-1}^1 K_1(y) y^2 \, dy.$$

*Lipschitz continuity*

It suffices to show that

$$\|\nabla \nabla^T \hat{\mathbb{M}}_n(w_1, \psi_1; x) - \nabla \nabla^T \hat{\mathbb{M}}_n(w_2, \psi_2; x)\| \leq L_n \|(w_1, \psi_1)^T - (w_2, \psi_2)^T\|, \tag{4}$$

where $L_n = O_P(1)$ and $(w_i, \psi_i) \in \Theta_n$, $i = 1, 2$. By taking

$$L_n = \sup_{(x,w,\psi)^T \in \tilde{\Theta}_{n,x}} \|(\nabla) \otimes (\nabla \nabla^T) \hat{\mathbb{M}}_n(w, \psi; x)\|,$$

where $\otimes$ is the Kronecker product, we may obtain (4) by the mean value theorem. One can show that all components of $(\nabla) \otimes (\nabla \nabla^T) \hat{\mathbb{M}}_n(w, \psi; x)$ are uniformly bounded in probability. For example, letting as above $\mathbf{v}_\psi = (\sin(\psi), \cos(\psi))$, then the third partial derivative for $w$ is

$$\partial_w^3 \hat{\mathbb{M}}_n(w, \psi; x)$$
$$= (nh)^{-2} \sum_{i_1, i_2 = 1}^n Y_{i_1, i_2} \langle (\mathbf{v}_\psi^T \otimes \nabla) \otimes (\nabla \nabla^T K)(h^{-1} \mathcal{R}_{-\psi}(\mathbf{p}(x) + wh\mathbf{e_2} - \mathbf{x_{i_1, i_2}}))(\mathbf{v}_\psi, \mathbf{v}_\psi)^T, \mathbf{v}_\psi^T \rangle.$$

Let $I(n, h)$ denote the set of indices for which the latter sum is not zero. Hence, $|I(n, h)| \leq 4n^2 h^2$ and therefore

$$\mathbb{E} \sup_{(x,w,\psi)^T \in \Theta_n} |\partial_w^3 \hat{\mathbb{M}}_n(w, \psi; x)| \leq 4(nh)^{-2} \|(\nabla \otimes \nabla \nabla^T) K\|_\infty \sum_{i_1, i_2 \in I(n,h)} \mathbb{E}|Y_{i_1, i_2}| = O(1).$$

The remaining partial derivatives are dealt with similarly.

*Deterministic Hessian matrix*

Note that $\nabla \nabla^T \mathbb{M}_n(0, \psi(x); x)$ has the same components as $\nabla \nabla^T \hat{\mathbb{M}}_n(0, \psi(x); x)$ besides the stochastic parts, which were denoted by $E_n$ in the proof steps above and which were stochastically negligible. Hence, $\nabla \nabla^T \mathbb{M}_n(0, \psi(x); x)$ has the same limit as $\nabla \nabla^T \hat{\mathbb{M}}_n(0, \psi(x); x)$ and this implies the convergence of the deterministic Hessian matrix $\mathbf{H}_n(x)$ against the same limit as the stochastic Hessian matrix $\hat{\mathbf{H}}_n(x)$.

□

*A.4. Properties of the score: proofs of Proposition 3 and of Lemma 5*

*Proof of Proposition 3.* For the first component of $\mathcal{S}_n^{\phi,\psi}(0,\psi(x);x)$, i.e. for

$$\mathcal{S}_n^{\phi}(0,\psi(x);x) = \partial_w(\hat{\mathbb{M}}_n(0,\psi(x);x) - \mathbb{M}_n(0,\psi(x);x))$$
$$= (nh)^{-2} \sum_{i_1,i_2=1}^{n} \varepsilon_{i_1,i_2} \langle (\nabla K)\left(h^{-1}\mathfrak{R}_{-\psi(x)}(\mathbf{p}(x) - \mathbf{x}_{\mathbf{i_1,i_2}})\right), (\sin\psi(x),\cos\psi(x))^T \rangle$$

use Lemma 12 with $f(\mathbf{z}) = [\nabla K(\mathbf{z}) \mid \nabla K(\mathbf{z})]$, $g_2 \equiv \mathbf{e}_1$, $g_3 \equiv (\sin\psi(x),\cos\psi(x))^T$, such that $j=1$, $r_1 = 2$, $r_2 = r_3 = 0$, and obtain that

$$\sup_{x \in I} \|\mathcal{S}_n^{\phi}(0,\psi(x);x)\| = O_P\left(\sqrt{\log(n)}(nh)^{-1}\right).$$

For the second component of $\mathcal{S}_n^{\phi,\psi}(0,\psi(x);x)$, that is for $\mathcal{S}_n^{\psi}(0,\psi(x);x)$, one can proceed analogously by using Lemma 12 with $f$ and $g_2$ as before, while $g_3(\mathbf{z}) = T(\mathbf{z};0,\psi(x)))$.

For the second claim recall from (25) that

$$\mathcal{S}_n^{\tau}(0,\psi(x);x) = \hat{\mathbb{M}}_n(0,\psi(x);x) - \mathbb{M}_n(0,\psi(x);x) = (nh)^{-2} \sum_{i_1,i_2=1}^{n} \varepsilon_{i_1,i_2} K(\mathbf{z}) \, d\mathbf{z}.$$

Thus, Lemma 12 with

$$g_1(\mathbf{z}) \equiv 1, \quad f(\mathbf{z}) = 4^{-1} \begin{pmatrix} K(\mathbf{z}) & K(\mathbf{z}) \\ K(\mathbf{z}) & K(\mathbf{z}) \end{pmatrix}, \quad g_2(\mathbf{z}) = g_3(\mathbf{z}) = (1,1)^T, j=1$$

such that $r_1 = 2$, $r_2 = r_3 = 0$, implies the assertion. □

*Proof of Lemma 5.* Let $x \in I$. We intend to make use of the Lindeberg-Feller Theorem (see, e.g. van der Vaart [5], Proposition 2.27). By Lemma 1

$$nh \begin{pmatrix} \mathcal{S}_n^{\phi}(0,\psi(x);x) \\ \mathcal{S}_n^{\psi}(0,\psi(x);x) \\ \mathcal{S}_n^{\tau}(0,\psi(x);x) \end{pmatrix} = \sum_{i_1,i_2=1}^{n} \mathbf{a}_{\mathbf{i_1,i_2}}(x)\varepsilon_{i_1,i_2},$$

with

$$\left(\mathbf{a}_{\mathbf{i_1,i_2}}(x)\right)_1 = (nh)^{-1} \langle (\nabla K)\left(h^{-1}\mathfrak{R}_{-\psi(x)}(\mathbf{p}(x) - \mathbf{x}_{\mathbf{i_1,i_2}})\right), (\sin(\psi(x)),\cos(\psi(x)))^T \rangle,$$
$$\left(\mathbf{a}_{\mathbf{i_1,i_2}}(x)\right)_2 = (nh^2)^{-1} \langle (\nabla K)\left(h^{-1}\mathfrak{R}_{-\psi(x)}(\mathbf{p}(x) - \mathbf{x}_{\mathbf{i_1,i_2}})\right), T(\mathbf{x}_{\mathbf{i_1,i_2}};0,\psi(x)) \rangle,$$
$$\left(\mathbf{a}_{\mathbf{i_1,i_2}}(x)\right)_3 = (nh)^{-1} K\left(h^{-1}\mathfrak{R}_{-\psi(x)}(\mathbf{p}(x) - \mathbf{x}_{\mathbf{i_1,i_2}})\right).$$

First, we consider the convergence of the covariance matrix,

$$\sum_{i_1,i_2=1}^{n} \operatorname{Cov}\left(\mathbf{a}_{\mathbf{i_1,i_2}}(x)\varepsilon_{i_1,i_2}\right) = \sigma^2 \sum_{i_1,i_2=1}^{n} \begin{pmatrix} \left(\mathbf{a}_{\mathbf{i_1,i_2}}(x)\right)_1^2 & \left(\mathbf{a}_{\mathbf{i_1,i_2}}(x)\right)_1\left(\mathbf{a}_{\mathbf{i_1,i_2}}(x)\right)_2 & \left(\mathbf{a}_{\mathbf{i_1,i_2}}(x)\right)_1\left(\mathbf{a}_{\mathbf{i_1,i_2}}(x)\right)_3 \\ \left(\mathbf{a}_{\mathbf{i_1,i_2}}(x)\right)_1\left(\mathbf{a}_{\mathbf{i_1,i_2}}(x)\right)_2 & \left(\mathbf{a}_{\mathbf{i_1,i_2}}(x)\right)_2^2 & \left(\mathbf{a}_{\mathbf{i_1,i_2}}(x)\right)_2\left(\mathbf{a}_{\mathbf{i_1,i_2}}(x)\right)_3 \\ \left(\mathbf{a}_{\mathbf{i_1,i_2}}(x)\right)_1\left(\mathbf{a}_{\mathbf{i_1,i_2}}(x)\right)_3 & \left(\mathbf{a}_{\mathbf{i_1,i_2}}(x)\right)_2\left(\mathbf{a}_{\mathbf{i_1,i_2}}(x)\right)_3 & \left(\mathbf{a}_{\mathbf{i_1,i_2}}(x)\right)_3^2 \end{pmatrix}$$
$$\to \sigma^2 \operatorname{diag}(V_\phi^S(x), V_\psi^S, V_\tau^S)).$$

*Step 1:*
We show

$$\sum_{i_1,i_2=1}^{n} \left(\mathbf{a}_{\mathbf{i_1,i_2}}(x)\right)_1^2 \to V_\phi^S(x).$$

Applying Lemma 10 with $g_1 \equiv 1$, $f(\mathbf{z}) = [\nabla K(\mathbf{z}) \mid \nabla K(\mathbf{z})] \in \mathbb{R}^{2 \times 2}$, $g_2 \equiv \mathbf{e}_1$, $g_3 \equiv (\sin\psi(x),\cos\psi(x))^T$, such that

12

$r_1 = 2$, $r_2 = r_3 = 0$, and $j = 2$, we obtain

$$\sum_{i_1,i_2=1}^{n} \left(\mathbf{a}_{\mathbf{i}_1,\mathbf{i}_2}(x)\right)_1^2 = \int_{[-1,1]^2} \langle \nabla K(\mathbf{z}), (\sin \psi(x), \cos \psi(x))^T \rangle^2 \, d\mathbf{z} + O(h) + O\left((nh)^{-1}\right)$$

$$= \sin^2(\psi(x)) \int_{[-1,1]^2} K^{(1,0)}(\mathbf{z})^2 \, d\mathbf{z} + 2\cos(\psi(x))\sin(\psi(x)) \int_{[-1,1]^2} K^{(1,0)}(\mathbf{z}) K^{(0,1)}(\mathbf{z}) \, d\mathbf{z}$$

$$+ \cos^2(\psi(x)) \int_{[-1,1]^2} K^{(0,1)}(\mathbf{z})^2 \, d\mathbf{z} + O(h) + O\left((nh)^{-1}\right)$$

$$= V_\phi^S(x) + O(h) + O\left((nh)^{-1}\right),$$

since $y \mapsto K_i^{(1)}(y) K_i(y)$ are odd functions, $i = 1, 2$.

*Step 2:*
We show
$$\sum_{i_1,i_2=1}^{n} \left(\mathbf{a}_{\mathbf{i}_1,\mathbf{i}_2}(x)\right)_2^2 \to V_\psi^S.$$

Applying Lemma 10 with $g_1 \equiv 1$, $f(\mathbf{z}) = [\nabla K(\mathbf{z}) \mid \nabla K(\mathbf{z})]$, $g_2 \equiv \mathbf{e}_1$, $g_3(\mathbf{z}) = T(\mathbf{z}; 0, \psi(x))$, such that $r_1 = 4$, $r_2 = 0$, $r_3 = 1$, $j = 2$, and using (2) lead to

$$\sum_{i_1,i_2=1}^{n} \left(\mathbf{a}_{\mathbf{i}_1,\mathbf{i}_2}(x)\right)_2^2 = \int_{[-1,1]^2} \langle \nabla K(z_1, z_2), (z_2, -z_1)^T \rangle^2 \, dz_1 dz_2 + O(h) + O\left((nh^2)^{-1}\right)$$

$$= V_\psi^S + O(h) + O\left((nh^2)^{-1}\right).$$

*Step 3:*
Our aim is verifying
$$\sum_{i_1,i_2=1}^{n} \left(\mathbf{a}_{\mathbf{i}_1,\mathbf{i}_2}(x)\right)_3^2 \to V_\tau^S.$$

Using Lemma 10 with

$$g_1(\mathbf{z}) \equiv 1, \quad f(\mathbf{z}) = 4^{-1} \begin{pmatrix} K(\mathbf{z}) & K(\mathbf{z}) \\ K(\mathbf{z}) & K(\mathbf{z}) \end{pmatrix}, \quad g_2(\mathbf{z}) = g_3(\mathbf{z}) = (1,1)^T, j = 2$$

such that $r_1 = 2$, $r_2 = r_3 = 0$, one has that

$$\sum_{i_1,i_2=1}^{n} \left(\mathbf{a}_{\mathbf{i}_1,\mathbf{i}_2}(x)\right)_3^2 = \int_{[-1,1]^2} K(\mathbf{z})^2 \, d\mathbf{z} + O(h) + O\left((nh)^{-1}\right) = V_\tau^S + O(h) + O((nh)^{-1}).$$

*Step 4:*
We verify
$$\sum_{i_1,i_2=1}^{n} \left(\mathbf{a}_{\mathbf{i}_1,\mathbf{i}_2}(x)\right)_1 \left(\mathbf{a}_{\mathbf{i}_1,\mathbf{i}_2}(x)\right)_2 \to 0.$$

Applying the more general part of Lemma 10 with $g_1 \equiv 1$, $f_1(\mathbf{z}) = f_2(\mathbf{z}) = [\nabla K(\mathbf{z}) \mid \nabla K(\mathbf{z})]$, $g_2 = g_4 \equiv \mathbf{e}_1$, $g_3(\mathbf{z}) =$

$\left( \sin \psi(x), \cos \psi(x) \right)^T$ and $g_5(\mathbf{z}) = T(\mathbf{z}; 0, \psi(x))$, such that $r_1 = 3$, $r_2 = r_3 = r_4 = 0$, $r_5 = 1$, and using (2) obtain

$$\sum_{i_1,i_2=1}^{n} \left(\mathbf{a}_{\mathbf{i_1},\mathbf{i_2}}(x)\right)_1 \left(\mathbf{a}_{\mathbf{i_1},\mathbf{i_2}}(x)\right)_2$$
$$= \int_{[-1,1]^2} \langle \nabla K(z_1,z_2), (\sin \psi(x), \cos \psi(x))^T \rangle \langle \nabla K(z_1,z_2), (z_2, -z_1)^T \rangle \, dz_1 dz_2$$
$$+ O(h) + O\left((nh^2)^{-1}\right)$$
$$= \sin(\psi(x)) \int_{[-1,1]^2} K^{(1,0)}(z_1,z_2)^2 z_2 \, dz_1 dz_2$$
$$- \sin(\psi(x)) \int_{[-1,1]^2} K^{(1,0)}(z_1,z_2) K^{(0,1)}(z_1,z_2) z_1 \, dz_1 dz_2$$
$$+ \cos(\psi(x)) \int_{[-1,1]^2} K^{(1,0)}(z_1,z_2) K^{(0,1)}(z_1,z_2) z_2 \, dz_1 dz_2$$
$$- \cos(\psi(x)) \int_{[-1,1]^2} K^{(0,1)}(z_1,z_2)^2 z_1 \, dz_1 dz_2 + O(h) + O\left((nh^2)^{-1}\right)$$
$$= O(h) + O\left((nh^2)^{-1}\right),$$

where the last line follows by the fact that $y \mapsto K_i^2(y) y$, $i = 1, 2$, are odd, so the first and fourth integral vanish, and $y \mapsto K_i^{(1)}(y) K_i(y)$, $i = 1, 2$, are odd as well, so that the second and third integral also vanish.

Step 5:

We verify
$$\sum_{i_1,i_2=1}^{n} \left(\mathbf{a}_{\mathbf{i_1},\mathbf{i_2}}(x)\right)_1 \left(\mathbf{a}_{\mathbf{i_1},\mathbf{i_2}}(x)\right)_3 \to 0.$$

Applying the more general part of Lemma 10 with $f_1(\mathbf{z}) = [\nabla K(\mathbf{z}) \mid \nabla K(\mathbf{z})]$, $f_2(\mathbf{z}) = 4^{-1} \begin{pmatrix} K(\mathbf{z}) & K(\mathbf{z}) \\ K(\mathbf{z}) & K(\mathbf{z}) \end{pmatrix}$, $g_1 \equiv 1$, $g_2 = \mathbf{e_1}$, $g_3(\mathbf{z}) = \left( \sin \psi(x), \cos \psi(x) \right)^T$, $g_4 = g_5 \equiv (1,1)^T$, such that $r_1 = 2$, $r_2 = r_3 = r_4 = r_5 = 0$, we get that

$$\sum_{i_1,i_2=1}^{n} \left(\mathbf{a}_{\mathbf{i_1},\mathbf{i_2}}(x)\right)_1 \left(\mathbf{a}_{\mathbf{i_1},\mathbf{i_2}}(x)\right)_3$$
$$= \int_{[-1,1]^2} \langle \nabla K(z_1,z_2), \left( \sin \psi(x), \cos \psi(x) \right)^T \rangle \times K(z_1,z_2) \, dz_1 dz_2 + O(h) + O\left((nh)^{-1}\right)$$
$$= O(h) + O\left((nh)^{-1}\right),$$

where the last equation is true, since $x \mapsto K_1^{(1)}(x) K_1(x)$ and $x \mapsto K_2^{(1)}(x) K_2(x)$ are both odd functions.

Step 6:

We show that
$$\sum_{i_1,i_2=1}^{n} \left(\mathbf{a}_{\mathbf{i_1},\mathbf{i_2}}(x)\right)_2 \left(\mathbf{a}_{\mathbf{i_1},\mathbf{i_2}}(x)\right)_3 \to 0.$$

Applying the more general part of Lemma 10 with $f_1(\mathbf{z}) = [\nabla K(\mathbf{z}) \mid \nabla K(\mathbf{z})]$, $f_2(\mathbf{z}) = 4^{-1} \begin{pmatrix} K(\mathbf{z}) & K(\mathbf{z}) \\ K(\mathbf{z}) & K(\mathbf{z}) \end{pmatrix}$, $g_1 \equiv 1$,

$g_2 = \mathbf{e_1}$, $g_3(\mathbf{z}) = T(\mathbf{z};0,\psi(x))$, $g_4 = g_5 \equiv (1,1)^T$, such that $r_1 = 3$, $r_2 = r_4 = r_5 = 0$, and $r_3 = 1$, we get that

$$\sum_{i_1,i_2=1}^{n} \left(\mathbf{a_{i_1,i_2}}(x)\right)_1 \left(\mathbf{a_{i_1,i_2}}(x)\right)_3$$
$$= \int_{[-1,1]^2} \langle \nabla K(z_1,z_2), (z_2,-z_1)^T \rangle \times K(z_1,z_2)\,dz_1 dz_2 + O(h) + O\left((nh)^{-1}\right)$$
$$= O(h) + O\left((nh)^{-1}\right),$$

where the last equation is true, since $x \mapsto K_2^2(x)x$ and $x \mapsto K_2^{(1)}(x)K_2(x)$ are both odd functions.

*Step 7:*

We check the Lindeberg condition: For any $\varepsilon > 0$

$$\sum_{i_1,i_2=1}^{n} \|\mathbf{a_{i_1,i_2}}(x)\|^2 \, \mathrm{E}\left(\varepsilon_{i_1,i_2}^2 \mathbf{1}_{\left\{\|\mathbf{a_{i_1,i_2}}(x)\| \, |\varepsilon_{i_1,i_2}| > \varepsilon\right\}}\right) = o(1).$$

where we assume without loss of generality that $\|\cdot\|$ is the Euclidean norm on $\mathbb{R}^3$. Notice that

$$\max_{1 \leq i_1,i_2 \leq n} \|\mathbf{a_{i_1,i_2}}(x)\|^2 \leq \max_{1 \leq i_1,i_2 \leq n} \left(\mathbf{a_{i_1,i_2}}(x)\right)_1^2 + \max_{1 \leq i_1,i_2 \leq n} \left(\mathbf{a_{i_1,i_2}}(x)\right)_2^2 + \max_{1 \leq i_1,i_2 \leq n} \left(\mathbf{a_{i_1,i_2}}(x)\right)_3^2$$
$$\leq 2\|\nabla K\|_\infty^2 \left(\left(n^2 h^2\right)^{-1} + \left(n^2 h^4\right)^{-1}\right) + \|K\|_\infty^2 \left(n^2 h^2\right)^{-1} = o(1),$$
(5)

and by Step 1 – 3 for any $x \in I$,

$$\sum_{i_1,i_2=1}^{n} \|\mathbf{a_{i_1,i_2}}(x)\|^2 = \sum_{i_1,i_2=1}^{n} \left(\mathbf{a_{i_1,i_2}}(x)\right)_1^2 + \sum_{i_1,i_2=1}^{n} \left(\mathbf{a_{i_1,i_2}}(x)\right)_2^2 + \sum_{i_1,i_2=1}^{n} \left(\mathbf{a_{i_1,i_2}}(x)\right)_3^2$$
$$\to V_\phi^S(x) + V_\psi^S + V_\tau^S.$$
(6)

Thus, by (5) and (6)

$$\sum_{i_1,i_2=1}^{n} \|\mathbf{a_{i_1,i_2}}(x)\|^2 \, \mathrm{E}\left(\varepsilon_{i_1,i_2}^2 \mathbf{1}_{\{\|\mathbf{a_{i_1,i_2}}(x)\| \, |\varepsilon_{i_1,i_2}|>\varepsilon\}}\right) \leq \max_{1\leq i_1,i_2\leq n} \mathrm{E}\left(\varepsilon_{i_1,i_2}^2 \mathbf{1}_{\{\|\mathbf{a_{i_1,i_2}}(x)\| \, |\varepsilon_{i_1,i_2}|>\varepsilon\}}\right) \sum_{i_1,i_2=1}^{n} \|\mathbf{a_{i_1,i_2}}(x)\|^2$$
$$= o(1).$$

Eventually, the Lindeberg-Feller theorem concludes the lemma. $\square$

## B. Supplement: Proofs for uniform confidence bands

*B.1. Rate of convergence of the Hessian matrix: proof of Lemma 6*

**Lemma 6.** *Let B be a compact subset of $\mathbb{R}$ and $\hat{A}_n, \hat{A} : B \to \mathbb{R}^{2\times 2}$ be matrix valued (random) functions with (stochastically) bounded marginals and*

$$\|\hat{A}_n - \hat{A}\|_\infty = O_P(r_n),$$

*for some real-valued sequence $r_n$. If $\hat{A}_n^{-1}(x)$ and $\hat{A}^{-1}(x)$ exist almost surely for every $x \in B$, then*

$$\|\hat{A}_n^{-1} - \hat{A}^{-1}\|_\infty = O_P(r_n).$$

*Proof.* First, note that for any $x \in B$ one has almost surely

$$\|\hat{A}_n^{-1}(x) - \hat{A}^{-1}(x)\| \leq \|\hat{A}^{-1}(x)\| \cdot \|\hat{A}_n(x) - \hat{A}(x)\| \cdot \|\hat{A}_n^{-1}(x)\|.$$

Since $\hat{A}_n$ and $A$ are stochastically bounded the inverse functions $\hat{A}_n^{-1}, A^{-1}$ are almost surely stochastically bounded. Therefore,
$$||\hat{A}_n^{-1} - \hat{A}^{-1}||_\infty \leq ||\hat{A}^{-1}||_\infty ||\hat{A}_n - \hat{A}||_\infty ||\hat{A}_n^{-1}||_\infty = O_P(r_n).$$

□

*Proof of Lemma 6.* By compactness of $I$ and Lipschitz continuity of cosine and sine it easily follows that
$$||\cos^2 \circ \hat{\psi}_n - \cos^2 \circ \psi||_\infty = O(||\hat{\psi}_n - \psi||_\infty), \quad ||\sin^2 \circ \hat{\psi}_n - \sin^2 \circ \psi||_\infty = O(||\hat{\psi}_n - \psi||_\infty),$$
so that by means of Proposition 1
$$||V_{\hat{\phi}}^H - V_{\phi}^H||_\infty = K_2^{(1)}(0) \sup_{x \in I} |\hat{\tau}_n(x) \cos^2(\hat{\psi}_n(x)) - \tau(x) \cos^2(\psi(x))|$$
$$\leq O(||\hat{\psi}_n - \psi||_\infty) + O(||\hat{\tau}_n - \tau||_\infty)$$
$$= O(||\hat{\psi}_n - \psi||_\infty) + O(||\hat{w}_n||_\infty) + O_P(h) + O((nh)^{-1}).$$

By Proposition 3 the preceding display is $O_P(h) + O_P(\sqrt{\log(n)}/nh)$. Similarly,
$$||V_{\hat{\phi}}^S - V_{\phi}^S||_\infty = O_P(\sqrt{\log(n)}/nh), \quad ||V_{\hat{\psi}}^H - V_{\psi}^H||_\infty = O_P(h) + O_P(\sqrt{\log(n)}/nh).$$

In the proof of Lemma 4 we have shown that for $\tilde{\mathbf{H}}_n(x) = \nabla \nabla^T \hat{\mathbb{M}}_n(0, \psi(x); x)$
$$||\tilde{\mathbf{H}}_n - \mathbf{H}||_\infty = O_P(h) + O_P((nh)^{-1}).$$

In addition, by the stochastic Lipschitz continuity of $\nabla \nabla^T \hat{\mathbb{M}}_n$ obtain
$$||\hat{\mathbf{H}}_n - \tilde{\mathbf{H}}_n||_\infty \leq O_P(1) ||(\hat{w}_n(\cdot), (\hat{\psi}_n(x) - \psi)(\cdot))^T||_\infty = O_P(\sqrt{\log(n)}/nh),$$
by means of Theorem 1. Combine the last two displays to get $||\hat{\mathbf{H}}_n - \mathbf{H}_n||_\infty = O_P(\sqrt{\log(n)}/nh)$. Application of Lemma 6 completes the proof. □

### B.2. Gaussian approximation and quantile approximation: proof of Lemma 7

For the proof we will need the following lemma.

**Lemma 7.** *Let $(X_n)_{n \in \mathbb{N}}$ and $(Y_n)_{n \in \mathbb{N}}$ be stochastic processes and each $Y_n$ has a continuous distribution. Assume that $(a_n)_{n \in \mathbb{N}}$ and $(b_n)_{n \in \mathbb{N}}$ are sequences of positive real numbers with*

1. *$b_n = o(1)$ and $a_n/b_n = o(1)$;*
2. *$|X_n - Y_n| = O_P(a_n)$ $n \to \infty$;*
3. *There exists a $\zeta_0 > 0$ with $\limsup_{n \to \infty} \mathbb{L}(Y_n, \zeta b_n) \leq \delta(\zeta)$ for any $\zeta \in (0, \delta(\zeta))$, where $\delta(\zeta) > 0$ and $\lim_{\zeta \to 0} \delta(\zeta) = 0$.*

*Then,*
$$\lim_{n \to \infty} P(X_n \leq q_\alpha(Y_n)) = \alpha.$$

*Proof of Lemma 7.* First, it holds that
$$\limsup_n P(X_n \leq q_\alpha(Y_n)) = \alpha.$$

16

Indeed, let $\zeta \in (0, \zeta_0)$, then

$$P(X_n \leq q_\alpha(Y_n)) \leq P(Y_n \leq q_\alpha(Y_n) + b_n\zeta) + P(|Y_n - X_n| \geq b_n\zeta)$$
$$\leq \alpha + P(q_\alpha(Y_n) \leq Y_n \leq q_\alpha(Y_n) + b_n\zeta) + P(|Y_n - X_n| \geq b_n\zeta)$$
$$\leq \alpha + \mathbb{L}(Y_n, b_n\zeta) + P(|Y_n - X_n| \geq b_n\zeta).$$

By considering $n \to \infty$ and then $\zeta \searrow 0$ yields the assertion above. Next,

$$\liminf_n P(X_n \leq q_\alpha(Y_n)) = \alpha,$$

which would complete the proof. Analogously as just,

$$P(X_n \leq q_\alpha(Y_n)) \geq P(Y_n \leq q_\alpha(Y_n) - b_n\zeta) - P(|Y_n - X_n| \geq b_n\zeta)$$
$$\geq \alpha - P(q_\alpha(Y_n) - b_n\zeta \leq Y_n \leq q_\alpha(Y_n)) - P(|Y_n - X_n| \geq b_n\zeta)$$
$$\geq \alpha - \mathbb{L}(Y_n, b_n\zeta) - P(|Y_n - X_n| \geq b_n\zeta).$$

Again, considering $n \to \infty$ and then $\zeta \searrow 0$ leads to the claim. □

*Proof of Lemma 7. Ad (i).*

On the one hand, by means of Lemma 8 with

$$g_1(\mathbf{z}) \equiv (V_\phi^S(x))^{-1/2}, \quad g_2 \equiv \mathbf{e}_1, \quad g_3(\mathbf{z}) \equiv (\sin \psi, \cos \psi)^T, \quad f(\mathbf{z}) = [\nabla K(\mathbf{z}) \mid \nabla K(\mathbf{z})],$$

and $w = 0$ and $\psi = \psi(x)$, such that $r_1 = 1$ and $r_2 = r_3 = 0$, gives us that for $n$ large enough

$$\|Z_{n,\phi}^\varepsilon - Z_{n,\phi}^W\|_\infty = O_P\Big(\frac{\sqrt{\log(n)}}{n^{1/2}h}\Big), \tag{7}$$

where

$$Z_{n,\phi}^W(x) = (V_\phi^S(x))^{-1/2} \int_{\mathbb{R}^2} \langle (\nabla K)\left(\mathcal{R}_{-\psi(x)}(\mathbf{p}(x)/h - \mathbf{z})\right), (\sin \psi(x), \cos \psi(x))^T \rangle dW(\mathbf{z}),$$

and $W$ is a suitable Wiener sheet on $\mathbb{R}^2$. On the other hand, using Lemma 8 with $\xi_{i_1,i_2}$ instead of $\varepsilon_{i_1,i_2}$ (this is fine, since $\xi_{i_1,i_2}$ has the same independence properties as $\varepsilon_{i_1,i_2}$ but stricter moment properties) and setting in that context $\sigma = 1$, $w = \hat{w}_n(x) = (\hat{\phi}_n(x) - \phi(x))/h$ and $\psi = \hat{\psi}_n(x)$, as well as

$$g_1(\mathbf{z}) \equiv (V_\phi^S(x))^{-1/2}, g_2 \equiv \mathbf{e}_1, g_3(\mathbf{z}) \equiv (\sin \psi, \cos \psi)^T f(\mathbf{z}) = [\nabla K(\mathbf{z}) \mid \nabla K(\mathbf{z})], \tag{8}$$

such that $r_1 = 1$ and $r_2 = r_3 = 0$, we obtain

$$\|Z_{n,\hat\phi}^\xi - Z_{n,\hat\phi}^W\|_\infty = O_P\Big(\frac{\sqrt{\log(n)}}{n^{1/2}h}\Big), \tag{9}$$

where

$$Z_{n,\hat\phi}^W(x) = \frac{1}{\sqrt{V_\phi^S(x)}} \int_{\mathbb{R}^2} \langle (\nabla K)\left(\mathcal{R}_{-\hat\psi_n(x)}(h^{-1}(x,\hat\phi_n(x))^T - \mathbf{z})\right), (\sin \hat\psi_n(x), \cos \hat\psi_n(x))^T \rangle dW(\mathbf{z}).$$

Define the function

$$F(\mathbf{z}; w, \psi) = g_1(h\mathbf{z})\langle f(\mathcal{R}_{-\psi}(\mathbf{p}(x)/h + w\mathbf{e}_2 - \mathbf{z}))g_2(h\mathbf{z}), g_3(h\mathbf{z})\rangle$$

and note that we suppressed the dependency of $g_1$ and $g_3$ on $\psi$ in the notation. With this, note that

$$Z_{n,\hat\phi}^W(x) = \int_{\mathbb{R}^2} F(\mathbf{z}; \hat{w}_n(x), \hat\psi_n(x)) \, dW(\mathbf{z}), \quad \text{and} \quad Z_{n,\phi}^W(x) = \int_{\mathbb{R}^2} F(\mathbf{z}; 0, \psi(x)) \, dW(\mathbf{z}).$$

17

Choosing $\delta = O_P(\sqrt{\log(n)}/nh)$, then the third part of Lemma 9 together with Theorem 1 implies
$$||Z^W_{n,\widehat{\phi}} - Z^W_{n,\phi}||_\infty = O_P(\log(n)/nh^2).$$

With this, (7) and (9) we have that
$$||Z^\xi_{n,\widehat{\phi}} - Z^\varepsilon_{n,\phi}||_\infty \leq ||Z^\xi_{n,\widehat{\phi}} - Z^W_{n,\widehat{\phi}}||_\infty + ||Z^W_{n,\widehat{\phi}} - Z^W_{n,\phi}||_\infty + ||Z^W_{n,\phi} - Z^\varepsilon_{n,\phi}||_\infty = O_P\Big(\frac{\sqrt{\log(n)}}{n^{1/2}h}\Big)$$

which shows the first part of the lemma.

*Ad (ii).*

Define $\mathcal{Z}^W_{n,\widehat{\phi}} = ||Z^W_{n,\widehat{\phi}}||_\infty$. The second part of Lemma 9 with the same $g_1, g_2, g_3$ and $f$ as in (8) shows that $\mathrm{E}(\mathcal{Z}^W_{n,\widehat{\phi}}) \leq C\sqrt{\log(n)}$, for some constant $C > 0$. Thus, Markov's inequality implies
$$\mathcal{Z}^W_{n,\widehat{\phi}} = O_P(\sqrt{\log(n)}),$$

showing the first claim in (ii) of this lemma. Next, let $\delta_n = o(\log(n)^{-1/2})$, then using Theorem 2.1 in Chernozhukov et al. [2] deduce that
$$\mathbb{L}\Big(\mathcal{Z}^W_{n,\widehat{\phi}}, \delta_n\Big) \leq 4\delta_n\Big(\mathrm{E}(\mathcal{Z}^W_{n,\widehat{\phi}}) + 1\Big) = o(1), \tag{10}$$

where we used for the last step the bound on the expected value. Next, for any $\zeta > 0$ and some real-valued sequence $b_n$,
$$\begin{aligned}\mathbb{L}(\mathcal{Z}^\xi_{n,\widehat{\phi}}, \zeta b_n) &\leq P\big(|\mathcal{Z}^\xi_{n,\widehat{\phi}} - \mathcal{Z}^W_{n,\widehat{\phi}}| > \zeta b_n\big) + \sup_{x \in I} P\big(|\mathcal{Z}^\xi_{n,\widehat{\phi}} - x| \leq \zeta b_n, |\mathcal{Z}^\xi_{n,\widehat{\phi}} - \mathcal{Z}^W_{n,\widehat{\phi}}| \leq \zeta b_n\big) \\ &\leq P\big(|\mathcal{Z}^\xi_{n,\widehat{\phi}} - \mathcal{Z}^W_{n,\widehat{\phi}}| > \zeta b_n\big) + \mathbb{L}(\mathcal{Z}^W_{n,\widehat{\phi}}, 2\zeta b_n).\end{aligned} \tag{11}$$

In view of (10), if $b_n = o(\log(n)^{-1/2})$ the Lévy concentration function in the preceding display is $o(1)$. Note that (9) implies by triangle inequality
$$||\mathcal{Z}^\xi_{n,\widehat{\phi}} - \mathcal{Z}^W_{n,\widehat{\phi}}||_\infty = O_P\Big(\frac{\sqrt{\log(n)}}{n^{1/2}h}\Big), \tag{12}$$

such that the first summand on the right-hand side of (11) can be made arbitrary small provided that $\frac{\log(n)}{n^{1/2}h} = o(1)$, which is implied by Assumption 3. Hence, if $b_n = o(\log(n)^{-1/2})$ holds, we obtain
$$\lim_{\zeta \to 0} \limsup_{n \to \infty} \mathbb{L}(\mathcal{Z}^\xi_{n,\widehat{\phi}}, \zeta b_n) = 0, \tag{13}$$

which shows the second claim of the lemma.

*Ad (iii).*

Lemma 7 shows $q_{1-\alpha}(\mathcal{Z}^W_{n,\widehat{\phi}}) \cong q_{1-\alpha}(\mathcal{Z}^\xi_{n,\widehat{\phi}})$, since (12) and (13) provide the assumptions of the lemma. Finally, the fourth part of Lemma 9 states that $q_{1-\alpha}(\mathcal{Z}^W_{n,\widehat{\phi}}) \cong \sqrt{\log(n)}$ which concludes the proof. $\square$

## C. Supplement: Auxiliary results for asymptotically handling various terms

Recall the definition of $j_\tau(\mathbf{z}) = \tau(z_1) 1_{[0,\phi(z_1)]}(z_2)$ for $\mathbf{z} = (z_1, z_2)^T \in [-1, 1]^2$ and $\mathbf{p}(x) = (x, \phi(x))^T$. In addition, for sake of brevity we write
$$\mathbf{p}_{w,h}(x) = \mathbf{p}(x) + hw\mathbf{e}_2 = (x, \phi(x) + hw)^T.$$

For all following statements in this section we assume that $(x, w, \psi)^T \in \Theta_n$, as stipulated in the notation.

Let $f : \mathbb{R}^2 \to \mathbb{R}^{2 \times 2}$ be a function with $f \in C^1(\mathbb{R}^2)$ having compact support in $[-1, 1]^2$ and bounded marginals. Moreover, let $g_1 : \mathbb{R}^2 \to \mathbb{R}$ be a function in $C^1(\mathbb{R}^2)$ which does not depend on the bandwidth $h$ and $g_2, g_3, g_4, g_5 : \mathbb{R}^2 \to \mathbb{R}^2$ be functions which have the following shape: For $i = 2, 3, 4, 5$ let either $g_i(\mathbf{z}) = \mathbf{C}_{m,g_i}(\mathbf{p}_{w,h}(x) - \mathbf{z})$, for some constant matrix $\mathbf{C}_{m,g_i} \in \mathbb{R}^{2 \times 2}$, or $g_i(\mathbf{z}) = \mathbf{C}_{g_i}$, for some constant vector $\mathbf{C}_{g_i} \in \mathbb{R}^2$. Both $\mathbf{C}_{m,g_i}$ and $\mathbf{C}_{g_i}$ can depend on $w$ or $\psi$. In addition, let $r_2, r_3, r_4, r_5 \in \{0, 1\}$ be such that $r_i = 0$ if and only if $g_i(\mathbf{z}) \equiv \mathbf{C}_{g_i}$, and $r_i = 1$ if and only if $g_i(\mathbf{z}) = \mathbf{C}_{m,g_i}(\mathbf{p}_{w,h}(x) - \mathbf{z})$, for $i = 2, 3, 4, 5$ respectively. Moreover, for $i = 2, 3, 4, 5$ let $\tilde{g}_i : \mathbb{R}^2 \to \mathbb{R}^2$ be defined by

$$\tilde{g}_i(\mathbf{z}) = \begin{cases} g_i(\mathbf{z}), & \text{if } g_i(\mathbf{z}) \equiv \mathbf{C}_{g_i}, \text{ or} \\ (-1)^{r_i} \mathbf{C}_{m,g_i}(\mathcal{R}_\psi \mathbf{z}), & \text{else.} \end{cases}$$

Note that these functions do not dependent on $h$ and

$$g_i(\mathbf{p}_{w,h}(x) - h\mathcal{R}_\psi \mathbf{z}) = h^{r_i} \tilde{g}_i(\mathbf{z}). \tag{14}$$

Without loss of generality assume that $n$ is large enough such that $h^{-1} > 1$.

Section C.1 and Section C.2 summarizes the assertions of the auxiliary results and the remaining sections provide their proofs if cumbersome.

*C.1. Properties of Gaussian processes*

Let $r_1 = 1 + r_2 + r_3$. Set for $(x, w, \psi)^T \in \Theta_n$

$$Z_n^\varepsilon(x, w, \psi) = \frac{1}{nh^{r_1}\sigma} \sum_{i_1, i_2 = 1}^n \varepsilon_{i_1, i_2} g_1(\mathbf{x}_{\mathbf{i_1}, \mathbf{i_2}}) \langle f(h^{-1}\mathcal{R}_{-\psi}(\mathbf{p}_{w,h}(x) - \mathbf{x}_{\mathbf{i_1}, \mathbf{i_2}})) g_2(\mathbf{x}_{\mathbf{i_1}, \mathbf{i_2}}), g_3(\mathbf{x}_{\mathbf{i_1}, \mathbf{i_2}}) \rangle \tag{15}$$

and

$$Z_n^W(x, w, \psi) := \frac{1}{h^{r_1 - 1}} \int_{\mathbb{R}^2} g_1(h\mathbf{z}) \langle f(\mathcal{R}_{-\psi}(\mathbf{p}(x)/h + w\mathbf{e}_2 - \mathbf{z})) g_2(h\mathbf{z}), g_3(h\mathbf{z}) \rangle \, dW(\mathbf{z}), \tag{16}$$

where $W$ is a Wiener sheet on $\mathbb{R}^2$.

**Gaussian approximation**

**Lemma 8.** *Under the Assumption 1 and Assumption 3, on an appropriate probability space there exists a Wiener Sheet $W$ on $\mathbb{R}^2$ for which*

$$\sup_{(x, w, \psi)^T \in \Theta_n} |Z_n^\varepsilon(x, w, \psi) - Z_n^W(x, w, \psi)| = O_P\left(\frac{\sqrt{\log(n)}}{n^{1/2}h}\right).$$

*In particular,*

$$\sup_{(x, w, \psi)^T \in \Theta_n} |Z_n^\varepsilon(x, w, \psi)| - \sup_{(x, w, \psi)^T \in \Theta_n} |Z_n^W(x, w, \psi)| = O_P\left(\frac{\sqrt{\log(n)}}{n^{1/2}h}\right).$$

The proof is provided in Section C.3.

**Moments and quantiles of Gaussian processes**

**Lemma 9.** *Suppose Assumption 3 is satisfied. Then the following holds for any Wiener Sheet $W$ on $\mathbb{R}^2$.*

19

1. For any $\zeta > C_0$, where $C_0 > 0$ is a finite constant uniform for $x, w$ and $\psi$, it holds

$$P\Big( \sup_{(x,w,\psi) \in \Theta_n} |Z_n^W(x,w,\psi)| > \zeta \Big) \leq C_1 \frac{\lambda_3(\Theta_n)\zeta}{h^3} \exp(-C_2 \zeta^2),$$

   where $C_1, C_2$ are finite absolute constants uniform for $x, w$ and $\psi$.

2. There exists a finite constant $C_3 > 0$ uniform for $x, w$ and $\psi$ such that for sufficiently small $h$ (sufficiently large $n$) holds

$$\mathrm{E} \sup_{(x,w,\psi) \in \Theta_n} |Z_n^W(x,w,\psi)| \leq C_3 \sqrt{\log(n)}.$$

3. There exists a finite constant $C_4 > 0$ uniform for $x, w$ and $\psi$ such that for any $\delta \in (0,1)$

$$\mathrm{E} \sup_{\theta_1, \theta_2 \in \Theta_n \,:\, \|\theta_1 - \theta_2\| \leq \delta} |Z_n^W(\theta_1) - Z_n^W(\theta_2)| \leq C_4 \, \delta \, h^{-1} \sqrt{\log(n)}.$$

4. For any $\alpha \in (0,1)$

$$q_\alpha \Big( \sup_{(x,w,\psi) \in \Theta_n} |Z_n^W(x,w,\psi)| \Big) \cong \sqrt{\log(n)}.$$

The proof is deferred to Section C.4.

*C.2. Asymptotics of components of the contrast function and their derivatives*

In order to control the smooth part of the image we will often use the following lemma, which is proven in Section C.5.

**Lemma 10.** *Let $r_1 = 2 + j(r_2 + r_3)$ and $j \in \{1,2\}$, then under Assumptions 2 – 4 it holds that*

$$S_n(x; w, \psi) := (n^2 h^{r_1})^{-1} \sum_{i_1, i_2 = 1}^{n} g_1(\mathbf{x_{i_1,i_2}}) \langle f\left(h^{-1} \mathcal{R}_{-\psi}\left(\mathbf{p}_{w,h}(x) - \mathbf{x_{i_1,i_2}}\right)\right) g_2(\mathbf{x_{i_1,i_2}}), g_3(\mathbf{x_{i_1,i_2}}) \rangle^j$$

$$= \int_{[-1,1]^2} g_1(\mathbf{p}_{w,h}(x) - h\mathcal{R}_\psi \mathbf{z}) \langle f(\mathbf{z}) \tilde{g}_2(\mathbf{z}), \tilde{g}_3(\mathbf{z}) \rangle^j \, d\mathbf{z} + O\big((nh^{r_1 - 1 - r_2 - r_3})^{-1}\big)$$

$$= g_1(\mathbf{p}_{w,h}(x)) \int_{[-1,1]^2} \langle f(\mathbf{z}) \tilde{g}_2(\mathbf{z}), \tilde{g}_3(\mathbf{z}) \rangle^j \, d\mathbf{z} + O(h) + O\big((nh^{r_1 - 1 - r_2 - r_3})^{-1}\big),$$

*uniformly for $x, w$ and $\psi$.*

*More generally, let $r_1 = 2 + r_2 + r_3 + r_4 + r_5$, and let $f_1$ resp. $f_2$ be functions with same properties as $f$. Then one has*

$$\tilde{S}_n(x; w, \psi) := (n^2 h^{r_1})^{-1} \sum_{i_1, i_2 = 1}^{n} g_1\left(\mathbf{x_{i_1,i_2}}\right) \langle f_1(h^{-1} \mathcal{R}_{-\psi}(\mathbf{p}_{w,h}(x) - \mathbf{x_{i_1,i_2}})) g_2(\mathbf{x_{i_1,i_2}}), g_3\left(\mathbf{x_{i_1,i_2}}\right) \rangle$$

$$\times \quad \langle f_2(h^{-1} \mathcal{R}_{-\psi}(\mathbf{p}_{w,h}(x) - \mathbf{x_{i_1,i_2}})) g_4(\mathbf{x_{i_1,i_2}}), g_5\left(\mathbf{x_{i_1,i_2}}\right) \rangle$$

$$= g_1\left(\mathbf{p}_{w,h}(x)\right) \int_{[-1,1]^2} \langle f_1(\mathbf{z}) \tilde{g}_2(\mathbf{z}), \tilde{g}_3(\mathbf{z}) \rangle \langle f_2(\mathbf{z}) \tilde{g}_4(\mathbf{z}), \tilde{g}_5(\mathbf{z}) \rangle \, d\mathbf{z} + O(h)$$

$$+ O\big((nh^{r_1 - 1 - r_2 - r_3})^{-1}\big) + O\big((nh^{r_1 - 1 - r_4 - r_5})^{-1}\big),$$

*uniformly for $x, w$ and $\psi$.*

To take care of the jump part of the image we need the following lemma, which proof is deferred to Section C.5.

20

**Lemma 11.** *Let $r_1 = 2 + r_2 + r_3$. Suppose 2 – 4 are satisfied, then it holds that*

$$J_n(x; w, \psi) := (n^2 h^{r_1})^{-1} \sum_{i_1, i_2=1}^{n} j_\tau(\mathbf{x}_{i_1, i_2}) \langle f(h^{-1} \mathcal{R}_{-\psi}(\mathbf{p}_{w,h}(x) - \mathbf{x}_{i_1, i_2})) g_2(\mathbf{x}_{i_1, i_2}), g_3(\mathbf{x}_{i_1, i_2}) \rangle$$

$$= \int_{h^{-1} \mathcal{R}_{-\psi}(\mathbf{p}(x) - [0,1]^2 \setminus epi(\phi)) + w(\sin(\psi), \cos(\psi))^T} \tau(x - (h \mathcal{R}_\psi \mathbf{z})_1) \langle f(\mathbf{z}) \tilde{g}_2(\mathbf{z}), \tilde{g}_3(\mathbf{z}) \rangle \, d\mathbf{z}$$

$$+ O((nh)^{-1})$$

$$= \tau(x) \int_{\mathbb{H}(\psi(x) - \psi) + w(\sin(\psi), \cos(\psi))^T} \langle f(\mathbf{z}) \tilde{g}_2(\mathbf{z}), \tilde{g}_3(\mathbf{z}) \rangle \, d\mathbf{z} + O(h) + O((nh)^{-1}),$$

*uniformly for $x, w$ and $\psi$, where $\mathbb{H}(\psi) = \mathcal{R}_\psi(\mathbb{R} \times [0, \infty))$, as defined in (1).*

For handling the error terms of the empirical contrast function, we formulate the following result.

**Lemma 12.** *Let $r_1 = 2 + r_2 + r_3$ then under the Assumptions 1 – 3 one obtains*

$$E_n(x; w, \psi) := (n^2 h^{r_1})^{-1} \sum_{i_1, i_2=1}^{n} \varepsilon_{i_1, i_2} \langle f(h^{-1} \mathcal{R}_{-\psi}(\mathbf{p}_{w,h}(x) - \mathbf{x}_{i_1, i_2})) g_2(\mathbf{x}_{i_1, i_2}), g_3(\mathbf{x}_{i_1, i_2}) \rangle = O_P(\sqrt{\log(n)} (nh)^{-1}),$$

*uniformly for $(x, w, \psi)$ over $\Theta_n$.*

*Proof.* Note that $nh/\sigma E_n(x; w, \psi) = Z_n^\varepsilon(x, w, \psi)$ and $\lambda_3(\Theta_n) = O(h^{-1})$. Moreover, let $\zeta > 0$ be a large constant specified below, then $\zeta \sqrt{\log(n)}/2\sigma > C_0$ for sufficiently large $n$. By Lemma 8 there exists a Wiener sheet $W$ on $\mathbb{R}^2$ such that

$$R_n = \sup_{(x, w, \psi)^T \in \Theta_n} |Z_n^W(x, w, \psi) - Z_n^\varepsilon(x, w, \psi)| = O_P(\sqrt{\log(n)}/n^{1/2} h).$$

Part one of Lemma 9 implies for some absolute constants $C_3, C_4$ that

$$P\left( \sup_{(x, w, \psi)^T \in \Theta_n} |E_n(x, w, \psi)| > \sqrt{\log(n)} (nh)^{-1} \zeta \right)$$

$$\leq P\left( \sup_{(x, w, \psi)^T \in \Theta_n} |Z_n^W(x, w, \psi)| > \zeta \sqrt{\log(n)}/2\sigma \right) + P\left( R_n > \zeta \sqrt{\log(n)}/2\sigma \right)$$

$$\leq C_3 \frac{n}{h^3} \exp(-C_4 \zeta^2 \log(n)) + o(1) = o(1),$$

where the last equation holds for sufficiently large $\zeta > 0$. $\square$

### C.3. Proof of Lemma 8

**Lemma 13** (Integration by parts for Wiener sheet integrals). *Let $[a_1, b_1] \times [a_2, b_2] \subset \mathbb{R}_+^2$. For a continuously differentiable function $f : \mathbb{R}^2 \to \mathbb{R}$ and a standard Wiener sheet on $[0, \infty)^2$ it holds that*

$$\int_{[a_1, b_1] \times [a_2, b_2]} f(z_1, z_2) \, dW(z_1, z_2) = \int_{[a_1, b_1] \times [a_2, b_2]} W(z_1, z_2) f^{(1,1)}(z_1, z_2) \, dz_1 dz_2$$

$$- \int_{[a_1, b_1]} W(z_1, b_2) f^{(1,0)}(z_1, b_2) + W(z_1, a_2) f^{(1,0)}(z_1, a_2) \, dz_1$$

$$- \int_{[a_2, b_2]} W(b_1, z_2) f^{(0,1)}(b_1, z_2) \, dz_2 + W(a_1, z_2) f^{(0,1)}(a_1, z_2) \, dz_2$$

$$+ W(b_1, b_2) f(b_1, b_2) - W(a_1, b_2) f(a_1, b_2) - W(b_1, a_2) f(b_1, a_2)$$

$$+ W(a_1, a_2) f(a_1, a_2).$$

21

*Proof of Lemma 8.* For sake of brevity write for $\mathbf{z} \in \mathbb{R}^2$,

$$F(\mathbf{z}) = g_1(\mathbf{p}_{w,h}(x) + h\mathcal{R}_\psi \mathbf{z}) \langle f(\mathbf{z}) g_2(\mathbf{p}_{w,h}(x) + h\mathcal{R}_\psi \mathbf{z}), g_3(\mathbf{p}_{w,h}(x) + h\mathcal{R}_\psi \mathbf{z}) \rangle. \tag{17}$$

Note that $F$ is a real-valued, bounded and continuous differentiable function with the same compact support as $f$, that is $[-1,1]^2$. Moreover, $F$ may depend on $w, \psi$ and $h$, but we will suppress this in the notation. One has

$$||F||_\infty = O(h^{r_2+r_3}), \quad \text{and} \quad ||F^{(1,1)}||_\infty = O(h^{r_2+r_3}), \tag{18}$$

uniformly for $x, w$ and $\psi$, due to the shapes of $g_i$, $i = 1, 2, 3$. This uniformity comes from the fact that $x$ and $\psi$ are element of compact sets and $w$ is such that $\mathbf{p}_{w,h}(x)$ is element of a compact subset. Furthermore,

$$F((h^{-1}\mathcal{R}_{-\psi}(\mathbf{p}_{w,h}(x) - \mathbf{z}))) = g_1(\mathbf{z}) \langle f\left(h^{-1}\mathcal{R}_{-\psi}(\mathbf{p}_{w,h}(x) - \mathbf{z})\right) g_2(\mathbf{z}), g_3(\mathbf{z}) \rangle,$$

so that

$$Z_n^\varepsilon(x, w, \psi) = \frac{1}{nh^{r_1}\sigma} \sum_{i_1,i_2=1}^n \varepsilon_{i_1,i_2} F(h^{-1}\mathcal{R}_{-\psi}(\mathbf{p}_{w,h}(x) - \mathbf{x}_{\mathbf{i_1},\mathbf{i_2}})).$$

*Step 1: Change summation order*

Define $S_{i_1,i_2} = \sum_{l=1}^{i_1} \sum_{k=1}^{i_2} \varepsilon_{l,k}$ for $(i_1, i_2) \in \{0, 1, \ldots, n\}^2$ and set $S_{i_1,0} = S_{0,i_2} = 0$ for all $(i_1, i_2) \in \{0, 1, \ldots, n\}^2$. Thus,

$$\varepsilon_{i_1,i_2} = S_{i_1,i_2} - S_{i_1-1,i_2} - S_{i_1,i_2-1} + S_{i_1-1,i_2-1}$$

and we get

$$Z_n^\varepsilon(x, w, \psi) = \frac{1}{nh^{r_1}\sigma} \sum_{i_1,i_2=1}^n F(h^{-1}\mathcal{R}_{-\psi}(\mathbf{p}_{w,h}(x) - \mathbf{x}_{\mathbf{i_1},\mathbf{i_2}}))(S_{i_1,i_2} - S_{i_1-1,i_2} - S_{i_1,i_2-1} + S_{i_1-1,i_2-1}).$$

Rewrite the process $Z_{0,n,h}$ as

$$Z_n^\varepsilon(x, w, \psi) = \frac{1}{nh^{r_1}\sigma} \sum_{i_1,i_2=1}^{n-1} \left( F(h^{-1}\mathcal{R}_{-\psi}(\mathbf{p}_{w,h}(x) - \mathbf{x}_{i_1+1,i_2+1})) - F(h^{-1}\mathcal{R}_{-\psi}(\mathbf{p}_{w,h}(x) - \mathbf{x}_{i_1,i_2+1})) \right.$$
$$\left. - F(h^{-1}\mathcal{R}_{-\psi}(\mathbf{p}_{w,h}(x) - \mathbf{x}_{i_1+1,i_2})) + F(h^{-1}\mathcal{R}_{-\psi}(\mathbf{p}_{w,h}(x) - \mathbf{x}_{\mathbf{i_1},\mathbf{i_2}})) \right) S_{i_1,i_2} + \frac{R_{n,1}(x, w, \psi)}{nh^{r_1}\sigma},$$

where

$$R_{n,1}(x, w, \psi) = -\sum_{i_1=1}^{n-1} \left( F(h^{-1}\mathcal{R}_{-\psi}(\mathbf{p}_{w,h}(x) - \mathbf{x}_{i_1+1,n})) - F(h^{-1}\mathcal{R}_{-\psi}(\mathbf{p}_{w,h}(x) - \mathbf{x}_{i_1,n})) \right) S_{i_1,n}$$
$$- \sum_{i_2=1}^{n-1} \left( F(h^{-1}\mathcal{R}_{-\psi}(\mathbf{p}_{w,h}(x) - \mathbf{x}_{n,i_2+1})) - F(h^{-1}\mathcal{R}_{-\psi}(\mathbf{p}_{w,h}(x) - \mathbf{x}_{n,i_2})) \right) S_{n,i_2}$$
$$+ F(h^{-1}\mathcal{R}_{-\psi}(\mathbf{p}_{w,h}(x) - x_{n,n})) S_{n,n}.$$

Because $(x, w, \psi)^T \in \Theta_n$ it holds that $\phi(x) + hw \in [h, 1-h]$. Additionally, all occurring design points in $R_{n,1}(x, w, \psi)$ are (nearly) edge-points such that for sufficiently large $n$ all terms in $R_{n,1}(x, w, \psi)$ vanish as the arguments in the inverse image of $F$ leave the compact support of $F$. This vanishing property holds uniformly over $\Theta_n$ for any $n$

sufficiently large. Thus, for sufficiently large $n$

$$Z_n^\varepsilon(x,w,\psi) = \frac{1}{nh^{r_1}\sigma} \sum_{i_1,i_2=1}^{n-1} \Big( F(h^{-1}\mathcal{R}_{-\psi}(\mathbf{p}_{w,h}(x) - \mathbf{x}_{i_1+1,i_2+1})) - F(h^{-1}\mathcal{R}_{-\psi}(\mathbf{p}_{w,h}(x) - \mathbf{x}_{i_1,i_2+1}))$$
$$- F(h^{-1}\mathcal{R}_{-\psi}(\mathbf{p}_{w,h}(x) - \mathbf{x}_{i_1+1,i_2})) + F(h^{-1}\mathcal{R}_{-\psi}(\mathbf{p}_{w,h}(x) - \mathbf{x}_{i_1,i_2})) \Big) S_{i_1,i_2}.$$

*Step 2: Approximation by a Wiener Sheet*

Introduce the process

$$Z_{n,1}^W(x,w,\psi) = \frac{1}{nh^{r_1}} \sum_{i_1,i_2=1}^{n-1} \Big( F(h^{-1}\mathcal{R}_{-\psi}(\mathbf{p}_{w,h}(x) - \mathbf{x}_{i_1+1,i_2+1})) - F(h^{-1}\mathcal{R}_{-\psi}(\mathbf{p}_{w,h}(x) - \mathbf{x}_{i_1,i_2+1}))$$
$$- F(h^{-1}\mathcal{R}_{-\psi}(\mathbf{p}_{w,h}(x) - \mathbf{x}_{i_1+1,i_2})) + F(h^{-1}\mathcal{R}_{-\psi}(\mathbf{p}_{w,h}(x) - \mathbf{x}_{i_1,i_2})) \Big) W(i_1,i_2),$$

where $W$ is a Wiener Sheet with

$$\sup_{1 \leq i_1,i_2 \leq n} |S_{i_1,i_2}\sigma^{-1} - W(i_1,i_2)| = O_P\left(n^{\frac{2-\delta}{2}}\sqrt{\log(n)}\right). \tag{19}$$

Such a Wiener Sheet exists, see Theorem 1 in [4], provided there exists some $\delta \in (0,1]$ such that $E|\varepsilon_{1,1}|^k < \infty$ for $k > 4/(2-\delta)$. Note that the preconditions of this lemma allow to choose $\delta = 1$. Let

$$A_{i_1,i_2} = [x_{i_1}, x_{i_1+1}) \times [x_{i_2}, x_{i_2+1}), \quad i_1, i_2 \in \{1, \ldots, n-1\},$$

where the right boundary is included if $i_1 = n-1$ or $i_2 = n-1$. This yields

$$|Z_n^\varepsilon(x,w,\psi) - Z_{n,1}^W(x,w,\psi)|$$
$$\leq \sup_{1 \leq i_1,i_2 \leq n} \frac{|S_{i_1,i_2}\sigma^{-1} - W(i_1,i_2)|}{nh^{r_1}} \sum_{i_1,i_2=1}^{n-1} |F(h^{-1}\mathcal{R}_{-\psi}(\mathbf{p}_{w,h}(x) - \mathbf{x}_{i_1+1,i_2+1})) - F(h^{-1}\mathcal{R}_{-\psi}(\mathbf{p}_{w,h}(x) - \mathbf{x}_{i_1,i_2+1}))$$
$$- F(h^{-1}\mathcal{R}_{-\psi}(\mathbf{p}_{w,h}(x) - \mathbf{x}_{i_1+1,i_2})) + F(h^{-1}\mathcal{R}_{-\psi}(\mathbf{p}_{w,h}(x) - \mathbf{x}_{i_1,i_2}))|$$
$$\leq \sup_{1 \leq i_1,i_2 \leq n} \frac{|S_{i_1,i_2}\sigma^{-1} - W(i_1,i_2)|}{nh^{r_1}} \sum_{i_1,i_2=1}^{n-1} \int_{A_{i_1,i_2}} h^{-2}|F^{(1,1)}(h^{-1}\mathcal{R}_{-\psi}(\mathbf{p}_{w,h}(x) - \mathbf{z}))|\,d\mathbf{z}$$
$$\leq \sup_{1 \leq i_1,i_2 \leq n} \frac{|S_{i_1,i_2}\sigma^{-1} - W(i_1,i_2)|}{nh^{r_1}} \int_{[-1,1]^2} |F^{(1,1)}(\mathbf{z})|\,d\mathbf{z} = O_P\left(n^{\frac{-1}{2}} h^{r_2+r_3-r_1}\sqrt{\log(n)}\right)$$
$$= O_P\left(n^{\frac{-1}{2}} h^{-1} \sqrt{\log(n)}\right),$$

uniformly for $x, w$ and $\psi$, where we used for the second last equality (18) and (19), while the last equality is due to the choice of $r_1$ (see the line above (15)).

*Step 3: Continuous approximation*

Next, introduce the process

$$Z_{n,2}^W(x,w,\psi) = \frac{1}{h^{r_1}} \int_{[0,1]^2} F(h^{-1}\mathcal{R}_{-\psi}(\mathbf{p}_{w,h}(x) - \mathbf{z}))\,dW(\mathbf{z}).$$

23

Hence, by integration by parts for the Wiener sheet (Lemma 13),

$$Z_{n,2}^W(x, w, \psi) = \frac{1}{h^{r_1}} \int_{[0,1]^2} h^{-2} F^{(1,1)}(h^{-1}\mathcal{R}_{-\psi}(\mathbf{p}_{w,h}(x) - z))W(\mathbf{z})\,d\mathbf{z} + \frac{1}{h^{r_1}} R_{n,3}(x, w, \psi),$$

where

$$\begin{aligned}
R_{n,3}(x, w, \psi) = &-\int_{[0,1]} h^{-1} W(z_1, 1) F^{(1,0)}(h^{-1}\mathcal{R}_{-\psi}(\mathbf{p}_{w,h}(x) - (z_1, 1)^T))\,dz_1 \\
&-\int_{[0,1]} h^{-1} W(1, z_2) F^{(0,1)}(h^{-1}\mathcal{R}_{-\psi}(\mathbf{p}_{w,h}(x) - (1, z_2)^T))\,dz_2 \\
&+ W(1,1)) F(h^{-1}\mathcal{R}_{-\psi}(\mathbf{p}_{w,h}(x) - (1,1)^T)).
\end{aligned}$$

For sufficiently large $n$ one has that $R_{n,3}(x, w, \psi) \equiv 0$ as the support of the partial derivatives of $F$ and $F$ itself will be exceeded resp. deceeded. This holds uniformly over $\Theta_n$ for any $n$ large enough. Furthermore, obtain by the scaling properties of the Wiener sheet

$$\begin{aligned}
Z_{n,1}^W(x, w, \psi) &= \frac{1}{nh^{r_1}} \sum_{i_1, i_2=1}^{n-1} W(i_1, i_2) \int_{A_{i_1, i_2}} h^{-2}(F^{(1,1)} h^{-1}\mathcal{R}_{-\psi}(\mathbf{p}_{w,h}(x) - z))\,d\mathbf{z} \\
&\stackrel{d}{=} \frac{1}{h^{r_1}} \sum_{i_1, i_2=1}^{n-1} W(i_1/n, i_2/n) \int_{A_{i_1, i_2}} h^{-2}(F^{(1,1)} h^{-1}\mathcal{R}_{-\psi}(\mathbf{p}_{w,h}(x) - z))\,d\mathbf{z}
\end{aligned}$$

and hence for large enough $n$

$$\begin{aligned}
&|Z_{n,1}^W(x, w, \psi) - Z_{n,2}^W(x, w, \psi)| \\
&\stackrel{d}{=} \frac{1}{h^{r_1}} |\sum_{i_1, i_2=1}^{n-1} h^{-2} \int_{A_{i_1, i_2}} F^{(1,1)}(h^{-1}\mathcal{R}_{-\psi}(\mathbf{p}_{w,h}(x) - \mathbf{z}))[W(\mathbf{z}) - W(i_1/n, i_2/n)]\,d\mathbf{z}|.
\end{aligned}$$

A modulus of continuity for the Wiener sheet is given in Theorem 1 in [3], that is for any $\delta \in (0, 1)$

$$\sup_{\|\mathbf{z}_1 - \mathbf{z}_2\|_\infty < \delta} |W(\mathbf{z}_1) - W(\mathbf{z}_2)| = O_P\left(\sqrt{\log(1/\delta)\delta}\right). \tag{20}$$

In addition, for a finite constant $C > 0$ which is uniform in $x, w$ and $\psi$ we have by substitution and (18) that

$$\sum_{i_1, i_2=1}^{n-1} \int_{A_{i_1, i_2}} h^{-2} |F^{(1,1)}(h^{-1}\mathcal{R}_{-\psi}(\mathbf{p}_{w,h}(x) - \mathbf{z}))|\,d\mathbf{z} \leq \int_{[-1,1]^2} |F^{(1,1)}(\mathbf{z})|\,d\mathbf{z} < Ch^{r_2+r_3}$$

and therefore by incorporating the equidistant design in (20) one has

$$|Z_{n,1}^W(x, w, \psi) - Z_{n,2}^W(x, w, \psi)| = O_P\left(\sqrt{\log(n)} n^{-1/2} h^{r_2+r_3-r_1}\right) = O_P\left(\sqrt{\log(n)} n^{-1/2} h^{-1}\right),$$

where the $O$-term is uniform in $x, w$ and $\psi$.

*Step 4: Extension of the support*

Approximate $Z_{n,2}^W(x,w,\psi)$ by $Z_n^W(x,w,\psi)$ as in (16). For this purpose, obtain by a substitution and the scaling properties of the Wiener sheet that

$$Z_{n,2}^W(x,w,\psi) \stackrel{d}{=} \frac{1}{h^{r_1-1}} \int_{[0,\frac{1}{h}]^2} F(\mathcal{R}_{-\psi}(\mathbf{p}_{w,h}(x)/h - \mathbf{z})) \, dW(\mathbf{z})$$

as well as

$$F(\mathcal{R}_{-\psi}(\mathbf{p}_{w,h}(x)/h - \mathbf{z})) = g_1(h\mathbf{z})\langle f(\mathcal{R}_{-\psi}(\mathbf{p}(x)/h + w\mathbf{e}_2 - \mathbf{z}))g_2(h\mathbf{z}), g_3(h\mathbf{z})\rangle.$$

The latter two displays lead to

$$Z_n^W(x,w,\psi) - Z_{n,2}^W(x,w,\psi) \stackrel{d}{=} \frac{1}{h^{r_1-1}} \int_{A_n} F(\mathcal{R}_{-\psi}(\mathbf{p}_{w,h}(x)/h - z)) \, dW(\mathbf{z}) =: R_{n,4}(x,w,\psi),$$

where $A_n = \mathbb{R}^2 \setminus [0, h^{-1}]^2$. For $\mathbf{z} = (z_1, z_2)^T \in A_n$ one has that $z_2 < 0$ or $z_2 > h^{-1}$. Moreover, as $(x,w,\psi)^T \in \Theta_n$, it holds that $\phi(x)/h + w \in [1, h^{-1} - 1]$ Thus, $\phi(x)/h + w - z_2$ leaves the support of the function $F$ for large enough $n$. This yields $R_{n,4}(x,w,\psi) \equiv 0$ for sufficiently large $n$ uniformly over $\Theta_n$.

*Step 5: Conclusion*

Summarizing all approximations steps above, we have for sufficiently large $n$

$$\begin{aligned}
&Z_n^\varepsilon(x,w,\psi) - Z_n^W(x,w,\psi) \\
&\leq |Z_n^\varepsilon(x,w,\psi) - Z_{n,1}^W(x,w,\psi)| + Z_{n,1}^W(x,w,\psi) - Z_n^W(x,w,\psi) \\
&\leq O_P\Big(\frac{\sqrt{\log(n)}}{n^{1/2}h}\Big) + |Z_{n,1}^W(x,w,\psi) - Z_{n,2}^W(x,w,\psi)| + Z_{n,2}^W(x,w,\psi) - Z_n^W(x,w,\psi) \\
&\leq O_P\Big(\frac{\sqrt{\log(n)}}{n^{1/2}h}\Big) + O_P\Big(\frac{\sqrt{\log(n)}}{\sqrt{n}h}\Big) + R_{n,4}(x,w,\psi) \\
&= O_P\Big(\frac{\sqrt{\log(n)}}{n^{1/2}h}\Big),
\end{aligned}$$

uniformly for $x, w$ and $\psi$. This concludes the first assertion of the lemma. The second assertion follows by the triangle inequality. $\square$

*C.4. Proof of Lemma 9*

By means of substitution one easily yields the following result.

**Lemma 14.** *For any $\varepsilon \in (0,1)$, $\rho < 1$ and $\eta > 0$ there exists a finite constant $C > 0$ depending only on $\rho$ and $\eta$ such that*

$$\int_0^\varepsilon x^{-\rho} \log(x^{-1})^\eta \, dx \leq C |\log(\varepsilon^{-1})|^\eta \varepsilon^{1-\rho}.$$

**Lemma 15.** *It holds that*

$$q_\alpha\Big(\sup_{(x,w,\psi)\in\Theta_n} |Z_n^W(x,w,\psi)|\Big) \cong \mathrm{E}\Big(\sup_{(x,w,\psi)\in\Theta_n} |Z_n^W(x,w,\psi)|\Big).$$

25

*Proof.* For sake of brevity write $\mathcal{Z}_n = \sup_{(x,w,\psi) \in \Theta_n} |Z_n^W(x,w,\psi)|$. Without loss of generality assume that $\mathrm{E}(\mathcal{Z}_n^2) \leq 1$. Use Lemma B.1 in [2] to obtain for any $\alpha \in (0,1)$

$$q_\alpha(\mathcal{Z}_n) \leq C \mathrm{E}(\mathcal{Z}_n),$$

for some constant $C > 0$. Let $C_1 = \sqrt{2|\log(1/\alpha)|}$, then Borell's inequality (see Proposition A.2.1 in van der Vaart and Wellner [6]) implies

$$P(\mathcal{Z}_n \leq \mathrm{E}(\mathcal{Z}_n) - C_1) \leq P(|\mathcal{Z}_n - \mathrm{E}(\mathcal{Z}_n)| \geq C_1) \leq \exp(-C_1^2/2) = \alpha.$$

This yields $q_\alpha(\mathcal{Z}_n) \geq c \mathrm{E}(\mathcal{Z}_n)$ for some constant $c > 0$. $\square$

*Proof of Lemma 9. First part: Maximal inequality*

We intend to use Proposition A.2.7 in van der Vaart and Wellner [6]. For this purpose, define the following semi-metric on $\Theta_n$

$$\rho^2\big((x_1,w_1,\psi_1)^T, (x_2,w_2,\psi_2)^T\big) = \mathrm{E}|Z_n^W(x_1,w_1,\psi_1) - Z_n^W(x_2,w_2,\psi_2)|^2.$$

Note that this semi-metric depends on $n$ and $h$, which we suppress in the notation. Let $F$ be as in (17) and note that

$$F(\mathfrak{R}_{-\psi}(\mathbf{p}(x)/h + w\mathbf{e}_2 - \mathbf{z})) = g_1(h\mathbf{z}) \langle f(\mathfrak{R}_{-\psi}(\mathbf{p}(x)/h + w\mathbf{e}_2 - \mathbf{z})) g_2(h\mathbf{z}), g_3(h\mathbf{z}) \rangle,$$

such that $Z_n^W(x,w,\psi) = \int_{\mathbb{R}^2} F(\mathfrak{R}_{-\psi}(\mathbf{p}(x)/h + w\mathbf{e}_2 - \mathbf{z})) \, dW(\mathbf{z})$. Recall that $F$ has compact support since $\mathrm{supp}(f) = [-1,1]$ and in addition $F$ is Lipschitz continuous on its support with uniform Lipschitz constant $L_F = O(h^{r_2 + r_3})$, compare to (18). Furthermore, the function

$$G(x,w,\psi) := \mathfrak{R}_{-\psi}(\mathbf{p}(x)/h + w\mathbf{e}_2 - \mathbf{z})$$

is Lipschitz continuous on $\Theta_n$ with Lipschitz constant $L_G = O(h^{-1})$, which is also uniform in $\mathbf{z}$.

Hence, by Ito-Isometry and the Lipschitz continuity of $F$ and $G$,

$$\rho^2\big((x_1,w_1,\psi_1)^T, (x_2,w_2,\psi_2)^T\big)$$
$$= \frac{1}{h^{r_1-1}} \int_{A_1 \cup A_2} |F(\mathfrak{R}_{-\psi_1}(\mathbf{p}(x_1)/h + w_1\mathbf{e}_2 - \mathbf{z})) - F(\mathfrak{R}_{-\psi_2}(\mathbf{p}(x_2)/h + w_2\mathbf{e}_2 - \mathbf{z}))|^2 \, d\mathbf{z}$$
$$\leq \frac{1}{h^{r_1-1}} L_F^2 L_G^2 \|(x_1,w_1,\psi_1)^T - (x_2,w_2,\psi_2)^T\|^2 \lambda_2(A_1 \cup A_2),$$

where $A_i$ denotes the bounded set $\{\mathbf{z} \in \mathbb{R}^2 : F(\mathfrak{R}_{-\psi_i}(\mathbf{p}(x_i)/h + w_i\mathbf{e}_2 - \mathbf{z})) \neq 0\}$ for $i = 1,2$ respectively. Furthermore, the semi-metric is bounded. Indeed, it easily follows since by choice $1 - r_1 + r_2 + r_3 = 0$ and because $F$ is bounded that $\rho^2\big((x_1,w_1,\psi_1)^T, (x_2,w_2,\psi_2)^T\big) \leq 2h^{1-r_1} \|F\|_\infty^2 \lambda_2(A_1 \cup A_2) \leq 8C_F$, for some constant $C_F > 0$. Thus, for any $(x_1,w_1,\psi_1), (x_2,w_2,\psi_2) \in \Theta_n$,

$$\rho\big((x_1,w_1,\psi_1)^T, (x_2,w_2,\psi_2)^T\big) \leq \sqrt{8C_F} \wedge Ch^{-1} \|(x_1,w_1,\psi_1)^T - (x_2,w_2,\psi_2)^T\|, \tag{21}$$

for some constant $C > 0$ uniform for $x, w$ and $\psi$.

With this, the diameter of $\Theta_n$ with respect to $\rho$ is bounded by $\mathrm{diam}_\rho(\Theta_n) \leq \sqrt{8C_F}$ and in addition, the number of balls of radius $r > 0$ in the semi-metric $\rho$ that cover $\Theta_n$ is not larger than $Cr^{-3}h^{-1}\lambda_3(\Theta_n)$. Moreover,

$$\sup_{(x,w,\psi)^T \in \Theta_n} \mathrm{E}|Z_n^W(x,w,\psi)|^2 \leq C_0$$

for some appropriate constant $C_0 > 0$, which can be chosen uniformly in $x, w$ and $\psi$. Applying Proposition A.2.7 in van der Vaart and Wellner [6] provides the first part of the lemma.

*Second part: Order of the moment*

For the second part use the second statement of Corollary 2.2.8 in van der Vaart and Wellner [6] and Lemma 14 together with a substitution to see that

$$E\Big(\sup_{(x,w,\psi)\in\Theta_n} |Z_n^W(x,0,\psi(x))|\Big) \leq E|Z_n^W(0,0,\psi(x))| + C_1 \int_0^{\sqrt{8}C_F} \sqrt{\log(Cy^{-1}h^{-1})}\,dy$$

$$\leq C_2 + C_3\sqrt{\log(n)},$$

where $C_i > 0$ are appropriate constants for $i = 1,2,3$, which can be chosen uniformly for $x, w$ and $\psi$. Note that for the last line of the preceding display we used that $\sqrt{\log(n)} \cong \sqrt{\log(h^{-1})}$, which is implied by the preconditions of this lemma, i.e. by Assumption 3.

*Third part: Order of the increments*

From (21) derive that

$$\rho\big((x_1,w_1,\psi_1)^T,(x_2,w_2,\psi_2)^T\big) \leq Ch^{-1}\|(x_1,w_1,\psi_1)^T - (x_2,w_2,\psi_2)^T\|,$$

for some finite constant $C > 0$ uniform for $x, w$ and $\psi$. The first statement of Corollary 2.2.8 in van der Vaart and Wellner [6] and leads to

$$E \sup_{\theta_1,\theta_2\in\Theta_n\,:\,\|\theta_1-\theta_2\|\leq\delta} |Z_n^W(\theta_1) - Z_n^W(\theta_2)|$$

$$\leq E \sup_{\theta_1,\theta_2\in\Theta_n\,:\,\rho(\theta_1,\theta_2)\leq Ch^{-1}\delta} |Z_n^W(\theta_1) - Z_n^W(\theta_2)|$$

$$\leq C_4 \int_0^{Ch^{-1}\delta} \sqrt{\log(Cy^{-1}h^{-1})}\,dy,$$

where $C_4$ is some finite absolute constant uniform in $x, w$ and $\psi$. Proceeding similarly as in the second part yields the assertion.

*Fourth part: Order of the quantile*

With the second part of this proof obtain $E\big(\sup_{(x,w,\psi)\in\Theta_n} |Z_n^W(x,w,\psi)|\big) \leq C\sqrt{\log(n)}$ for some constant $C > 0$. By Sudakov's inequality (see Proposition A.2.5 in van der Vaart and Wellner [6]),

$$c\sqrt{\log(n)} \leq E\Big(\sup_{(x,w,\psi)\in\Theta_n} |Z_n^W(x,w,\psi)|\Big)$$

for some constant $c > 0$. This concludes the fourth statement of the lemma in view of Lemma 15. □

*C.5. Proofs of Lemmas 10 and 11*

We use $C > 0$ as a generic constant which can vary at every appearance.

*Proof of Lemma 10.* We only prove the first representation, as the second can be derived analogously. For sake of brevity let us write $S_n$ for $S_n(x;w,\psi)$ and

$$\check{F}(\mathbf{z};x,w,\psi,h) = \langle f\big(h^{-1}\mathcal{R}_{-\psi}(\mathbf{p}_{w,h}(x) - \mathbf{z})\big) g_2(\mathbf{z}), g_3(\mathbf{z})\rangle.$$

By Riemann-sum approximation

$$S_n = h^{-r_1}\int g_1(\mathbf{z})\check{F}(\mathbf{z};x,w,\psi,h)^j\,d\mathbf{z} + O\Big(\big(nh^{r_1-1-r_2-r_3}\big)^{-1}\Big),$$

27

uniformly for $x, w$ and $\psi$. To see this, write $I(n,h)$ for the index set for which the sum in $S_n$ is not zero and notice that $|I(n,h)| \leq 4n^2h^2$, due to the equidistant design and due to $\mathrm{supp}(f) = [-1,1]$. Let

$$A_{i_1,i_2} = [x_{i_1}, x_{i_1+1}) \times [x_{i_2}, x_{i_2+1}), \quad i_1, i_2 \in \{1, \ldots, n-1\},$$

where the right boundary is included if $i_1 = n-1$ or $i_2 = n-1$. Notice that $|x^2 - y^2| \leq 2C|x-y|$ for $x, y \in A$, where $A$ is a compact subset in $\mathbb{R}$ and $C = \sup A$. Thus, for any $\mathbf{y}_1, \mathbf{y}_2 \in [0,1]^2$,

$$|g_1(\mathbf{y}_1)\check{F}(\mathbf{y}_1; x, w, \psi, h)^2 - g_1(\mathbf{y}_2)\check{F}(\mathbf{y}_2; x, w, \psi, h)^2| \leq C|\sqrt{g_1(\mathbf{y}_1)}\check{F}(\mathbf{y}_1; x, w, \psi, h) - \sqrt{g_1(\mathbf{y}_2)}\check{F}(\mathbf{y}_2; x, w, \psi, h)|,$$

for some suitable constant $C > 0$. So it suffices to consider the case $j = 1$ by controlling the error term, say $Err$, in the Riemann-sum approximation. As the product of $g_1$ and $\check{F}$ is a $C^1$-function, the mean value theorem leads to

$$Err \leq \left(n^2 h^{r_1}\right)^{-1} \sum_{i_1, i_2 \in I(n,h)} |\sup_{\mathbf{y}_1 \in A_{i_1,i_2}} g_1(\mathbf{y}_1)\check{F}(\mathbf{y}_1; x, w, \psi, h) - \inf_{\mathbf{y}_2 \in A_{i_1,i_2}} g_1(\mathbf{y}_2)\check{F}(\mathbf{y}_2; x, w, \psi, h)|$$

$$\leq \frac{n^2 h^2}{n^3 h^{r_1}} \sup_{i_1, i_2 \in I(n,h)} \sup_{\mathbf{y} \in A_{i_1,i_2}} \|\nabla g_1(\mathbf{y})\check{F}(\mathbf{y}; x, w, \psi, h)\|,$$

since by the equidistant design $\sup_{i_1, i_2 \in I(n,h)} \sup_{\mathbf{y}_1, \mathbf{y}_2 \in A_{i_1,i_2}} \|\mathbf{y}_1 - \mathbf{y}_2\| \leq n^{-1}$. Considering the special representations of $g_2$ and $g_3$ we obtain by the chain rule

$$\sup_{i_1, i_2 \in I(n,h)} \sup_{\mathbf{y} \in A_{i_1,i_2}} \|\nabla g_1(\mathbf{y})\check{F}(\mathbf{y}; x, w, \psi, h)\| = O(h^{r_2+r_3-1}),$$

uniformly for $x, w$ and $\psi$. This uniformity comes from the fact that $x$ and $\psi$ are elements of compact sets and $w$ is such that $\mathbf{p}_{w,h}(x)$ is element of a compact subset. Therefore, $Err = O\left((nh^{r_1-1-r_2-r_3})^{-1}\right)$ uniformly over $x, w$ and $\psi$ as claimed. With the substitution

$$\mathbf{z} \mapsto h^{-1}\mathcal{R}_{-\psi}\left(\mathbf{p}_{w,h}(x) - \mathbf{z}\right) \tag{22}$$

one has

$$S_n = h^{2-r_1} \int_{[-1,1]^2} g_1\left(\mathbf{p}_{w,h}(x) - h\mathcal{R}_\psi \mathbf{z}\right) \check{F}(\mathbf{p}_{w,h}(x) - h\mathcal{R}_\psi \mathbf{z}; x, w, \psi, h)^j \, d\mathbf{z} + O\left((nh^{r_1-1-r_2-r_3})^{-1}\right)$$

$$= \int_{[-1,1]^2} g_1(\mathbf{p}_{w,h}(x) - h\mathcal{R}_\psi \mathbf{z}) \langle f(\mathbf{z})\tilde{g}_2(\mathbf{z}), \tilde{g}_3(\mathbf{z}) \rangle^j \, d\mathbf{z} + O\left((nh^{r_1-1-r_2-r_3})^{-1}\right),$$

uniformly for $x, w$ and $\psi$, where we used the choice of $r_1$ as well as (14) in the last line. Additionally, using the smoothness of $g_1$ obtain for $\mathbf{z} \in [-1,1]^2$ that

$$|g_1(\mathbf{p}_{w,h}(x) - h\mathcal{R}_\psi \mathbf{z}) - g_1(\mathbf{p}_{w,h}(x))| \leq Ch,$$

where $C > 0$ is a constant, which can be chosen uniformly for $x$, $w$ and $\psi$ as well. Hence, from the previous two displays one derives

$$S_n = g_1\left(\mathbf{p}_{w,h}(x)\right) \int_{[-1,1]^2} \langle f(\mathbf{z})\tilde{g}_2(\mathbf{z}), \tilde{g}_3(\mathbf{z}) \rangle^j \, d\mathbf{z} + O(h) + O\left((nh^{r_1-1-r_2-r_3})^{-1}\right),$$

uniformly for $x, w$ and $\psi$. □

*Proof of Lemma 11.* For sake of brevity write $J_n$ for $J_n(x; w, \psi)$ and

$$G(\mathbf{z}; x, w, \psi, h) := j_\tau(\mathbf{z}) \langle f(h^{-1}\mathcal{R}_{-\psi}(\mathbf{p}_{w,h}(x) - \mathbf{z}))g_2(\mathbf{z}), g_3(\mathbf{z}) \rangle.$$

Let $A_{i_1,i_2}$ and $I(n,h)$ be as in the proof of Lemma 10. Additionally, let $E(n,h)$ denote the set of indices in the sum

of $J_n$, for which the design points intersect with the curve $y = \phi(x)$ and $I^*(n,h) = I(n,h) \setminus E(n,h)$ be the set of the remaining indices for which the sum is not zero. Notice that $|E(n,h)| = O(nh)$, due to smoothness of $\phi$ and $|I^*(n,h)| \leq |I(n,h)| \leq 4n^2h^2$. By Riemann-sum approximation

$$J_n = h^{-r_1} \int G(\mathbf{z}; x, w, \psi, h) \, d\mathbf{z} + O\big((nh)^{-1}\big),$$

uniformly for $x, w$ and $\psi$. Indeed, the error term $Err$ in the Riemann-sum approximation can be bounded as follows

$$Err \leq \big(n^2 h^{r_1}\big)^{-1} \sum_{i_1, i_2 \in I^*(n,h)} |\sup_{\mathbf{y}_1 \in A_{i_1, i_2}} G(\mathbf{y}_1; x, w, \psi, h) - \inf_{\mathbf{y}_2 \in A_{i_1, i_2}} G(\mathbf{y}_2; x, w, \psi, h)|$$

$$+ \big(n^2 h^{r_1}\big)^{-1} \sum_{i_1, i_2 \in E(n,h)} |\sup_{\mathbf{y}_1 \in A_{i_1, i_2}} G(\mathbf{y}_1; x, w, \psi, h) - \inf_{\mathbf{y}_2 \in A_{i_1, i_2}} G(\mathbf{y}_2; x, w, \psi, h)|$$

$$=: Err_{(1)} + Err_{(2)}.$$

Note that $G(\cdot; x, w, \psi, h)$ is a $C^1$ function on the design squares $A_{i_1, i_2}$ for $i_1, i_2 \in I^*(n,h)$. Hence, we can proceed for $Err_{(1)}$ as in the proof of Lemma 10 to derive that $Err_{(1)} \leq C(nh)^{-1}$, due to choice of $r_1$, where the constant $C > 0$ can be chosen uniformly for $x$, $w$ and $\psi$. Additionally, by the special form of $g_2$ and $g_3$

$$\sup_{i_1, i_2 \in E(n,h)} \sup_{\mathbf{y} \in A_{i_1, i_2}} \|g_i(\mathbf{y})\| \leq Ch^{r_i}, \quad i = 2, 3,$$

where the constant $C > 0$ can be chosen uniformly for $x$, $w$ and $\psi$. With this it easily follows, $Err_{(2)} \leq C(nh)^{-1}$ such that $Err \leq C(nh)^{-1}$ for some constant $C > 0$ which is uniform in $x$, $w$ and $\psi$. Let

$$\mathbb{H}_{x, \psi, h} = h^{-1} \mathcal{R}_{-\psi} \big(\mathbf{p}(x) - [0,1]^2 \setminus \mathrm{epi}(\phi)\big).$$

With the substitution (22) and (14) obtain

$$J_n = h^{2-r_1} \int_{h^{-1} \mathcal{R}_{-\psi}(\mathbf{p}_{w,h}(x) - [0,1]^2 \setminus \mathrm{epi}(\phi))} \tau\big(x - (h\mathcal{R}_\psi \mathbf{z})_1\big) G(\mathbf{p}_{w,h}(x) - h\mathcal{R}_\psi \mathbf{z}; x, w, \psi, h) \, d\mathbf{z} + O\big((nh^{r_1 - 1 - r_2 - r_3})^{-1}\big)$$

$$= \int_{\mathbb{H}_{x,\psi,h} + w(\sin(\psi), \cos(\psi))^T} \tau\big(x - (h\mathcal{R}_\psi \mathbf{z})_1\big) \langle f(\mathbf{z}) \tilde{g}_2(\mathbf{z}), \tilde{g}_3(\mathbf{z}) \rangle \, d\mathbf{z} + O\big((nh^{r_1 - 1 - r_2 - r_3})^{-1}\big)$$

$$= \tau(x) \int_{\mathbb{H}(\psi(x) - \psi) + w(\sin(\psi), \cos(\psi))^T} \langle f(\mathbf{z}) \tilde{g}_2(\mathbf{z}), \tilde{g}_3(\mathbf{z}) \rangle \, d\mathbf{z} + R_n + O\big((nh^{r_1 - 1 - r_2 - r_3})^{-1}\big),$$

uniformly for $x, w$ and $\psi$, where

$$R_n = \int_{\mathbb{H}_{x,\psi,h} + w(\sin(\psi), \cos(\psi))^T} \big[\tau\big(x - (h\mathcal{R}_\psi \mathbf{z})_1\big) - \tau(x)\big] \langle f(\mathbf{z}) \tilde{g}_2(\mathbf{z}), \tilde{g}_3(\mathbf{z}) \rangle \, d\mathbf{z}$$

$$+ \int_{\mathbb{H}(\psi(x) - \psi) + w(\sin(\psi), \cos(\psi))^T \triangle \mathbb{H}_{x,\psi,h} + w(\sin(\psi), \cos(\psi))^T} \tau(x) \langle f(\mathbf{z}) \tilde{g}_2(\mathbf{z}), \tilde{g}_3(\mathbf{z}) \rangle \, d\mathbf{z}$$

$$=: R_{n,1} + R_{n,2}.$$

By smoothness assumptions on $\tau$ and the compact support of $f$ it holds $|R_{n,1}| \leq Ch$, where $C > 0$ can be chosen uniformly for $x$, $w$ and $\psi$. Moreover,

$$\lambda_2\big(\mathbb{H}(\psi(x) - \psi) \triangle \mathbb{H}_{x,\psi,h}\big) \leq Ch^2$$

and $C > 0$ is independent of $w$ and can be chosen uniformly as $x$ and $\psi$ take values in a compact subset of $\mathbb{R}^2$. Taking the compact support of $f$ into account, leads to

$$|R_{n,2}| \leq \sup_{\mathbf{z} \in [-1,1]^2} \|f(\mathbf{z}) \cdot \tilde{g}_2(\mathbf{z})\| \cdot \|\tilde{g}_3(\mathbf{z})\| \cdot \|\tau(\mathbf{z})\| \, \lambda_2\big(\mathbb{H}(\psi(x) - \psi) \triangle \mathbb{H}_{x,\psi,h}\big) = O(h^2),$$




\begin{thebibliography}{}

\bibitem[\protect\citeauthoryear{Bengs, Eulert, and Holzmann}{Bengs
  et~al.}{2018}]{beh2018supplement}
V. Bengs, M.~Eulert, H.~Holzmann,
Supplement to: Asymptotic confidence sets for the jump curve in
  bivariate regression problems,
Technical Report, 2019.

\bibitem[\protect\citeauthoryear{Bickel and Rosenblatt}{Bickel and
  Rosenblatt}{1973}]{bickel1973some}
P.J. Bickel, M.~Rosenblatt,
On some global measures of the deviations of density function
  estimates,
Ann. Statist. 1 (1973) 1071--1095.

\bibitem[\protect\citeauthoryear{Bissantz, D{\"u}mbgen, Holzmann, and
  Munk}{Bissantz et~al.}{2007}]{bissantz2007}
N. Bissantz, L.~D{\"u}mbgen, H.~Holzmann, A.~Munk,
Nonparametric confidence bands in deconvolution density estimation,
J. R. Stat. Soc. Ser. B (Stat. Methodol.) 69 (2007) 483--506.

\bibitem[\protect\citeauthoryear{Chernozhukov, Chetverikov, and
  Kato}{Chernozhukov et~al.}{2014}]{chernozhukov2014anti}
V. Chernozhukov, D.~Chetverikov, K.~Kato,
Anti-concentration and honest, adaptive confidence bands,
Ann. Statist. 42 (2014) 1787--1818.

\bibitem[\protect\citeauthoryear{Delaigle, Hall, and Jamshidi}{Delaigle
  et~al.}{2015}]{delaigle2015confidence}
A. Delaigle, P.G.~Hall, F.~Jamshidi, 
Confidence bands in non-parametric errors-in-variables regression,
J. R. Stat. Soc. Ser. B (Stat. Methodol.) 77 (2015) 149--169.

\bibitem[\protect\citeauthoryear{Eubank and Speckman}{Eubank and
  Speckman}{1993}]{eubank1993confidence}
R. Eubank, P.~Speckman,
\newblock Confidence bands in nonparametric regression,
J. Amer. Statist. Assoc. 88 (1993) 1287--1301.

\bibitem[\protect\citeauthoryear{Garlipp and M{\"u}ller}{Garlipp and
  M{\"u}ller}{2007}]{garlipp2007robust}
T. Garlipp, C.~M{\"u}ller,
Robust jump detection in regression surface,
Sankhy{\=a 69 (2007) 55--86.

\bibitem[\protect\citeauthoryear{Gijbels, Hall, and Kneip}{Gijbels
  et~al.}{2004}]{gijbels2004interval}
I. Gijbels, P.G.~Hall, A.~Kneip (2004).
Interval and band estimation for curves with jumps,
J. Appl. Probab. 41 (2004) 65--79.

\bibitem[\protect\citeauthoryear{{Gin{\'e}} and {Nickl}}{{Gin{\'e}} and
  {Nickl}}{2010}]{gine2010confidence}
E. {Gin{\'e}}, R.~{Nickl}, 
Confidence bands in density estimation,
Ann. Statist. 38 (2010) 1122--1170.

\bibitem[\protect\citeauthoryear{Kang and Qiu}{Kang and
  Qiu}{2014}]{kang2014jump}
Y. Kang, P.~Qiu, 
Jump detection in blurred regression surfaces,
Technometrics 56 (2014) 539--550.

\bibitem[\protect\citeauthoryear{Korostelev and Tsybakov}{Korostelev and
  Tsybakov}{1993}]{korostelev1993minimax}
A. Korostelev, A.~Tsybakov
Minimax Methods for Image Reconstruction,
Springer, New York, 1993.

\bibitem[\protect\citeauthoryear{Loader}{Loader}{1996}]{loader1996change}
C. Loader, 
Change point estimation using nonparametric regression,
Ann. Statist. 24 (1996) 1667--1678.

\bibitem[\protect\citeauthoryear{Mammen and Polonik}{Mammen and
  Polonik}{2013}]{mammenpolonik2013}
Mammen, E. and W.~Polonik (2013).
\newblock Confidence regions for level sets.
\em Journal of Multivariate Analysis\/}~{\em 122\/}(C), 202--214.

\bibitem[\protect\citeauthoryear{M{\"u}ller}{M{\"u}ller}{1992}]{muller1992change}
H.-G. M{\"u}ller,
Change-points in nonparametric regression analysis,
Ann. Statist. 20 (1992) 737--761.

\bibitem[\protect\citeauthoryear{M{\"u}ller and Song}{M{\"u}ller and
  Song}{1994}]{muller1994maximin}
H.-G. M{\"u}ller, K.~Song,
Maximin estimation of multidimensional boundaries,
J. Multivariate Anal. 50 (1994) , 265--281.

\bibitem[\protect\citeauthoryear{M{\"u}ller and Stadtm{\"u}ller}{M{\"u}ller and
  Stadtm{\"u}ller}{1999}]{muller1999discontinuous}
H.-G. M{\"u}ller, U.~Stadtm{\"u}ller, 
Discontinuous versus smooth regression,
Ann. Statist. 27 (1999) 299--337.

\bibitem[\protect\citeauthoryear{Munk, Bissantz, Wagner, and Freitag}{Munk
  et~al.}{2005}]{munk2005difference}
A. Munk, N.~Bissantz, T.~Wagner, G.~Freitag,
On difference-based variance estimation in nonparametric regression
  when the covariate is high dimensional,
J. R. Stat. Soc. Ser. B (Stat. Methodol.) 67 (2005) 19--41.
  
\bibitem[\protect\citeauthoryear{Neumann and Polzehl}{Neumann and
  Polzehl}{1998}]{neumann1998simultaneous}
M.H. Neumann, J.~Polzehl,
Simultaneous bootstrap confidence bands in nonparametric regression,
J. Nonparamet. Statist. 9 (1998) 307--333.

\bibitem[\protect\citeauthoryear{Porter and Yu}{Porter and
  Yu}{2015}]{porter2015regression}
J. Porter, P.~Yu,
Regression discontinuity designs with unknown discontinuity points:
  Testing and estimation,
  J. Econometrics 189 (2015)132--147.

\bibitem[\protect\citeauthoryear{Proksch, Bissantz, and Dette}{Proksch
  et~al.}{2015}]{proksch2015confidence}
K. Proksch, N.~Bissantz, H.~Dette, 
Confidence bands for multivariate and time dependent inverse
  regression models, Bernoulli 21 (2015) 144--175.

\bibitem[\protect\citeauthoryear{Qiao and Polonik}{Qiao and
  Polonik}{2016}]{qiaopolonik2016}
W. Qiao, W.~Polonik, 
Theoretical analysis of nonparametric filament estimation,
Ann. Statist. 44 (2016) 1269--1297.

\bibitem[\protect\citeauthoryear{Qiu}{Qiu}{1997}]{qiu1997nonparametric}
P. Qiu, 
 Nonparametric estimation of jump surface,
Sankhy{\=a}  59 (1997) 268--294.

\bibitem[\protect\citeauthoryear{Qiu}{Qiu}{2002}]{qiu2002nonparametric}
P. Qiu, 
A nonparametric procedure to detect jumps in regression surfaces,
J. Comput. Graph. Statist. 11 (2002) 799--822.

\bibitem[\protect\citeauthoryear{Qiu}{Qiu}{2005}]{qiu2005image}
P. Qiu,
Image Processing and Jump Regression Analysis,
Wiley, New York, 2005.

\bibitem[\protect\citeauthoryear{Qiu and Yandell}{Qiu and
  Yandell}{1997}]{qiu1997jump}
P. Qiu, B.~Yandell, 
Jump detection in regression surfaces,
J. Comput. Graph. Stat. 6 (1997) 332--354.

\bibitem[\protect\citeauthoryear{Seijo and Sen}{Seijo and
  Sen}{2011}]{seijo2011change}
E. Seijo, B.~Sen,
Change-point in stochastic design regression and the bootstrap,
Ann. Statist. 39 (2011) 1580--1607.

\bibitem[\protect\citeauthoryear{van~der Vaart}{van~der
  Vaart}{2000}]{van2000asymptotic}
A.W. van~der Vaart,
Asymptotic Statistics,
Cambridge University Press, Cambridge, 2000.

\bibitem[\protect\citeauthoryear{Wang}{Wang}{1995}]{wang1995jump}
Y. Wang, 
Jump and sharp cusp detection by wavelets,
Biometrika 82 (1995) 385--397.

\bibitem[\protect\citeauthoryear{Wang}{Wang}{1998}]{wang1998change}
Y. Wang, 
Change curve estimation via wavelets,
J. Amer. Statist. Assoc.   93 (1998) 163--172.

\bibitem[\protect\citeauthoryear{Wu and Chu}{Wu and Chu}{1993}]{wu1993kernel}
J. Wu, C.~Chu,
Kernel-type estimators of jump points and values of a regression
  function, Ann. Statist. 21 (1993) 1545--1566.

\end{thebibliography}
\end{document}